\documentclass{siamart0516}
\usepackage[margin=1in]{geometry}
\usepackage{amsmath,booktabs}
\usepackage{amssymb}
\usepackage{amsfonts} 
\usepackage[normalem]{ulem}
\usepackage[normal]{subfigure}

\providecommand{\U}[1]{\protect\rule{.1in}{.1in}}

\renewcommand\thesection{\arabic{section}}
\renewcommand\thesubsection{\thesection.\arabic{subsection}}
\renewcommand\thesubsubsection{\thesubsection.\arabic{subsubsection}}
\renewcommand\theequation{\thesection.\arabic{equation}}
\renewcommand{\eqref}[1]{(\ref{#1})}
\usepackage[numbers,sort&compress]{natbib}

\newcommand{\eps}{{\displaystyle \varepsilon}}
\newcommand{\bsub}{\begin{subequations}}
\newcommand{\esub}{\end{subequations}$\!$}
\newcommand{\ds}[0]{\displaystyle}
\newcommand{\bs}[0]{\boldsymbol}
\newcommand{\littleoh}{ \mbox{{\scriptsize $\mathcal{O}$}}}
\newcommand{\bigoh}{\mathcal{O}}

\newcommand{\bx}{\mathbf{x}}

\newcommand{\by}{\mathbf{y}}

\newcommand{\vs}{\bs{\sigma}}

\newcommand{\bxi}{\boldsymbol{\xi}}

\newcommand{\vx}{\mathbf{x}}

\newcommand{\wdisk}{{\overline{w}}}

\newcommand{\Om}{\Omega}
\newcommand{\pOm} {\partial \Omega}

\newcommand{\Gama}{\Gamma_a}
\newcommand{\Gamr}{\Gamma_r}

\newcommand{\al}[1]{\textcolor{black}{#1}}
\newcommand{\ajb}[1]{\textcolor{black}{#1}}

\newcommand{\A}{{\mathcal{A}}}
\newcommand{\B}{{\mathcal{B}}}
\newcommand{\Da}{\mathcal{D}_a}
\newcommand{\E}{{\mathcal{E}}}
\newcommand{\cH}{{\mathcal{H}}}
\newcommand{\J}{{\mathcal{J}}}
\newcommand{\cR}{R} 
\newcommand{\cS}{{\mathcal{S}}}
\newcommand{\cT}{{\mathcal{T}}}

\newcommand{\tilt}{\omega}

\newcommand{\hn} {{\hat{\textbf{n}}}}
\newcommand{\flux} {\rho}
\newcommand{\Flux} {P}

\newcommand{\tp}{2 \pi}
\newcommand{\rtp}{\frac{1}{2 \pi}}
\newcommand{\D}{D}  
\newcommand{\sone}{\sigma_1}
\newcommand{\stwo}{\sigma_2}
\newcommand{\surfrac}{\sigma}
\newcommand{\KMCProb}{p_M}
\newcommand{\sD}{\alpha}
\newcommand{\JS}{\mathcal{J}}
\newcommand{\fluxplane}{\mathcal{J}}
\newcommand{\Q}{Q}
\begin{document}

\title{Kinetic Monte Carlo methods for three-dimensional diffusive capture problems in exterior domains}

\date{\today} 
\author{
Andrew J. Bernoff\thanks{Dept.~of Mathematics, Harvey Mudd College, Claremont, CA, 91711, USA. {\tt bernoff@g.hmc.edu} } \and 
Alan E. Lindsay\thanks{ Dept.~of Applied and Computational Math \& Statistics, University of Notre Dame, Notre Dame, Indiana, 46656, USA. {\tt a.lindsay@nd.edu}}
}

\baselineskip=15pt

\maketitle

\begin{abstract}  
Cellular scale decision making is modulated by the dynamics of signalling molecules and their diffusive trajectories from a source to small absorbing sites on the cellular surface. Diffusive capture problems \ajb{which model this process}
are computationally challenging due to \ajb{their} complex geometry and  \ajb{mixed} boundary conditions together with intrinsically long transients that occur before a particle is captured. This paper reports on a particle-based Kinetic Monte Carlo (KMC) method that provides rapid accurate simulation of arrival statistics for (i) a half-space bounded by a surface with a finite collection of absorbing traps and (ii) the domain exterior to a convex cell again with absorbing traps. We validate our method by replicating classical results and \ajb{verifying some} newly developed boundary homogenization theories and matched asymptotic expansions on capture rates. In the case of non-spherical domains, we describe a new shielding effect in which geometry can play a role in sharpening cellular estimates on the directionality of diffusive sources.
\end{abstract}


\label{firstpage}



\pagestyle{myheadings}
\markboth{A.~J.~Bernoff, A.~E.~Lindsay.}{KMC for 3D Exterior Domains}

\section{Introduction}\label{sec:intro}

\begin{figure}[htbp]
\centering
\subfigure[Three dimensional diffusion to a reflecting plane with absorbing surface sites.]{\label{Fig:schematic_a} \includegraphics[width=0.55\textwidth]{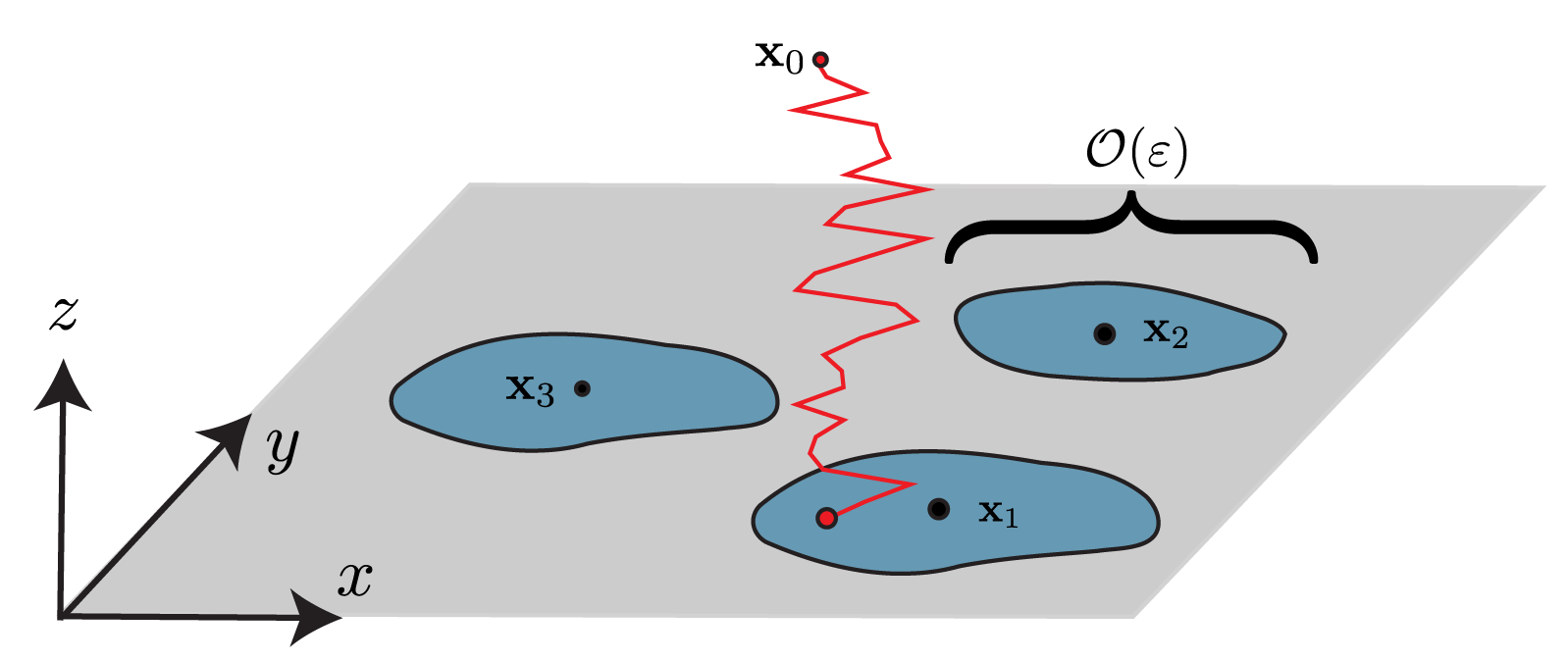}}\hspace{0.5in}
\subfigure[Diffusion to a triangulated convex surface with absorbing surface sites.]{\label{Fig:schematic_b} \includegraphics[width=0.315\textwidth]{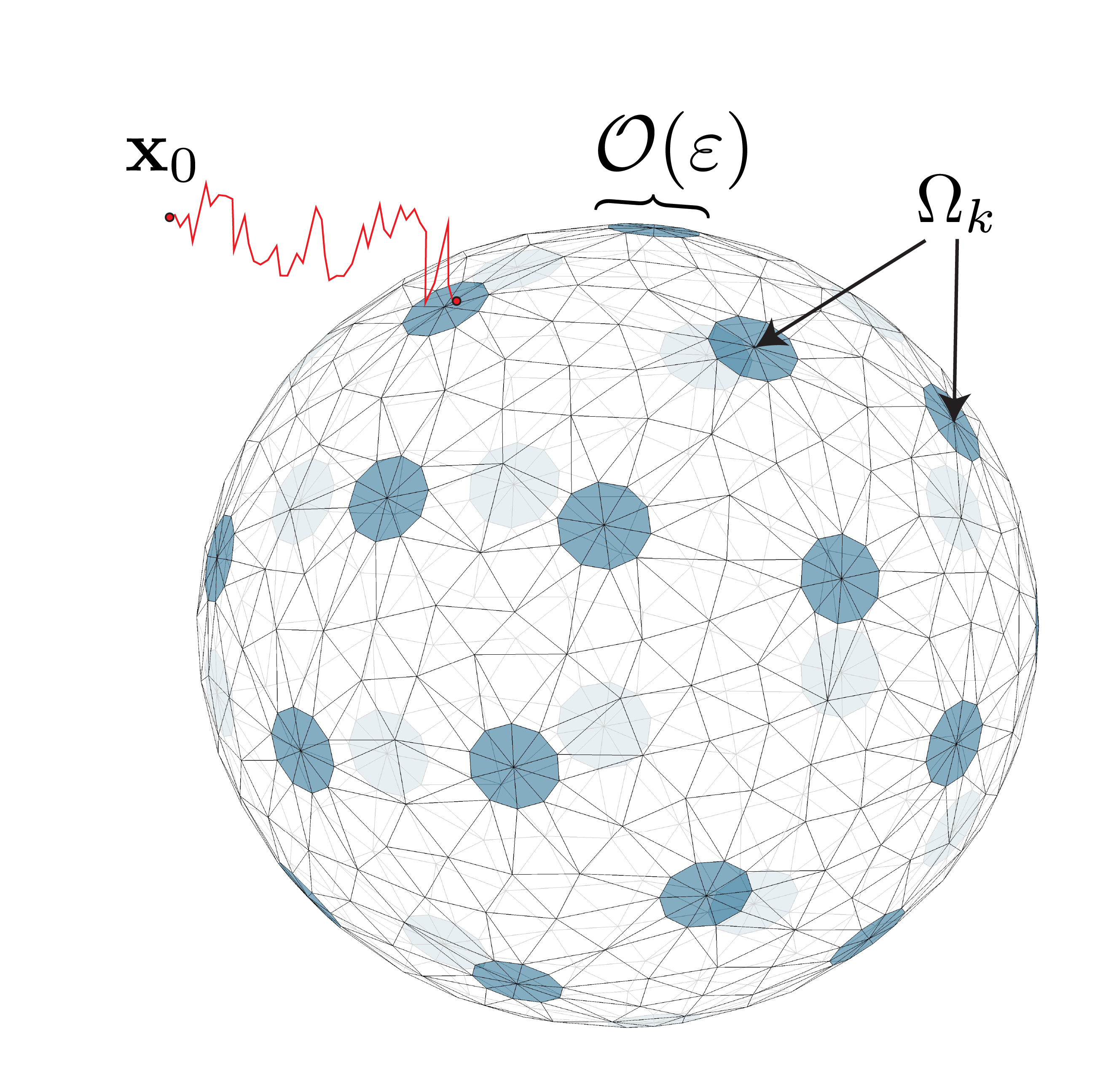}}
\caption{Schematic of three dimensional Brownian motion to a plane (left) and a convex surface (right) with absorbing surface sites. The absorbing sites are labeled $\Omega_k$ with centers $\bx_k$ and the initial particle location is $\bx_0$. We introduce and validate a numerical method to solve diffusion problems in both these scenarios. \ajb{The biologically relevant limit for this problem is when the typical pore size, denoted here by $\epsilon$, tends to zero.}
In the half-space geometry of panel (a), we derive expressions for the time-dependent fluxes to each of the absorbing sites in the limit as $\eps\to0$.  \ajb{For} convex geometries \ajb{illustrated}  in panel (b), we explore homogenized solutions to describe the capture rate.
 \label{Fig:schematic} }
\end{figure}

We consider the problem of computing the arrival time distributions of diffusing particles to absorbing sites arranged on planar and convex surfaces as shown in Fig.~\ref{Fig:schematic}. Related problems in the diffusive transport of cargo and chemical signaling are central in many biological phenomena and engineered systems \cite{ZSchuss2013,Wei2011,benichou2014first,schuss2012narrow,RH,app10186543,redner2001guide,FPPA2014,NewbyBressloff2013,Bressloff2024,Grebenkov2023}. For a particle released from location $\bx_0$, the central quantity of interest is the dynamic fluxes to each absorbing site together with the dependence on the number and spatial  configuration of these sites. The principal contribution of this work  is an efficient numerical method\footnote{Data and relevant code for this research work are stored in GitHub: \url{https://github.com/alanlindsay/3DKMC} and have been archived within the Zenodo repository: \url{https://doi.org/10.5281/zenodo.13997929}} to rapidly determine these quantities in the convex three dimensional geometries shown in Fig.~\ref{Fig:schematic}.

The general problem takes the form of a diffusion equation where for $\bx \in \Omega \subset \mathbb{R}^3$, the quantity $p(\bx,t;\bx_0)$ is the probability density that a particle originating at $\bx_0$ is free at time $t$ and position $\bx$. This distribution solves
\bsub\label{eqn:IntroP}
\begin{gather}
\label{eqn:IntroP_a} \frac{\partial p}{\partial t} = D \Delta p, \qquad \bx \in \Om, \quad t>0; \qquad 
p(\bx,0) = \delta (\bx-\bx_0), \quad \bx \in \Om;\\[5pt]
\label{eqn:IntroP_b} p = 0 \quad \mbox{on} \quad \bx \in \Gama; \qquad 
D\nabla p \cdot \hn = 0 \quad \mbox{on} \quad \bx \in \Gamr ,
\end{gather}
\esub
where $D$ is the diffusivity of the particle, the domain $\Omega$ is a subset 
of $\mathbb{R}^3$ whose boundary, $\pOm$ is partitioned into an absorbing set 
of pores, $\Gama$, whose impermeable complement $\Gamr$ is reflecting. The absorbing set is a union of $N$ non-overlapping pores $\Gama = \cup_{k=1}^N\Omega_k$.  The pores may be any shape, however on the plane we will usually consider circular pores and for domains exterior to a three-dimensional body we approximate smooth surfaces with convex polyhedra so the pores tend to be collections of polygonal surface facets (most often triangles).  We choose $\hn$, the normal to the surface $\pOm$, to point into the bulk. For a particle released from location $\bx_0$, the central quantity is the flux \ajb{density (sometimes known as the \emph{first passage time distribution})} 
through the $k^{th}$ absorbing site
$$\J_k(t;\bx_0)=D\int_{\Omega_k} \nabla p \cdot \hn\, dS.$$
Other important quantities of interest, such as the \ajb{flux density } to targets, total capture \ajb{probability} by $\Gama$ or splitting probabilities of individual receptors $\Omega_k$, can be obtained in terms of these fluxes.

We develop a particle based \emph{Kinetic Monte Carlo} (KMC) method which solves \eqref{eqn:IntroP} through sampling of the Brownian motion
\begin{equation}\label{introBM}
d \bx = \sqrt{2D}\, d\textbf{W}_t, \qquad \bx(0) = \bx_0.
\end{equation}
where $d\textbf{W}_t$ is the increment of a Wiener process \cite{Risken1992,redner2001guide}. Rapid and accurate solution of \eqref{introBM} is accomplished by combining three exact solutions of the heat equation for carefully chosen geometries. First, projection from a point to a plane. Second, projection from a point exterior to a sphere to its surface. Third, projection from the center of a hemisphere to its surface. In Sec.~\ref{sec:KMCalg}, we describe how appropriate combinations of these projection operators allow for exact simulation of \eqref{introBM} and hence sampling of \eqref{eqn:IntroP} for convex bodies with combinations of absorbing and reflecting portions. 

The KMC method developed here belongs to a class of meshless methods \cite{MAURO2014,GROSS2022,HWANG20101089,SCHUMM2023,PlunkettLawley2024} that can bypass obstacles inherent to traditional solution methods (e.g.~finite element, finite difference, and boundary integral methods), when resolving singularities of the surface flux at the interface of Neumann and Dirichlet components. This benefit of the KMC method is tempered by relatively slow convergence of accuracy; for a given statistic (such as the flux into a given pore over a unit of time) if $M$ particles are captured, the error scales as $M^{-\frac12}$ as $M\to\infty$. Fortunately particles can be simulated independently facilitating parallel computation for very large numbers of particles (millions, billions, or more) as required for the desired accuracy.
 
Mixed boundary values problems such as \eqref{eqn:IntroP} are surprisingly resilient to enquiry from classical solution methods and closed form solutions have been developed only in the case of steady state solutions with one \cite{sneddon1966mixed} or two \cite{Strieder08,Strieder12} circular pores. Boundary homogenization \cite{berez2006,muratov,Grebenkov2023} seeks to remediate these limitations by replacing the configuration $\Gamma_a$ with a single Robin condition $D\nabla p \cdot \hn = \kappa p$ over $\partial\Omega$. The problem is then reduced to determining the single parameter $\kappa$ which best represents the general configuration, usually by replicating certain global quantities, such as the total capture rate or capacitance. This process has been successfully applied to approximate diffusive transport over periodic \cite{muratov,belyaev,LBS2018,Bruna15,PlunkettLawley2024} and finite arrangements \cite{berez2014,Berez2012,Manhart2022} of surface sites, however, it does not easily yield the distribution across individual sites $\Omega_k$, which are needed in many biological applications such as inferring the direction of chemical cues \cite{LLM2020,Lindsay2023a,Lindsay2023b} through splitting probabilities or understanding how receptor clustering modulates immune signalling \cite{Care2011}.

While the contribution of this work is principally the development of KMC methods, in the case of well separated pores on the plane (cf.~Fig.~\ref{Fig:schematic_a}), we obtain a new result in the form of a two-term asymptotic expansion for the flux to each pore. This limit is most relevant for biological applications where reactive receptors occupy a small portion of the cellular surface. If we consider the absorbing set $\Gama$ to be $N$ non-overlapping absorbing pores in the plane $\pOm$ with centers given by coordinates $\bx_k= (x_k,y_k,0)$, we have that
\begin{equation}\label{eq:domains_intro}
\Gamma_a = \bigcup_{k=1}^N  \A_k \left (\frac{\bx -\bx_k} {\eps} \right ), \qquad \bx = (x,y,0).
\end{equation}
Here $\A_k(\bx)$ is  a closed bounded set representing the shape of the $k^\text{th}$ pore \ajb{ if it were centered at the origin. The parametrized family of homothetic pores, $\A_k \left (\frac{\bx -\bx_k} {\eps} \right )$, are of the same shape as $\A_k(\bx)$ but scaled by $\eps$
and centered at $\bx_k$}. In section \ref{sec:asyplane} we derive that as $\eps\to0$, \ajb{the flux density into the $k^{th}$ pore is asymptotically}
\begin{equation}\label{eq:MainPlanarResult}
\J_k(t;\bx_0)  =  \frac{\eps c_k}{\sqrt{4\pi D t^3}} e^{-\frac{R_k^2}{4Dt}} \Big[
  1 - \eps \left ( \frac { c_k}{ R_k} - \frac{ c_k R_k} {2D t}  \right ) \\
+  \eps \sum_{\substack{j =1\\ j\neq k} }^N c_je^{\frac{R_k^2-(R_j+d_{jk})^2}{4Dt} }
 \left  (\frac 1 { R_j}+ \frac 1 {d_{jk}} \right ) \Big] +\bigoh(\eps^3),
\end{equation}
where $d_{jk} = |\bx_j-\bx_k|$, $R_k = |\bx_0 - \bx_k|$ and $c_k$ is the capacitance of the \ajb{unscaled} $k^{th}$ absorber\footnote{\ajb{Capacitance scales as dimension; if the capacitance of $\A_k(\bx)$ is $c_k$ then the capacitance of $\A_k\left (\frac{\bx -\bx_k} {\eps} \right )$ is $\eps c_k$.}}. In the case of a circular absorber of radius $\eps a_k$, it is known that $\eps c_k = 2\eps a_k/\pi$. For \ajb{absorbers} of general shape $\A_k$, our numerical algorithm can calculate $c_k$ to high precision. The expression \eqref{eq:MainPlanarResult} can be used to derive many important quantities relating to capture statistics, including the \ajb{flux density} $\J(t;\bx_0) = \sum_{k=1}^{N} \J_k(t;\bx_0)$ or the steady state flux to receptors \ajb{from a continuous source}, also known as the \emph{splitting probabilities},  $\Q_k(\bx_0) = \int_{0}^{\infty} \J_k(\tilt;\bx_0) \, d\tilt$. We remark that \eqref{eq:MainPlanarResult} is one of only a few \ajb{known} closed form expressions for the full distribution of arrival times in a narrow escape problem \cite{IN2013,LTS2016}.

\begin{table}

\al{
    \begin{tabular}{ccc}
    Symbol & Description & Probabilistic interpretation\\
    \midrule
    $\J_k(t)$ & $\begin{array}{c} \text{Flux density to the} \\ k^{th} \text{ receptor at time } t.\end{array}$ & $\begin{array}{c} \mathcal{J}_k(t)\,dt \text{ is the probability a particle} \\ \text { reaches } k^{th} \text{ receptor in } (t,t+dt).\end{array}$\\[10pt]
    $\J(t) = \sum_{k=1}^N \J_k(t)$ & $\begin{array}{c} \text{Flux density to the }\\ \text{absorbing set at time } t. \end{array}$ & $\begin{array}{c} \J(t)\,dt \text{ is the probability a particle}\\ \text{  reaches the absorbing set in } (t,t+dt).\end{array} $ \\[10pt]
    $q_k (t)  = \int_0^{t}\J_k(\tau)d\tau$ & $\begin{array}{c} \text{Cumulative flux to}\\
    \text{the } k^{th} \text{receptor.}\end{array}$ & $\begin{array}{c} \text{Probability a particle is absorbed by}  \\ \text{the } k^{th} \text{ receptor by time } t.\end{array} $\\[10pt]
    $F(t) = \int_0^t \J(\tau)d \tau$ & $\begin{array}{c} \text{Cumulative flux to the}\\ \text{absorbing set.} \end{array}$ & $\begin{array}{c} \text{Probability that a particle} \\ \text{ reaches the absorbing set by time }t.\end{array}$\\[10pt]
     $\Q_k(\bx_0) = \lim_{t\to \infty} q_k(t)$ & $\begin{array}{c} \text{Splitting probability to} \\  \text{the } k^{th} \text{ receptor.}\end{array}$ & $\begin{array}{c} \text{ Probability a particle released at a point } \\ \text {$\bx_0$ is absorbed by the }  k^{th} \text{ receptor.}\end{array}$\\
    \bottomrule
    \end{tabular}
    }
    \vspace{10pt}
    \caption{\al{The quantities calculated by both analytical approximations of the governing equation \eqref{eqn:IntroP} and discrete simulation by the KMC method. Our methods exploit the established connections \cite{Risken1992,redner2001guide} between the Langevin equation \eqref{introBM} and the associated PDE \eqref{eqn:IntroP}. The quantities listed all depend upon the initial release point, $\bx_0$, however, this is only referred to explicitly for the splitting probability $\Q_k(\bx_0)$.}}
\end{table}

The scenario where absorbing sites are arranged on a sphere, or a more general surface, is an important generalization of the planar problem which allows for consideration of geometric effects in capture statistics. For a sphere with circular pores, boundary integral methods have been developed to compute the steady solution to high precision \cite{LBS2018,KAYE2020}. Homogenization for the sphere was performed \cite{LWB2017} on the steady state in the scenario \al{of $N$ pores with common radius $a$ in the limit} $N\to\infty$, $a\to0$ with fixed absorption fraction $\surfrac = \frac{Na^2}{4}$. This gives the effective boundary condition
\begin{equation}\label{eq:IntroRobinBC}
D\nabla p \cdot \hn = \kappa p, \quad |\bx| = 1, \qquad \kappa = \frac{4D\surfrac}{\pi a}\left[1- \frac{4}{\pi} \sqrt{\surfrac} + \frac{a}{\pi}\log(4e^{-1/2}\sqrt{\surfrac})\right]^{-1}.
\end{equation}
The radially symmetric problem which results from coupling \eqref{eqn:IntroP_a} with \eqref{eq:IntroRobinBC} is solvable and yields the \ajb{flux density}
\begin{equation}\label{eq:IntroHomogSphereFlux}
\J(t) = \frac{\kappa}{R} 
e^{-\frac{\left(R -1\right)^2}{4 \D t}}
\left[ \frac{1}{\sqrt{\pi  \D  t}}
-\mathrm{erfc}\! \left(\beta \right)
 {\mathrm e}^{\beta^2}  \left(\kappa/D +1\right)\right],
\end{equation}
where $R=|\bx_0|$ and $\beta= \frac{R -1}{2 \sqrt{Dt}}+\left(\kappa/D +1\right) \sqrt{Dt}$. The new computations presented in this work verify that beyond predicting the steady state capture statistics, the homogenized flux \eqref{eq:IntroHomogSphereFlux} is a remarkably accurate predictor of absorption to $\partial\Omega$ over almost all timescales, with the exception of the exponentially small initial regime as $t\to0^{+}$. In Sec.~\ref{sec:resultsCube} we investigate a similar \ajb{homogenization} process for the well known computational challenge problem of determining the capacitance of the cube. We show that replacing the cube with a simpler spherical geometry of equivalent capacitance again yields a remarkably accurate capture rate at all timescales, except as $t\to0^{+}$.

In a final example we consider a family of oblate ellipsoids with two \ajb{flat} circular pores located at the north and south poles. The equator radius varies; in one limit the surface is a circular cylinder while in the limit of large oblateness the skirt between the two pores becomes a wide disc-like barrier. \al{Ratiometric sensing is a proposed biological mechanism for cells to fix the direction of signaling sources through the differential fluxes between surface pores \cite{Lew2019,Ismael2016,Lakhani2017,Bumsoo2019}. In this geometry we demonstrate that a wide skirt can block the diffusion of particles and enhance the ratiometric sensing mechanism for directional inference.}

\section{Kinetic Monte Carlo Methods}\label{sec:KMC}

Monte Carlo simulations provide a valuable tool for numerically estimating the distribution of capture times of diffusing particles for problems such as \eqref{eqn:IntroP} and have been used extensively \cite{BatBerShv2003,berez2004,berez2006,berez2014,Northrup1988,Litwin1980,Opplestrup2006,Oppelstrup2009,CLNQ22}. 
In its simplest form, a Monte Carlo method simulates the diffusive (Brownian) motion of a particle as a sequence of small displacements of randomly chosen orientation which terminates when the particle transits an absorbing surface. The algorithm is repeated for many particles (millions or even billions) to sample the capture time distribution. These Monte Carlo methods are hampered by a set of problems. First, the adoption of a fixed stepsize introduces an error at that lengthscale. Second, in capture problems such as \eqref{eqn:IntroP} with fat-tailed distributions a significant fraction of realizations undergo long excursions before they are captured, particularly when the domain is unbounded and/or the pores are small. Third, for exterior domains a finite proportion of particles are never captured by a pore and one must decide when/if a particle has escaped capture.

Kinetic Monte Carlo (KMC) methods \cite{BatBerShv2003} split the diffusion process into steps, where each step corresponds to a diffusion problem on a simpler geometry that can be solved analytically. For example, in the case of a half-space above a planar boundary, the distribution of the time and location of the first impact on the plane can be solved analytically and numerically sampled replacing the simulation of long excursions of the particles with a single calculation. Early work used these ideas in $N$-body simulations of kinetic gases \cite{Opplestrup2006} and chemical reactions \cite{Wu2006}.  

In Section \ref{subsec:KMCProp} we describe three KMC propagators used to simulate portions of the diffusion process in various geometries. Section \ref{sec:KMCalg} describe how these propagators are assembled to simulate diffusion in a half-space bounded by a plane that is reflecting except for a finite collection of absorbing pores (Sec.~\ref{subsec:KMCplane}) and to a convex polyhedron with whose faces are either reflecting or absorbing (Sec.~\ref{subsec:KMCconvexpoly}). In practice, we demonstrate this method on triangulations representing convex bodies such as spheres and ellipsoids.

\subsection{KMC Propagators}\label{subsec:KMCProp}

Our KMC method  constructs arbitrary Brownian paths \eqref{introBM} from an initial location to planar surfaces (Fig.~\ref{fig:KMCplane}) and convex polyhedra (Fig.~\ref{fig:KMCpoly})
via three exactly solvable diffusion problems (propagators) that are described below.

\subsubsection{KMC Projector I: Propagating from a point in a half-space to the bounding plane}
\label{subsubsec:KMCHalf}
\begin{figure}[tbp]
\centering
\includegraphics[width = 0.95\textwidth]{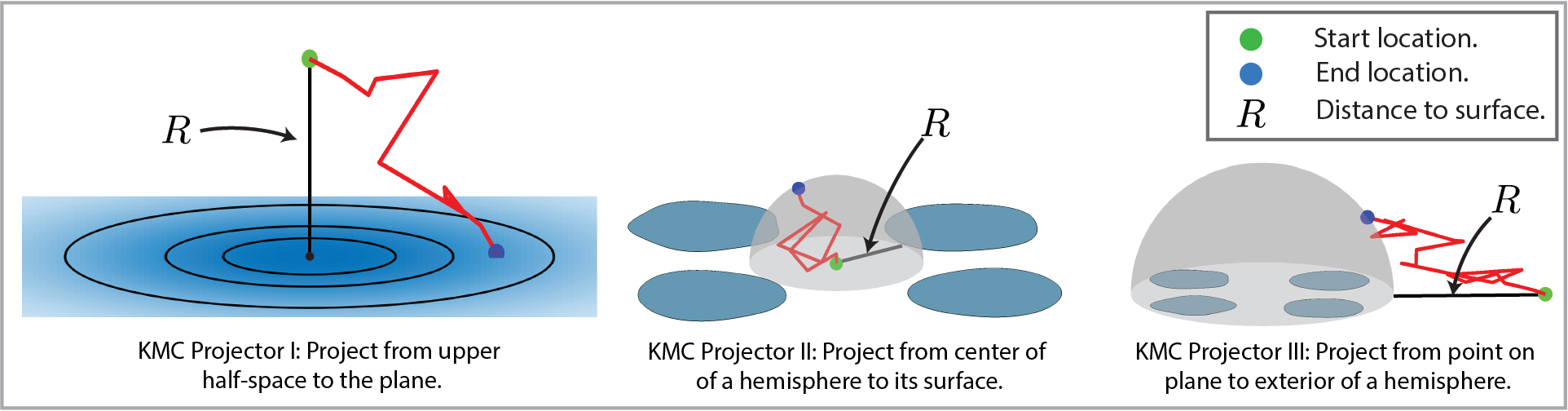}
\caption{Schematic of the three KMC projectors in the simulation of diffusion from the upper half-space to a plane with a finite set of absorbers.  KMC Projector I described in Sec.~\ref{subsubsec:KMCHalf} projects to the plane. KMC Projector II described in Sec.~\ref{subsubsec:KMCSphere} projects from the plane to the surface of the hemisphere of radius $R$ determined by the distance to \ajb{the} closest absorber. KMC Projector III described in Sec.~\ref{subsubsec:KMCReinsert} projects from the plane to a hemisphere that encloses all absorbers.\label{fig:KMCplane}}
\end{figure}

\begin{figure}[tbp]
    \centering
    \includegraphics[width=0.9\linewidth]{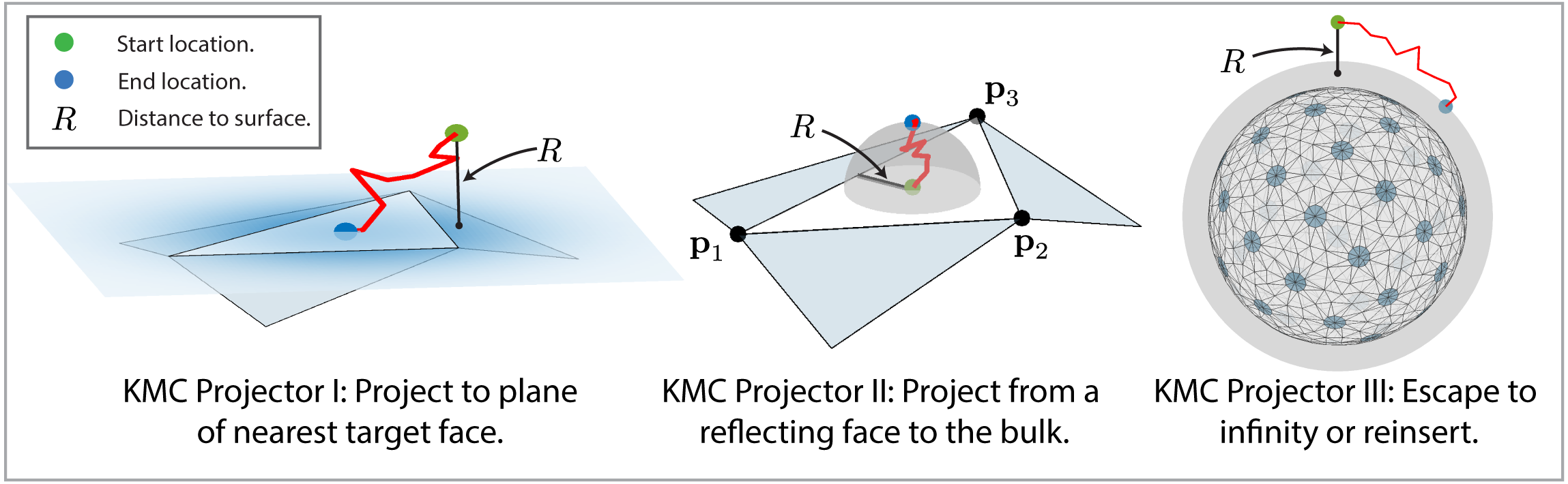}
    \caption{The three projection stages for the simulation of diffusing particles to convex surfaces with absorbing and reflecting patches. KMC Projector I described in Sec.~\ref{subsubsec:KMCHalf} projects to the plane with the largest signed distance $R$. KMC Projector II described in Sec.~\ref{subsubsec:KMCSphere} \ajb{projects a particle on a face to the largest hemisphere whose base is a circle of radius $R$ contained within the face and centered at the particle.} KMC Projector III described in Sec.~\ref{subsubsec:KMCReinsert} projects to a sphere that encloses the target. }
    \label{fig:KMCpoly}
\end{figure}

For particles diffusing toward a set of absorbers on a plane we can propagate a particle forward from the bulk to the bounding plane (see Fig.~\ref{fig:KMCplane}) using an exact solution derived via the method of images.
Similarly, for particles in the bulk exterior to a convex polyhedron one can identify a half-space contained within the bulk whose bounding plane contains one of the polyhedron's faces (see Fig.~\ref{fig:KMCpoly}). Here we describe how to propagate the particle to this bounding plane which it will impact almost\footnote{The term \emph{almost} here is used in the sense of measure/probability theory - the particle impacts the plane with probability one, although there is a set of zero probability events which would allow the particle to escape to infinity.} certainly.

Consider a coordinate system with the origin at the point closest to the particle on the bounding plane for which the half-space $\Omega$ corresponds to $z>0$. 
The arrival time and position distribution for the particle's first impact on this plane is computed from a density $p(x,y,z,t)$ which satisfies the diffusion equation (\ref{eqn:IntroP}) in the bulk with a delta function source at $\vx_0=(0, 0, z_0)$ and an absorbing boundary at $z=0$. In the framework of (\ref{eqn:IntroP}), $\Gamma_a$ is the entire bounding plane on which $p(x,y,0,t)=0$.

The solution is constructed via the method of images and the known free space Green's function \cite{CJ}, 
\begin{equation} 
\label{PlaneFluxDensity}
p(x,y,z,t) = \frac{1}{(4\pi D t)^{3/2}} \left[ e^{ -\frac{x^2+y^2+ (z-z_0)^2}{4 D t} } - e^{  -\frac{x^2+y^2+ (z+z_0)^2}{4 D t} } \right].
\end{equation}
The flux density through the boundary, $\fluxplane(t)$, is the probability distribution function (PDF) of transit times to the plane,
\begin{equation}
\label{PlaneFlux}
\fluxplane(t) =\iint_{z=0} D \nabla p \cdot \hat{\bf n} \, dA =
\int_{x=-\infty}^\infty \int_{y=-\infty}^\infty D p_z(x,y,0,t)  ~ dy \, dx= \frac{z_0}{2 \sqrt{\pi D}\, t^{3/2}} e^{-\frac{z_0^2}{4D t}}\, ,
 \end{equation}
with an associated cumulative distribution function (CDF)  given by
\begin{equation}\label{eq:DistPlane}
 \Flux_T(t) = \int_0^t \fluxplane (\tilt) ~d\tilt =\int_0^t \frac{z_0}{2\sqrt{\pi D }\, \tilt^{3/2}} e^{-\frac{z_0^2}{4D \tilt}} ~d\tilt =  \text{erfc} \left( \frac{z_0}{2\sqrt{Dt}} \right).
 \end{equation}
This distribution of transit times, $t_*$, from the bulk to the plane can then be sampled by choosing a uniform random variable\footnote{\ajb{Uniform random variables in {\sc Matlab} exclude the interval endpoints, which is reflected when we say $\nu \in (0,1)$. This conveniently avoids certain potential infinities such as when $\nu=0$ or $\nu=1$ in \eqref{PlaneImpact}. }}, $\nu$, from the interval $(0,1)$ and letting\footnote{\al{In our implementation, we use the {\sc Matlab} function {\tt erfinv} to generate samples of $\text{erfc}^{-1}(\nu)$.}}
\begin{equation}
\label{PlaneImpact}
 t_* =\frac 1 {4D} \left [\frac {z_0}  {\text{erfc}^{-1}(\nu) }\right ] ^2, \qquad \nu \in (0,1).
\end{equation}
The spatial distribution of the particle flux \ajb{density} at the arrival time $ t_*$, $\flux_{XY}(x_*,y_*;t_*)$, is given by
\begin{align}
\label{PlaneImpactPoint} 
\flux_{XY}(x_*,y_*;t_*) &\equiv  
\frac{1} {\fluxplane(t_*)}
\left [
D \nabla p \cdot \hat{\bf n}
\right ]_{t=t_*} 
=
\frac{1} {\fluxplane(t_*)}
{D p_z(x_*,y_*,0,t_*)}
\\
&=\frac{1}{4\pi D t_*}e^{ -\frac{x_*^2+y_*^2}{4Dt_*}}= \left [\frac{1}{\sqrt{4\pi D t_*}}e^{ -\frac{x_*^2}{4D t_*}} \right ] \left [\frac{1}{\sqrt{4\pi Dt_*}}e^{ -\frac{y_*^2}{4Dt_*}} \right ] =\flux_X(x_*;t_*) \flux_Y(y_*;t_*), \nonumber
\end{align}
which is the product of two Gaussian (normal) distributions. Accordingly, $x_*$ and $y_*$ are both drawn from $\mathcal{N}(0,2Dt_*)$, the normal distribution with mean zero and variance $2Dt_*$.

To summarize, one first determines the time that elapses to impact, $t_*$, via \eqref{PlaneImpact}, followed by the horizontal displacements, $x_*$ and $y_*$, via \eqref{PlaneImpactPoint}. The particle is then propagated forward in time and displaced to the bounding plane accordingly.

\subsubsection{KMC Projector II: Propagating from the center of a hemisphere to its surface}
\label{subsubsec:KMCSphere}

If a particle has impacted a reflecting portion of the surface, the next step is to propagate it back into the bulk. We can propagate it forward to the surface of a hemisphere, $\cH$, of radius $R$ centered at the impact point 
(See Fig.~\ref{fig:KMCplane} and Fig.~\ref{fig:KMCpoly} for the half-space and convex polyhedral versions respectively). The value of $R$ is chosen as the radius of largest circle that can be inscribed on the surface that remains within the reflecting portion \ajb{of the plane}.\footnote{In practice one needs to bound the value of $R$ from below (we use ten times double precision machine \ajb{epsilon}, $\epsilon_M$) - Monte Carlo methods are notorious for exploring all edge cases and a particle can get stuck on the edge of an absorber if it is close enough that $R$ rounds to zero.}
This corresponds to solving the diffusion equation (\ref{eqn:IntroP}) in a hemispherical domain 
$$\Omega = \left \{\
(x,y,z)\ |\ x^2+y^2+z^2 < R^2, \ z>0 
\right \},$$
with a reflecting circular base and absorbing surface ${\cH}$,
$$\Gamma_a = {\cH} 
\equiv
\left \{\
(x,y,z)\ |\ x^2+y^2+z^2 =R^2, \ z>0 
\right \} 
\qquad 
\Gamma_r = 
\left \{\ (
x,y,0)\ |\ x^2+y^2 < R^2
\right \}.
$$
The particle is initially infinitesimally above the origin, that is
$\vx_0=(0, 0,z_0)$ for $0<z_0\ll1$.

The joint distribution for the first exit time and exit location on the hemisphere can be deduced by noting the equivalent problem on the full sphere has a solution that is radial,  $p(x,y,z,t) =p(r,t)$ for $r = \sqrt{x^2+y^2+z^2}$,  thus insuring that the reflection symmetry on the surface is satisfied.  The separable solution (cf.~\cite{Litwin1980,LBS2018}) yields the CDF of exit times,
 \begin{align}
 \label{eq:CDFFull}
 \Flux_\cT(\tau) & 
 = 1+ 2 \sum_{n=1}^{\infty}(-1)^n e^{-n^2 \tau }, \qquad \tau = \frac{\pi^2D t}{R^2}.
\end{align}
The series \eqref{eq:CDFFull} converges quickly for large $\tau$, but more slowly when $\tau$ is small. The following theta function identity (derived via the Poisson summation formula)  
remedies this issue (cf.~\cite[Ch.~4]{Stein2003}),
\begin{equation}\label{eq:hemi_cdf1}
\sum_{n=-\infty}^\infty e^{-\pi q (n+a)^2 } = \sum_{n=-\infty}^\infty q^{-1/2} e^{-\pi n^2 / q} e^{2\pi i n a} .
\end{equation}
Applying the identity \eqref{eq:hemi_cdf1}, with $a = \frac{1}{2}$ and $q=\pi/ \tau$, to \eqref{eq:CDFFull} yields that
\begin{equation}\label{eq:hemi_cdf_extra}
\Flux_\cT(\tau) =  2\sqrt{\frac{\pi}{\tau}} \sum_{n=0}^\infty e^{-\pi^2 (n+\frac{1}{2})^2 /\tau},
\end{equation}
which converges rapidly for $\tau$ small. To sample an arrival time to the sphere, we draw a uniform random number $\xi\in (0,1)$ and numerically solve the equation
\begin{equation}
\label{HemiImpact}
\Flux_\cT(\tau_*) = \xi,
\end{equation}
for $\tau_*$. After rescaling, this yields the exit time $t_*=R^2 \tau_*/D\pi^2$. The values of $\Flux_\cT(\tau)$ are precomputed and tabulated for computational efficiency (using \eqref{eq:CDFFull} for $\tau \ge 1$  and \eqref{eq:hemi_cdf_extra} for $\tau < 1$)  and the value of $\tau_*$ is determined by linear interpolation\footnote{In general we choose enough points to bound the relative error here to one part in $10^6$.} unless $\xi$ is close to unity in which case the asymptotic approximation $\tau \sim  \ln [2/(1-\xi)] + \bigoh \left ((1-\xi)^3 \right )$ is used.

Once an exit time has been determined, an exit point on the hemisphere, $\cH$, can be chosen isotropically (due to the spherical symmetry).
As the surface area element satisfies 
$$dS = R^2 \sin \phi ~d\theta~ d \phi = R ~d \theta ~dz, $$
one can select a pair of random variables, $\eta \in (0, 2 \pi)$ and $\zeta \in (0,1)$ and an associated exit point 
\begin{equation}\label{eq:HemiSpherePoint}
(x_*,y_*,z_*) = R\, \Big( \sqrt{1-\zeta^2}\, \cos  \eta,  \sqrt{1-\zeta^2} \, \sin  \eta, \zeta \Big)
\qquad \eta \in (0, 2 \pi), \quad \zeta \in (0,1)
.
\end{equation}
To summarize, one first determines the time that elapses to impact, $t_*$, via \eqref{HemiImpact}.
The exit point $(x_*,y_*,z_*)$ on the hemisphere, $\cH$, is then chosen isotropically via \eqref{eq:HemiSpherePoint}. The particle is then propagated forward in time and displaced to the bounding hemisphere accordingly.

\subsubsection{KMC Projector III: Propagating from a point exterior to a sphere to its surface}
\label{subsubsec:KMCReinsert}
In three dimensional exterior domain problems, we must account for the finite probability that any particle can escape to infinity. We account for this with a KMC projector such that particles wandering sufficiently far from the absorbers are either propagated \ajb{inward} to a sphere (or hemisphere) bounding the target or marked as having escaped to infinity.

\begin{figure}[htbp]
\centering
\subfigure[Schematic of reinsertion.]{\includegraphics[width = 0.38\textwidth]{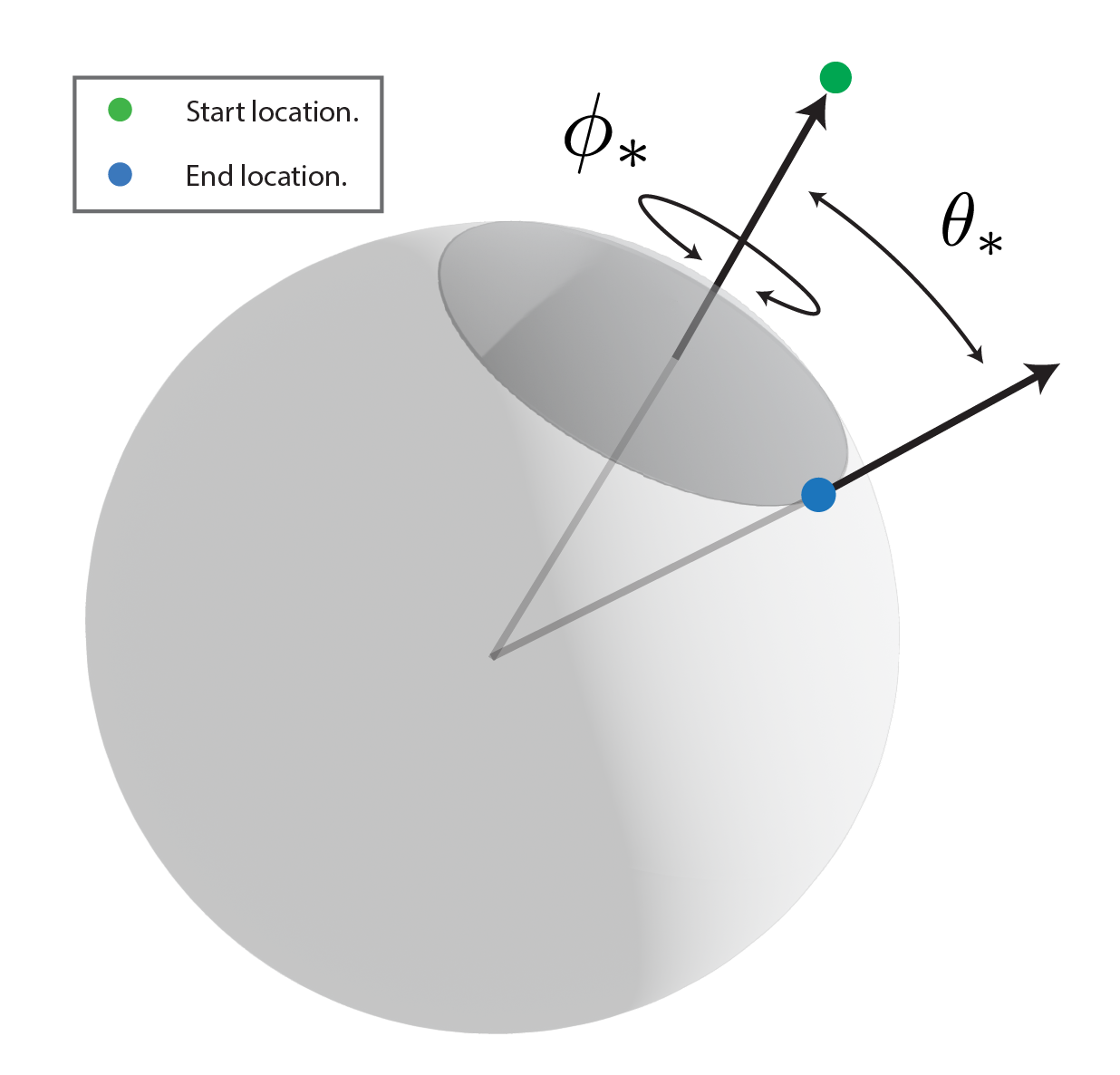}\label{fig:reinsertion_a}}\qquad
\subfigure[Curves of $P_{\Theta}( \theta; t_{\ast})$ for several $t_{\ast}$.]{\includegraphics[width = 0.425\textwidth]{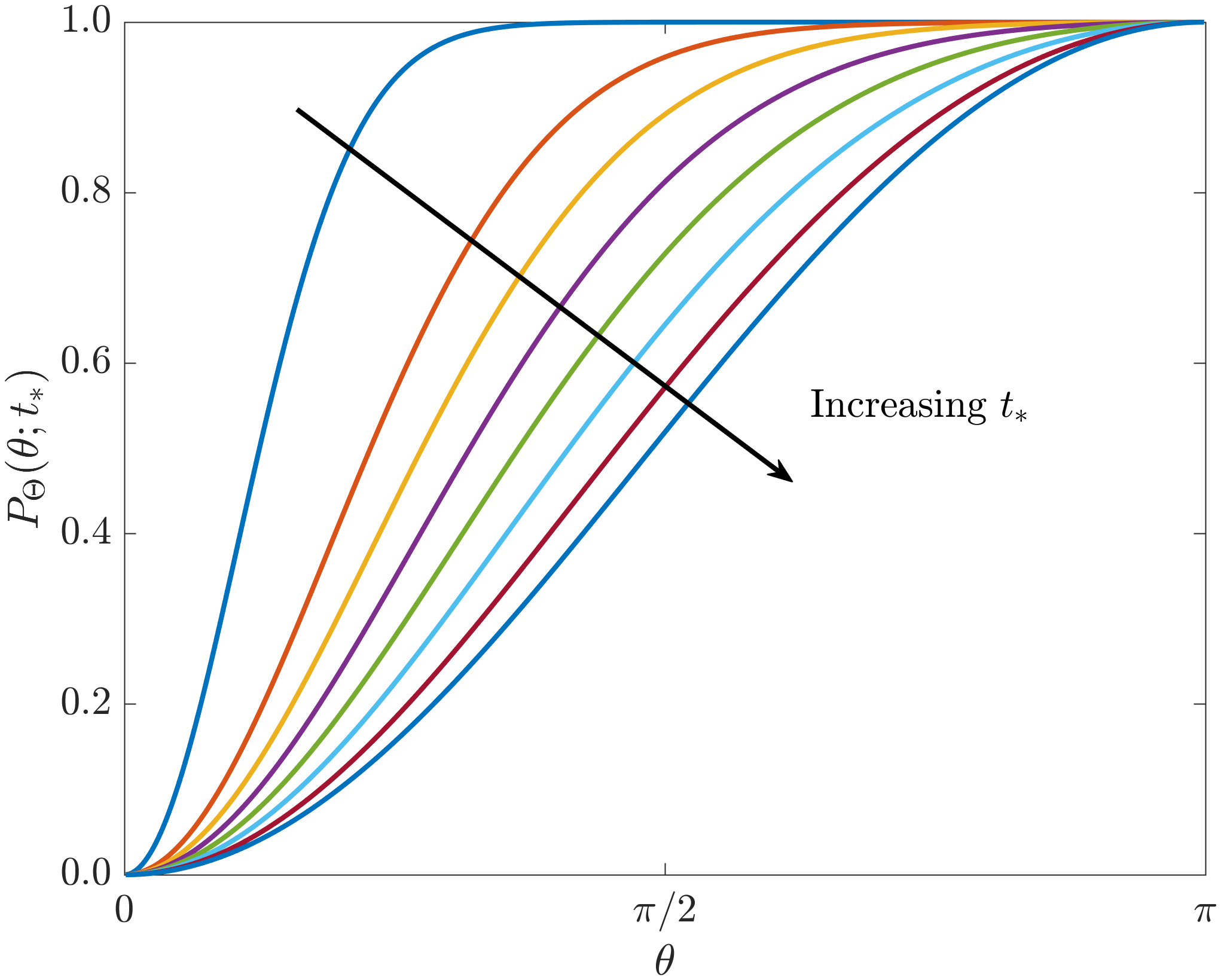}\label{fig:reinsertion_b}}
\caption{Schematic of reinsertion. Left: Consider a sphere \al{$\cS_a$ of radius $R_a$} containing all the absorbing targets and a particle (green point) external to the sphere at a distance $\al{R_0>R_a}$ from the center. The particle will escape to infinity with probability $p_{\textrm{escape}}=\al{(1-R_0/R_a)}$. Otherwise it will impact the sphere $\cS_a$ (blue point). We can describe the impact point via spherical coordinates on $\cS_a$ where $\theta_*$ measures the \ajb{polar angle, that is the} declination from the release point and $\phi_*$ is the azimuthal angle.  Right: Solution of the CDF \eqref{KMCReinsertAngle} for several values of \ajb{$t_{\ast}$,} the length of the sojourn from release to impact. The distribution of  \ajb{azimuthal} angle ($\phi_*$) is isotropic due to symmetry. However, the distribution of \ajb{polar angles} goes from being sharply peaked near zero (for small $t_{\ast}$) to uniformly distributed as $t_{\ast}\to \infty$. The  uniform distribution on the sphere corresponds to 
$P_{\Theta}( \theta; t_{\ast})=\frac12 (1-\cos\theta)$. \label{fig:app:angle_dist}  }
\label{fig:reinsertion}
\end{figure}

Suppose that the absorbers in our problem are contained within a sphere \al{$\cS_a$} of radius $\al{R_a}$ and a particle is initially is at $(x_0,y_0,z_0)$ with $\al{R_0} = \sqrt{x_0^2+y_0^2+z_0^2}$ and $\al{R_0 > R_a}$ (cf.~Fig.~\ref{fig:KMCpoly}).
The particle will eventually either escape to infinity (with probability $\al{1-R_0/R_a}$) or strike the surface of the bounding sphere \al{$\cS_a$} (with probability $\al{R_0/R_a}$). We call this propagation \emph{reinsertion} and, after a rotation, rescaling, and possible translation 
the reinsertion time $t_*$ and point $(x_*,y_*,z_*)$ can be determined by solving the equivalent problem of a particle initially on the polar axis at a point $(0,0 ,R)$ with $R>1$ external to the unit sphere as illustrated in Fig.~\ref{fig:reinsertion}. 

Similarly, there is an analogous problem in the half-space geometry as shown in Fig.~\ref{fig:KMCplane}. Here we choose a hemisphere that encloses all the absorbers and use the method of images; extend the hemisphere to a full sphere and if a particle is reinserted to the lower hemisphere below the bounding plane it is reflected to the image point on the upper hemisphere.

In Appendix \ref{sec:spherearrival} we calculate the PDF ($\JS(t$)) and CDF ($P_T(t)$)  for first arrival on the sphere to be
\begin{equation}
\JS(t) = 
\frac{R-1}{2R\sqrt{\pi D} t^{\frac32}} \exp\left[ \frac{-(R-1)^2}{4Dt}\right], \qquad P_T(t) = \int_0^t \JS(\tilt)\, d\tilt = \frac{1}{R} \mbox{erfc} \left[\frac{R-1}{2\sqrt{Dt}} \right].
\end{equation}
We remark that $\int_0^{\infty} \JS(\tilt) \, d \tilt = R^{-1}$ so that the probability of capture is not unity, but inversely proportional to $R$, the ratio of the distance to the sphere's center to the radius of the sphere. To model this process we draw a uniform random variable $\mu \in (0,1)$. If $\mu \in [R^{-1},1)$ the particle escapes to infinity and no further sampling is done. Otherwise, the arrival time, $t_*$, is
\begin{equation}
\label{KMCSphereReinsertTime}
t_* =  \frac{1}{4D} \left[ \frac{R-1}{\mbox{erfc}^{-1}(R\mu)} \right]^2.
\end{equation}
Next we obtain the arrival location $(\theta_*,\phi_*)$ on the sphere (cf.~Fig.~\ref{fig:reinsertion}) which is conditioned on the arrival time $t_*$. Heuristically, we  expect in the case of short arrival times that distribution of arrival locations is tightly focused on the north pole (which is closest to the initial location). On the other hand, if the particle undergoes a lengthy sojourn before eventual capture, the memory of the initial location is diminished and the distribution of arrival locations tends towards uniform. For a particle originating above the north pole, the distribution of the azimuthal angle $\phi_* \in [0,2\pi)$ is uniform. To sample the time dependent distribution in the polar angle (for which $\theta_* \in (0,\pi)$), we use the CDF of the \al{conditional} distribution \eqref{eq:JointSp},
\begin{equation}
\label{KMCReinsertAngle}
P_{\Theta}(\theta;t_*)   =\frac{1}{2} [1- \cos \theta]+ \frac{1}{2}\sum_{n=1}^{\infty} \frac{\chi_n(t_*)}{{\JS(t_*)}}  [ P_{n-1}(\cos\theta) - P_{n+1}(\cos\theta)].
\end{equation}
Here $P_n$ is  the $n^{th}$ Legendre polynomial and the functions $\chi_n(t)$ are described in detail in Appendix \ref{sec:spherearrival}. 
One can obtain the  angle $\theta_*$ by drawing a uniform random number $\nu \in (0,1)$ and solving for $P_{\Theta}(\theta_*;t_*,R) =\nu$.

To review, one strategy for obtaining an angle $\theta_*$, is to draw a uniform random number $ \nu \in(0,1)$ and solve $P_{\Theta}(\theta_*;t_*) = \nu $ \al{for $\theta_*$}. However, the sampling of this angular distribution is a time consuming element of the KMC algorithm because for each $t_*$ the summation in \eqref{KMCReinsertAngle} must be computed and sampled. 

\al{An alternative strategy is to tabulate a large number of time-angle pairs $(t_{\ast},\theta_{\ast})$ and draw upon them as needed during simulation. In order that a single tabulation of impact times and locations are effective for all initial radii $R_0$, we fix the ratio between the distance to the initial \emph{launch} point (green dot in Fig.~\ref{fig:reinsertion_a}) and the final \emph{landing} sphere (blue dot in Fig.~\ref{fig:reinsertion_a}). Hence, for a fixed value $R$, a particle initially at $|\bx| = R_0$ is always reinserted to a landing sphere of radius $R_0/R$ such that $R_0/R>R_a$. The potential for escape or reinsertion is sampled for any point $\bx\in\Omega$ such that $|\bx|> R\times \al{R_a}$. In the case of escape, the particle is removed from the simulation. Otherwise it is reinserted to a sphere of radius $|\bx|/R$. The corresponding time of the particle is incremented by $t_{\ast}|\bx|^2/(R^2)$.}

In practice, \al{we choose $R_a$ to be the smallest sphere bounding all absorbing faces} and the parameter value $R=3$ for which allows a high likelihood of escape ($p_{\textrm{escape}}=2/3$)
while also being likely to displace the particle significantly from the polar axis if it is reinserted. 
We tabulate $t_*$ and $\theta_*$ for an equally spaced grid of $\mu \in (0,R^{-1})$ and $\nu \in (0,1)$ (typically in a $400 \times 400$ grid). 
In our implementation we first determine if a particle escapes to infinity (that is if $\mu \in [R^{-1},1)$). 
Otherwise we randomly \al{select} an entry in this tabulation to choose the time ($t_*$) and polar angle ($\theta_*$) of the reinsertion. The azimuthal angle ($\phi_*$) is selected randomly from $(0,2\pi)$.
While this suffices for our needs, the accuracy could be improved by using interpolation on the polar angle grid and using \eqref{KMCSphereReinsertTime} to obtain a more accurate value of $t_*$.

\section{KMC Algorithms for diffusion to a plane with absorbing surface sites and to a convex polyhedron with some absorbing faces} \label{sec:KMCalg}

In this section we describe how to assemble the propagators derived in the previous section into an effective KMC algorithm in these geometries.

\subsection{A KMC method for diffusion to a plane with absorbing surface sites}\label{subsec:KMCplane}

Here we describe a three-dimensional Kinetic Monte Carlo method where particles are free to move in the half-space above the plane. The surface is reflecting except for a finite set of compact absorbing regions. Initially we identify a \ajb{disc, $\Da$, on the plane containing these absorbing sets.} 

\begin{center}
\textbf{KMC Algorithm for Diffusion to a Plane with Absorbing Surface Sites}
 \end{center}
 
\emph {Stage I: Projection from the bulk to the bounding plane.} Starting from the bulk, the particle is propagated forward to its first impact on the bounding plane \al{(Fig.~\ref{fig:KMCplane}, left)}. The location and time of impact are drawn from the exact distributions given in Sec.~\ref{subsubsec:KMCHalf}. \al{There are three possibilities at the next stage.}

\begin{enumerate}
\item[(i)] \al{The particle impacts the absorbing set. The time and location are recorded and the algorithm halts.}
\item[(ii)] \al{The particle is  close to the absorbing set}\ajb{, that is it is within the disc $\Da$}\al{. We apply Stage IIa (details below) where the particle is reinserted to the bulk by projecting to the interior surface of a hemisphere centered at the impact location and with a radius given by the shortest distance to the absorbing set (Fig.~\ref{fig:KMCplane}, center)}.
\item[(iii)] \al{The particle is  far from the absorbing set}\ajb{, that is it is outside the disc $\Da$}\al{. We apply Stage IIb (details below) to decide whether the particle escapes to infinity or is reinserted to the exterior surface of a hemisphere whose base contains the absorbing set (Fig.~\ref{fig:KMCplane}, right)}.
\end{enumerate}

\emph{Stage IIa: Projection from plane to bulk.}
\ajb{If the particle is on the bounding plane and within the disc $\Da$ that contains the absorbing sets,} calculate the distance $R$ to the absorbing set and propagate the particle forward to a random location on the hemisphere with radius $R$.  The arrival time to the hemisphere is drawn from a known exact distribution, given in Sec.~\ref{subsubsec:KMCSphere}. As the particle is \ajb{now} in the bulk, we next proceed back to Stage I.

\emph{Stage IIb: Escape or reinsertion.} If the particle is \ajb{on the bounding plane but outside the disc $\Da$ that contains the absorbing sets, KMC Projector III is used and it either escapes to infinity or is reinserted to a hemisphere whose base contains the absorbing set (cf.~Fig.~\ref{fig:KMCplane}, right). The probability of escape/reinsertion is sampled as discussed in Sec.~\ref{subsubsec:KMCReinsert}. If escape occurs, the particle is marked as having escaped and removed from the simulation. If escape does not occur, the particle is now in the bulk and we next proceed back to Stage I.}

The three stages in this process are illustrated  in Fig.~\ref{fig:KMCplane}. For a particular realization, the method alternates between Stage I (during which a particle may be absorbed) and Stage IIa or Stage IIb (during which a particle may escape to infinity).
This method can be applied to pores of general geometry; collections of circular absorbers are simplest as the calculation of the signed distance to the attracting set $\Gama$ is straightforward. This identification of the closest absorbing set scales as (number of particles) $\times$ (the number of absorbing sets) and for a large number of absorbing sets this is the most time consuming portion of the calculation.

\ajb{As discussed in Sec.~\ref{subsubsec:KMCReinsert},} for the reinsertion portion of the algorithm we have found that choosing $\cR=3$ works well \ajb{(which means we must choose the disc $D_a$ with a radius three times as large as the smallest disc that contains the absorbing set).}
In this case $1-1/\cR=2/3$ of the particles that reach this stage escape.

As a finite proportion of the particles are either captured or escape after passing through Stage I and Stage II of this algorithm, the number of surviving particles drops geometrically. Consequently, the computation scales naively as the number of particles and the independence of each particle allows for ready parallelization.

If the absorbing sets extend to infinity, such as a striped or a doubly periodic set of absorbers, the algorithm can be simplified to alternate between Stage I and Stage IIa. Details can be found in our previous work  \cite{LBS2018}.

\subsection{A KMC method for diffusion to a convex polyhedron with some absorbing faces}\label{subsec:KMCconvexpoly}

Here we describe a three-dimensional Kinetic Monte Carlo method where particles are external to a convex polyhedron whose faces are either absorbing or reflecting.  The target polyhedron could have a small number of faces, such as a cube (which we examine below in Sec.~\ref{sec:resultsCube}) or a generated triangulation approximating a sphere (Sec.~\ref{sec:sphere}), an ellipsoid (Sec.~\ref{sec:Ellipsoid}), or other convex surface.

\ajb{When the particle is in the bulk sufficiently far from the target polyhedron we can use KMC Projector III (cf. Sec.~\ref{subsubsec:KMCReinsert}) to either propagated inward to a sphere bounding the target or mark the particle as having escaped to infinity. We define a ball $\B_a$  with a radius three times as large as the smallest sphere that contains the absorbing set (i.~e.~$\cR=3$) and particles outside this sphere are subject to the escape/reinsertion propagator. Of the particles that reach this stage, $1-1/\cR=2/3$ escape.}

As compared to the half-space problem examined earlier, a new challenge is determining which face on the polyhedron to target with the KMC method. The target face plane we select is the one with the largest positive signed distance to the particle.  Here the positivity ensures that the particle and the target polyhedron are on opposite sides of the dividing face plane.  Choosing the largest distance maximizes the length of the sojourn the planar propagator will take. While this step is conceptually simple it can be the most time consuming as it scales with the (number of faces) $\times$ (number of particles). 

Once the particle is propagated to the face plane, if it actually lies in the target face it is either absorbed (if the face is absorbing) or propagated to the surface of the largest hemisphere whose base is contained in the face. The particle still lies above the same target face, so it is propagated down to the face plane again and tested to see if it still is within the target face (which occurs roughly $80 \%$ of the time). Essentially the particle executes a walk on the target face until it escapes at which point it \ajb{has returned to the bulk and subsequently} either escapes to infinity or is projected onto a new target face plane.

\begin{center}
\textbf{KMC Algorithm for Diffusion to a Convex Polyhedron with Some Absorbing Faces}
 \end{center}
 
 \emph{Stage I: Escape or reinsertion.} If the particle is sufficiently far from the polyhedron \ajb{(that is outside the ball $\B_a$)}, it undergoes reinsertion \al{and the possibility of escape is sampled using KMC projector III} as described in Sec.~\ref{subsubsec:KMCReinsert} \al{(see Fig.~\ref{fig:KMCpoly})}. Otherwise proceed directly to Stage II. If escape occurs the particle is marked as having escaped and removed from the simulation. If escape does not occur, the particle is reinserted to a sphere that encloses the convex polyhedron (cf.~Fig.~\ref{fig:KMCpoly}, right). As the particle now is in the bulk near the target, we proceed to Stage II.

 \emph {Stage II: Select a target face plane on the polyhedron.} Select the face plane that has the largest signed distance to a particle; this is the target face plane with an associated target face. Proceed to Stage III.

\emph {Stage III: Projection from the bulk to the target face plane.} Starting from the bulk, the particle is propagated forward to the first impact on the plane that contains the target face. The location and time are drawn from exact distributions, given in Sec.~\ref{subsubsec:KMCHalf}. Proceed to Stage IV.

\emph{Stage IV: Check if the particle is within the target face.}

\begin{enumerate}
\item[(i)] If the particle is outside the target face, return to Stage I.
\item[(ii)] If the particle is within the target face, and the face is absorbing, the time and the absorbing face are recorded and the algorithm halts. 
\item[(iii)] If the face is reflecting, the particle is reinserted into the bulk by the hemispherical projector, given in Sec.~\ref{subsubsec:KMCSphere}. \al{The hemisphere is centered on the impact point with radius chosen to make the largest hemisphere whose base is completely contained in the face (see Fig.~\ref{fig:KMCpoly} center).} Now, repeat Stage III \al{and project onto the same face. This guarantees that the particle leaves the vicinity of the current face before the algorithm returns to Stage I.}
\end{enumerate}

\section{Results}\label{sec:results}

In this section we will report on some examples of diffusive capture problems for collections of pores on a plane and the analogous problem for a convex surface. Our goal here is to validate the KMC methods described in the previous section by comparing to published static and dynamic results and also to exhibit the breadth of problems that can be investigated with these numerical methods.

\subsection{Arrival time distribution to pores on a plane}\label{sec:Planar}

For a planar surface with a single circular pore or two circular pores of equal radius analytical solutions are available for the capacitance problem which can be used to validate our method. The capacitance $C>0$ is a unique scalar that reflects the ability to hold an electric charge and is determined by the shape of $\partial\Omega$ together with the applied boundary conditions \cite{OKMascagni2004}. If particles are released uniformly on a hemisphere of radius $R$ which encloses the entire geometry (and whose base contains the absorbing pores), then $C = R p_{\ast}$ where $p_{\ast}$ is the capture probability which we can estimate using the KMC method.

For the time-dependent dynamics we find close agreement with our asymptotic and homogenization formulae described in Appendix \ref{sec:asy}. We consider an example with six pores to show the versatility of both the KMC method and our asymptotic approximations.

\subsubsection{Single circular pore} 
\begin{figure}[htbp]
\centering
\subfigure[Error in $C$ for $N=1$ pore.]{\includegraphics[width = 0.31\textwidth]{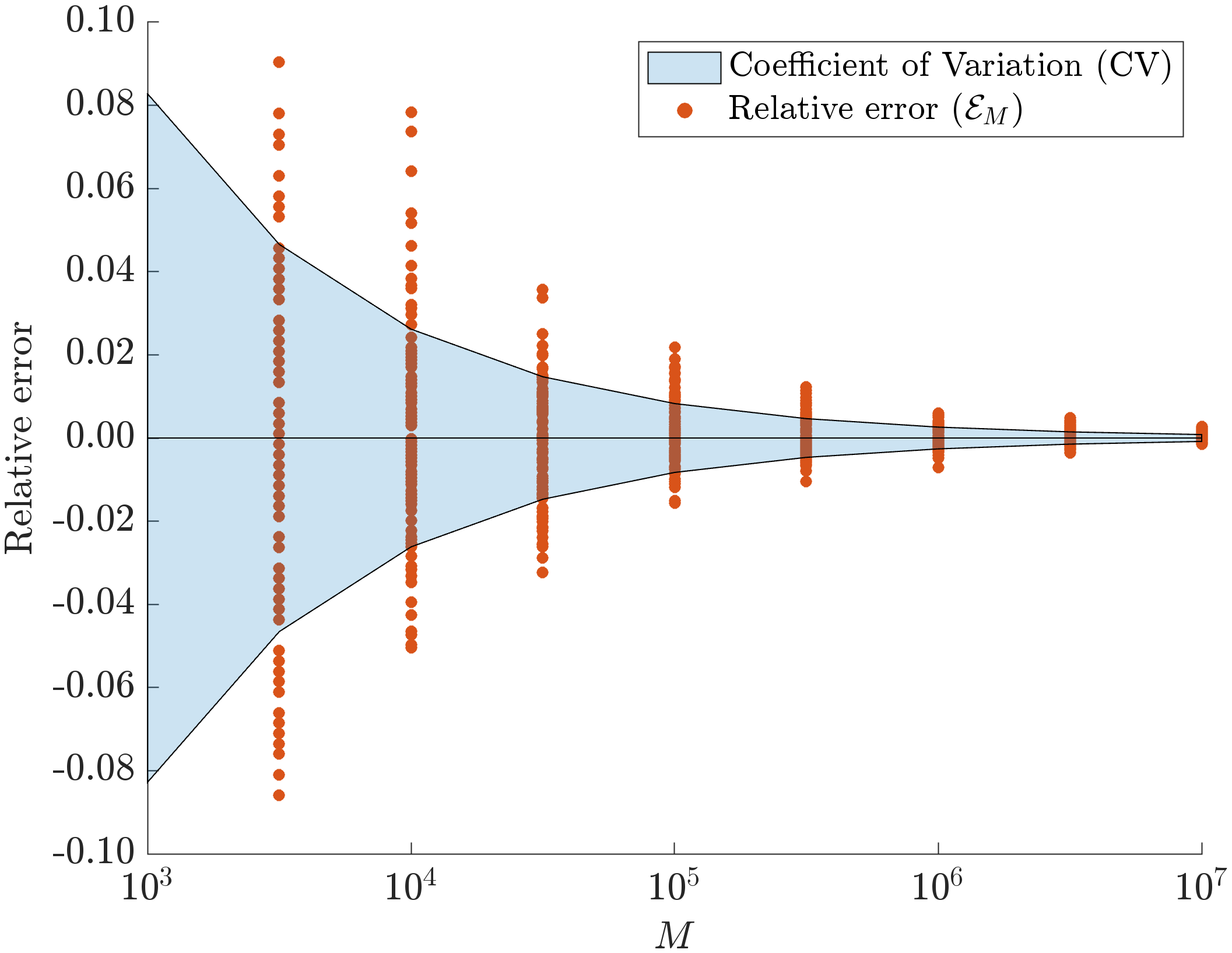} \label{fig:planarOnePore}} \quad
\subfigure[Single pore capture distribution.]{\includegraphics[width = 0.31\textwidth]{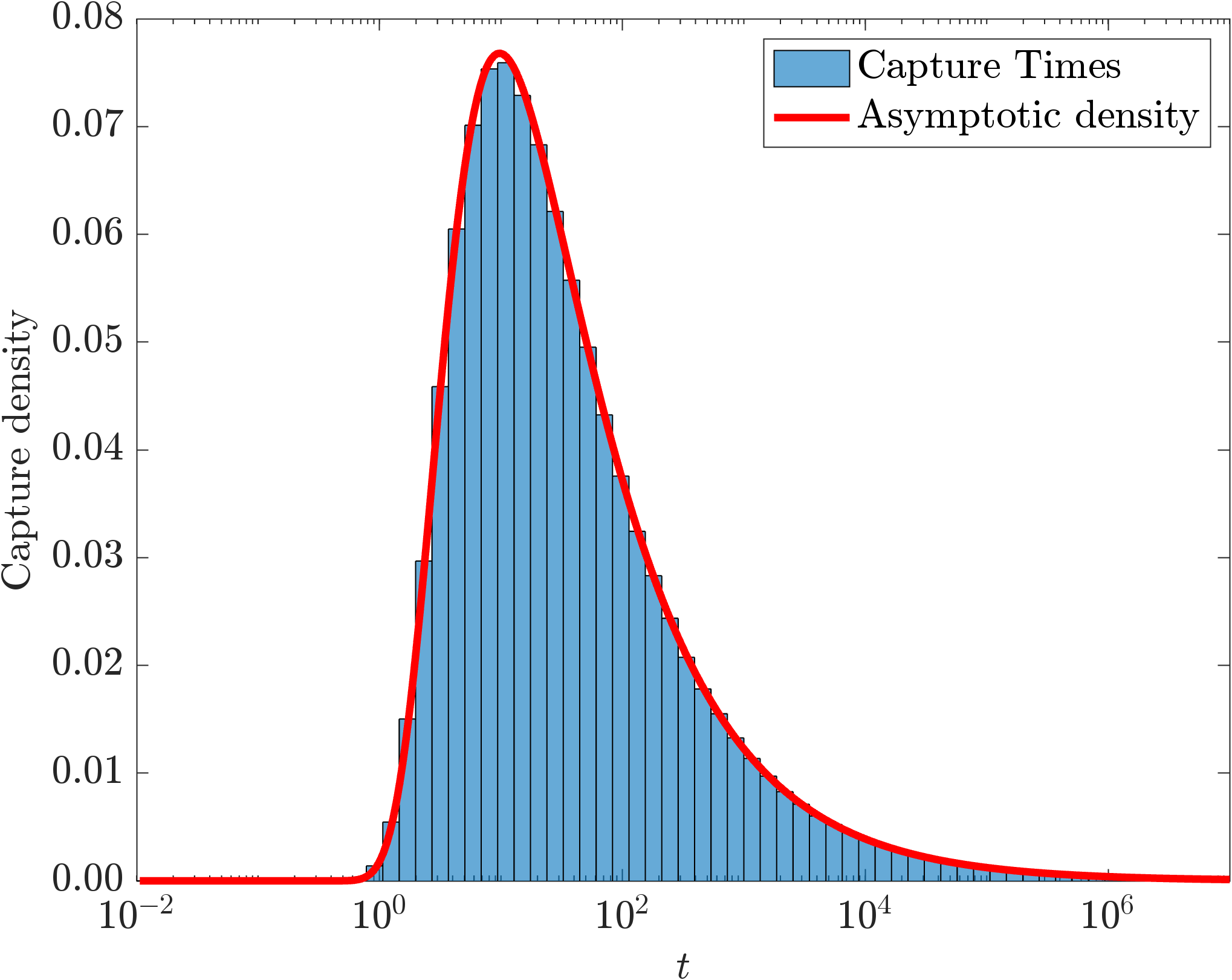}\label{fig:planarOnePorePDF}} \quad
\subfigure[Capacitance for $N=2$ pores.]{\includegraphics[width = 0.31\textwidth]{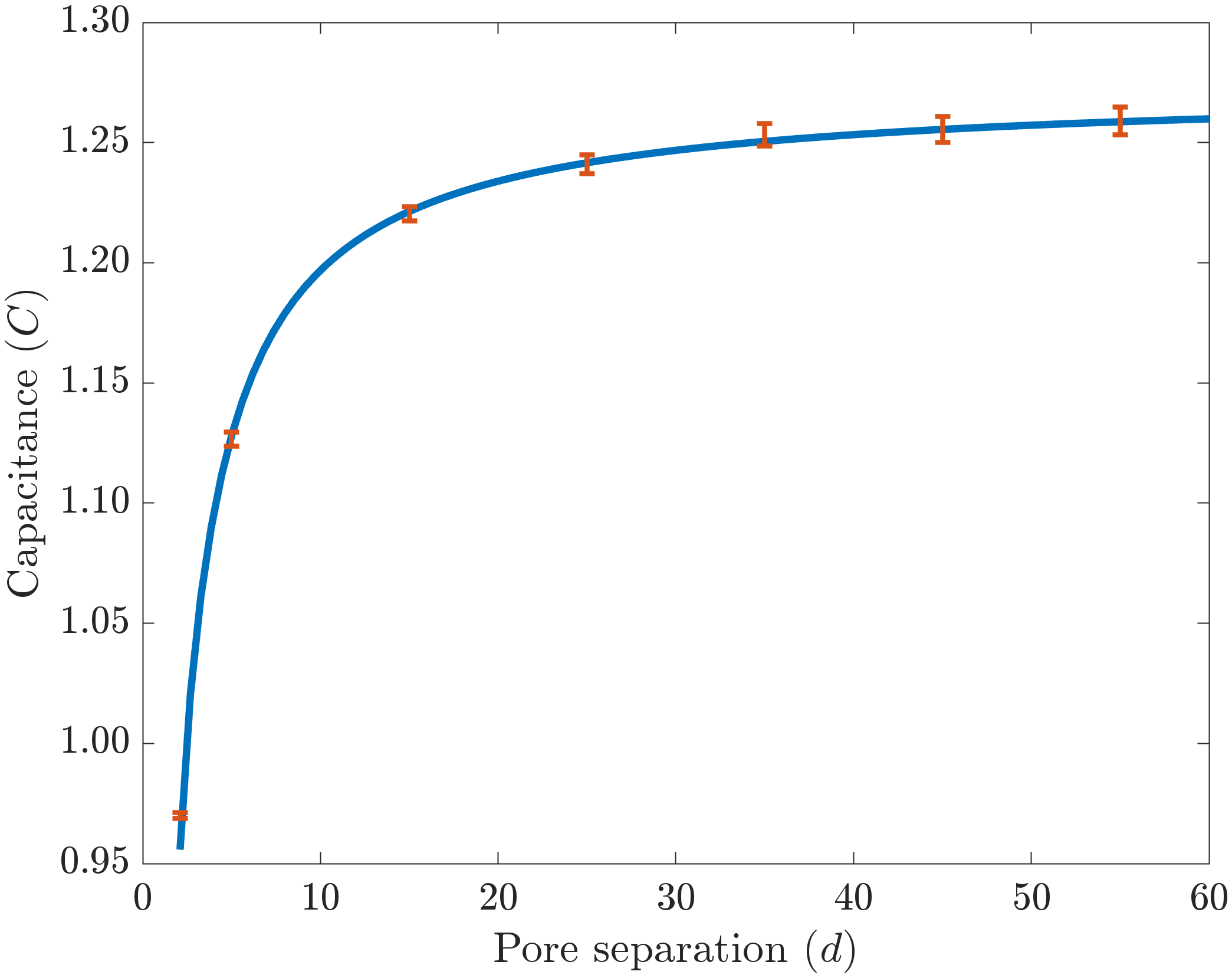}\label{fig:planarTwoPore}}
\caption{Capacitance and dynamics of capture by absorbing circular pores on the plane. Panel (a): Convergence of estimates for the capacitance of a single pore based on KMC trajectories for $M$ particles. Panel (b): Logarithmic histogram of capture times for a single ($N=1$) circular pore 
with $M= 10^7$ trajectories of diffusivity $D=1$ initialized from $\bx_0=(0,0,5)$. The red curve is the asymptotic approximation obtained from (\ref{eq:asyCDFone}).
Panel (c): Agreement between the two pore capacitance formula (solid) \eqref{eq:2PoreCap} and KMC (red) over various pore separations, \al{$d$}, for $M=10^8$ particles. Error bars indicate CV \eqref{CV} obtained from KMC simulations.}
\end{figure}

We first validate our method by considering the capacitance of a single circular pore of unit radius on the plane. The known capacitance is $C=c_0= 2/\pi$ and we confirm this value by computing  KMC trajectories for $M$ particles initialized uniformly on a hemisphere of radius $R=5$. The probability of capture is exactly
$p_* = C/R$ and the coefficient of variation of $p_*$ (for $M$ particles) is given by 
\begin{equation}
\text{CV} = 
\frac{\text{standard deviation}}
{\text{expected value}} = 
\frac{1}{p_*}
\sqrt{\frac{p_*(1-p_*)}{M}} =
\sqrt{\frac{(1-p_*)}{p_*M}} .
\label{CV}
\end{equation}
We estimate the probability of capture, $\KMCProb$, as the ratio of the number of capture particles to the total number of particles and then compare the relative error $\E_M = (\KMCProb-p_*)/p_*$ to the coefficient of variation.
In Fig.~\ref{fig:planarOnePore} bootstrap resampling with $100$ replications is used to generate additional estimates as a function of $M$ and confirm the expected $\bigoh(M^{-\frac12})$ convergence of our KMC method. 
For the time-dependent problem, we derived the asymptotic estimate
(\ref{eq:asyCaptureOne}) for the PDF of arrival times in Sec.~\ref{sec:asyplane}
\begin{equation}\label{eq:asyCaptureOne_Main}
\J(t)   \sim  \frac{1}{\sqrt{D} (\pi t)^{3/2}  } e^{-\frac{R^2}{4Dt} }  \left [ 1 - \left( \frac {2}{\pi R} -\frac{R} {D \pi t} \right ) \right ],
\end{equation}
where $R = |\bx_0|$. The associated CDF for this density (\ref{eq:asyIntCaptureOne}) was also derived yielding
\begin{equation}
\label{eq:asyCDFone}
F(t) \sim  \frac{1}{ R} \mbox{erfc} \left (\frac{R}{2 \sqrt{Dt}} \right )
+ \frac{1}{R}
\frac{e^{-\frac{R^2}{4Dt} } }{\sqrt{\pi Dt}} \ .
\end{equation}
These expression are valid in the limit when $R$ is much larger than the pore size and agree well with our KMC results as seen in Fig.~\ref{fig:planarOnePorePDF}. This histogram needs a small note of explanation; to see the details of the distribution including the exponentially small initial captures and the slow decaying algebraic tail a logarithmic scale is appropriate. The histogram bins are chosen to be equally spaced in $\log_{10} (t)$ and the red curve represents the asymptotic estimate of the capture count in each bin determined via the CDF \eqref{eq:asyCDFone}.

\subsubsection{Two circular pores of equal radius} Strieder employed bi-polar coordinates to obtain expressions for the capacitance of two pores on the plane \cite{Strieder08,Strieder12}. See also equation (2.23) of \cite{BL2018}. For two circular pores of unit radius and separation $d\in(2,\infty)$, the capacitance has the series solution
\begin{equation}\label{eq:2PoreCap}
C_{\textrm{Str}}(d) =  \frac{4}{\pi} \left[ 1 - \frac{2}{\pi d} + \frac{4}{\pi^2 d^2} -  \frac{2(12+\pi^2)}{3\pi^2 d^3} + \frac{16(3+\pi^2)}{3\pi^4 d^4} - \frac{4(120+ 70\pi^2 + 3\pi^4)}{15\pi^5 d^5}\right] + \mathcal{O}(d^{-6}) , \quad d \to\infty.
\end{equation}
In Fig.~\ref{fig:planarTwoPore} we plot favorable agreement between the two pore capacitance formula \eqref{eq:2PoreCap} and KMC simulations based on $M= 10^8$ trajectories.

\subsubsection{A cluster of six circular pores} 
\begin{figure}[htbp]
    \centering
    \subfigure[Capture distributions.]{\includegraphics[width = 0.475\textwidth]{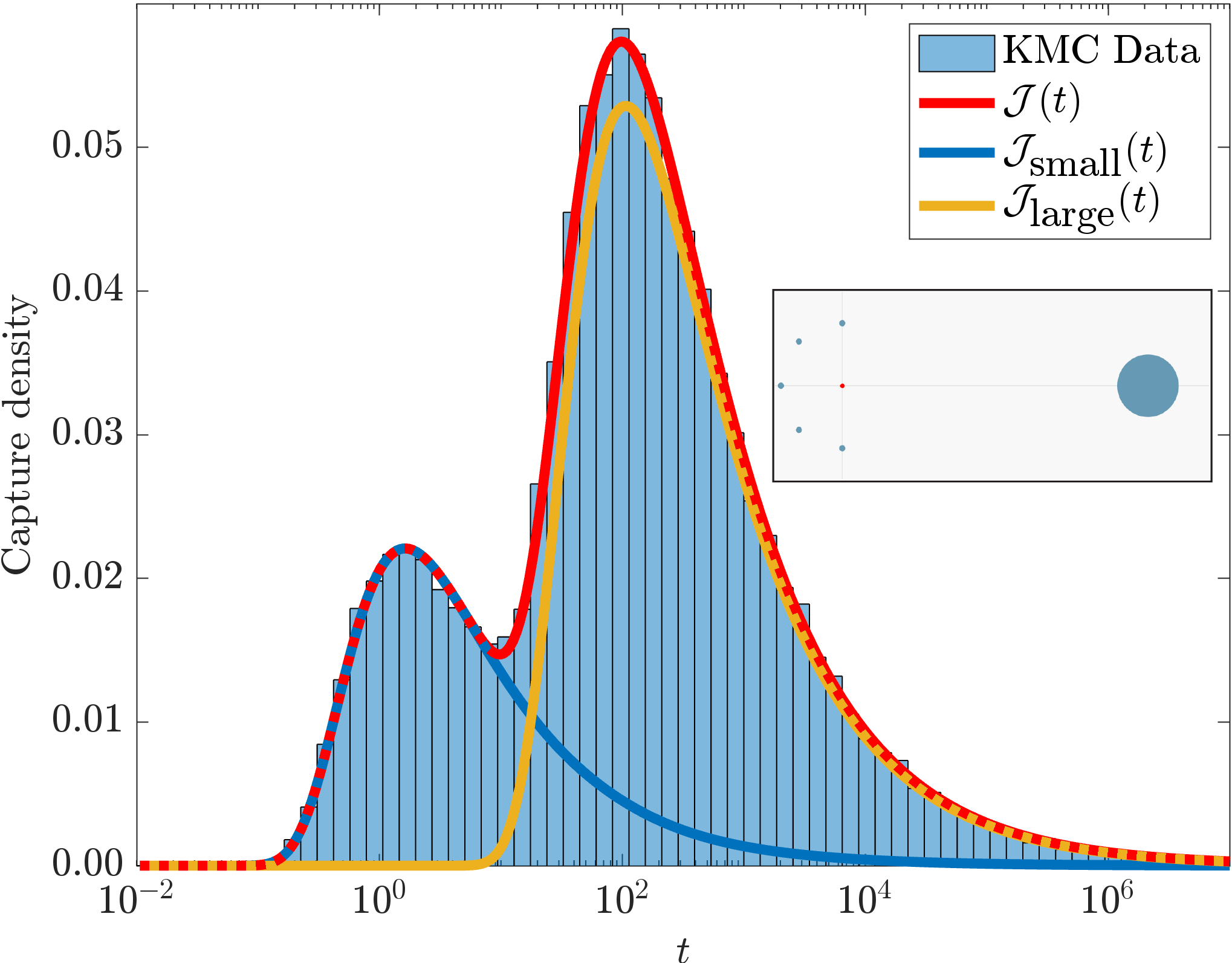}\label{fig:planarSixPore_a}}\qquad
    \subfigure[Capture fractions.]{\includegraphics[width = 0.475\textwidth]{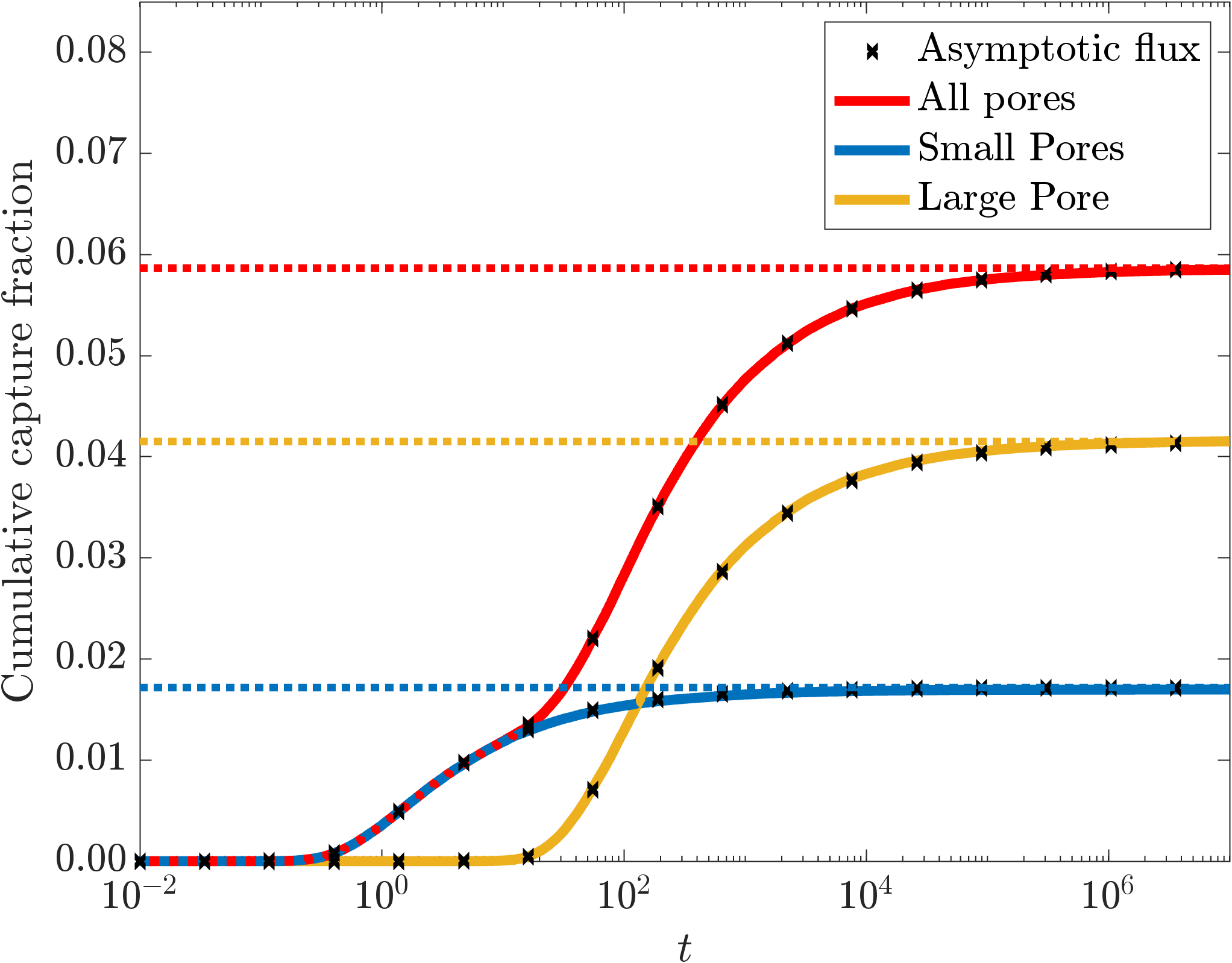}\label{fig:planarSixPore_b}}
    \caption{The 6-pore example \eqref{eq:6pore} where the combination of small and large pores generates a bimodal capture time distribution. Data from $M=10^7$ trajectories of diffusivity $D=1$. Panel (a):  Logarithmic histogram of capture times for the small pores combined, $\J_{\textrm{small}}(t)$; the single large pore, $\J_{\textrm{large}}(t)$, and \ajb{all pores} combined $\J(t)$. Schematic of the trapping configuration is shown inset. 
    Panel (b): The corresponding cumulative capture distribution as predicted from KMC data (solid lines) and the asymptotic approximation \eqref{asy:CDF}. The limiting capture fraction (dashed line) is predicted from the splitting probabilities \eqref{eq:splittingProbs}.}
    \label{fig:planarSixPore}
\end{figure}

Here we consider a six pore example with centers and radii given by
\bsub\label{eq:6pore}
\begin{gather}
\bx_k  = (\cos \theta_k, \sin \theta_k,0), \quad \theta_k = \frac{\pi}{2} + (k-1) \frac{\pi}{4}, \quad r_k = 0.01, \quad k = 1,2,3,4,5.\\[4pt]
\bx_6 = (15,0,0), \qquad r_6 = 1.0, \qquad \bx_0 = (0,0,0).
\end{gather}
\esub
The rationale behind the construction of this example is to explore the competition between several small and near pores to a single large but distant pore. In Fig.~\ref{fig:planarSixPore} we plot the arrival time distributions for $M=10^7$ KMC trajectories released from the origin with diffusivity $D=1$. 

In Fig.~\ref{fig:planarSixPore_a} we plot the combined flux \ajb{density} $\J(t)$ as predicted by KMC data (histogram) and asymptotics \eqref{eq:MainPlanarResult} (solid red) with very good agreement observed. The competition for flux from the trapping array gives rise to a bi-modal capture distribution. To explore this further, we plot the capture rates $\J_{\textrm{small}}(t) = \sum_{k=1}^5\J_k(t)$ and $\J_{\textrm{large}}(t) = \J_6(t)$ to the small pores and the large pore respectively. These curves indicate that the earlier short peak represents capture at the small pores while the later large peak represents capture at the single large pore. To understand the combined capture fraction at each set of absorbers, we plot in Fig.~\ref{fig:planarSixPore_b} the cumulative capture fractions associated with the distributions in Fig.~\ref{fig:planarSixPore_a}. In addition, we plot (dashed lines) the splitting probabilities $\Q_{\textrm{small}}(\bx_0)= \sum_{k=1}^5\Q_k(\bx_0)$ and $\Q_{\textrm{large}}(\bx_0)= \Q_6(\bx_0)$ representing the cumulative capture fraction at the five small pores and single large pore, respectively. In Appendix \ref{eq:split_final} we derive asymptotic expressions for the splitting probabilities $\Q_k(\bx_0)$. The result and its reduction for the values in \eqref{eq:6pore}, are given by 
\begin{equation}\label{eq:splittingProbs} 
\nonumber \Q_k(\bx_0) \sim \frac{c_k}{|\bx_0-\bx_k|} -  \sum_{\substack{j =1\\ j\neq k} }^N \frac{c_j c_k}{|\bx_j - \bx_k | |\bx_0- \bx_j |} 
=\frac{2}{\pi}\frac{r_k}{|\bx_k|} - \frac{4}{\pi^2}\sum_{\substack{j =1\\ j\neq k} }^6 \frac{r_j r_k}{|\bx_j - \bx_k | |\bx_j |}.
\end{equation}
For the parameters \eqref{eq:6pore}, we calculate that 
\[
\frac{\Q_{\textrm{large}}(\bx_0)}{\Q_{\textrm{small}}(\bx_0)}\approx 2.4175,
\]
indicating that the single large pore captures more than twice that of the  five small pores combined.

\subsection{Arrival time distribution to cube}\label{sec:resultsCube}

A long standing computational challenge problem has been to approximate the capacitance of the cube \cite{WINTLE200451,Higgins51,OKMascagni2004} for which no closed form expression is known to exist. Traditional finite element and difference approximations are hampered by the unbounded domain and the corner singularities.

\begin{figure}[htbp]
\centering
\subfigure[Convergence of capacitance with $M$.]{\includegraphics[width = 0.45\textwidth]{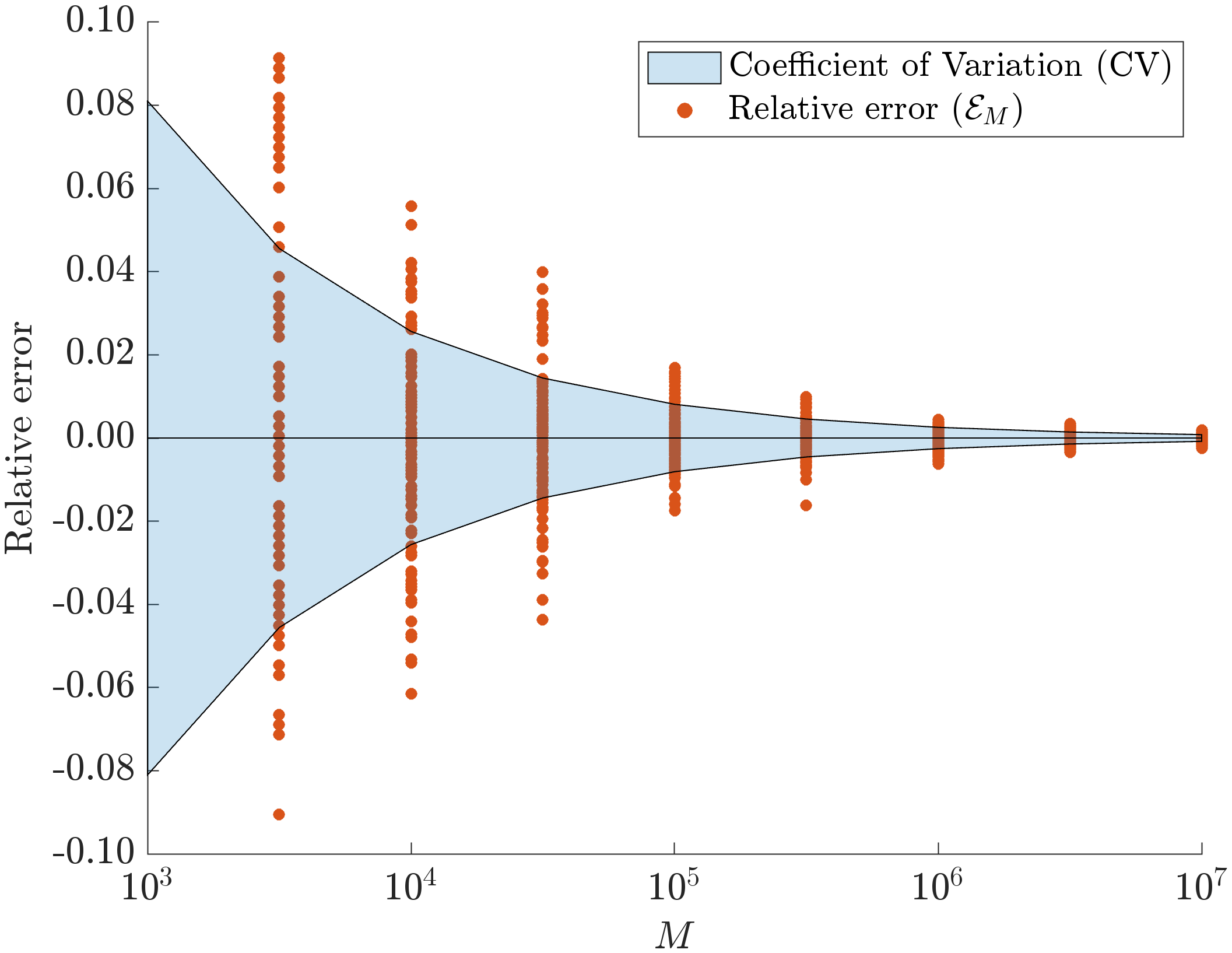} \label{fig:cubeerror}}\qquad
\subfigure[CDF \eqref{eq:CDF_Cube} of cube capture times compared to the effective sphere.]{\includegraphics[width = 0.45\textwidth]{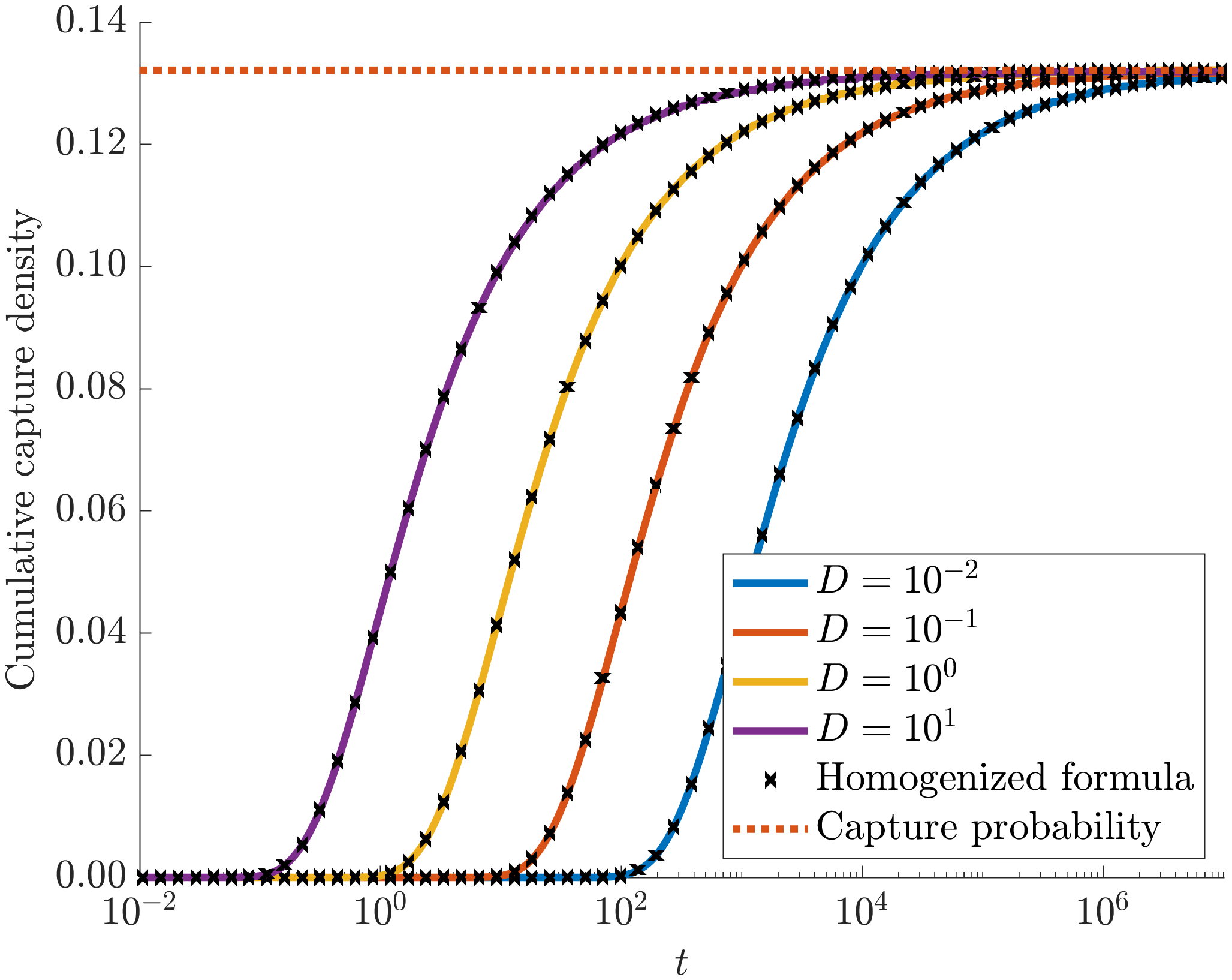} \label{fig:cubepdf}}
\subfigure[Logarithmic histogram of cube capture times compared to the effective sphere.]{\includegraphics[width = 0.45\textwidth]{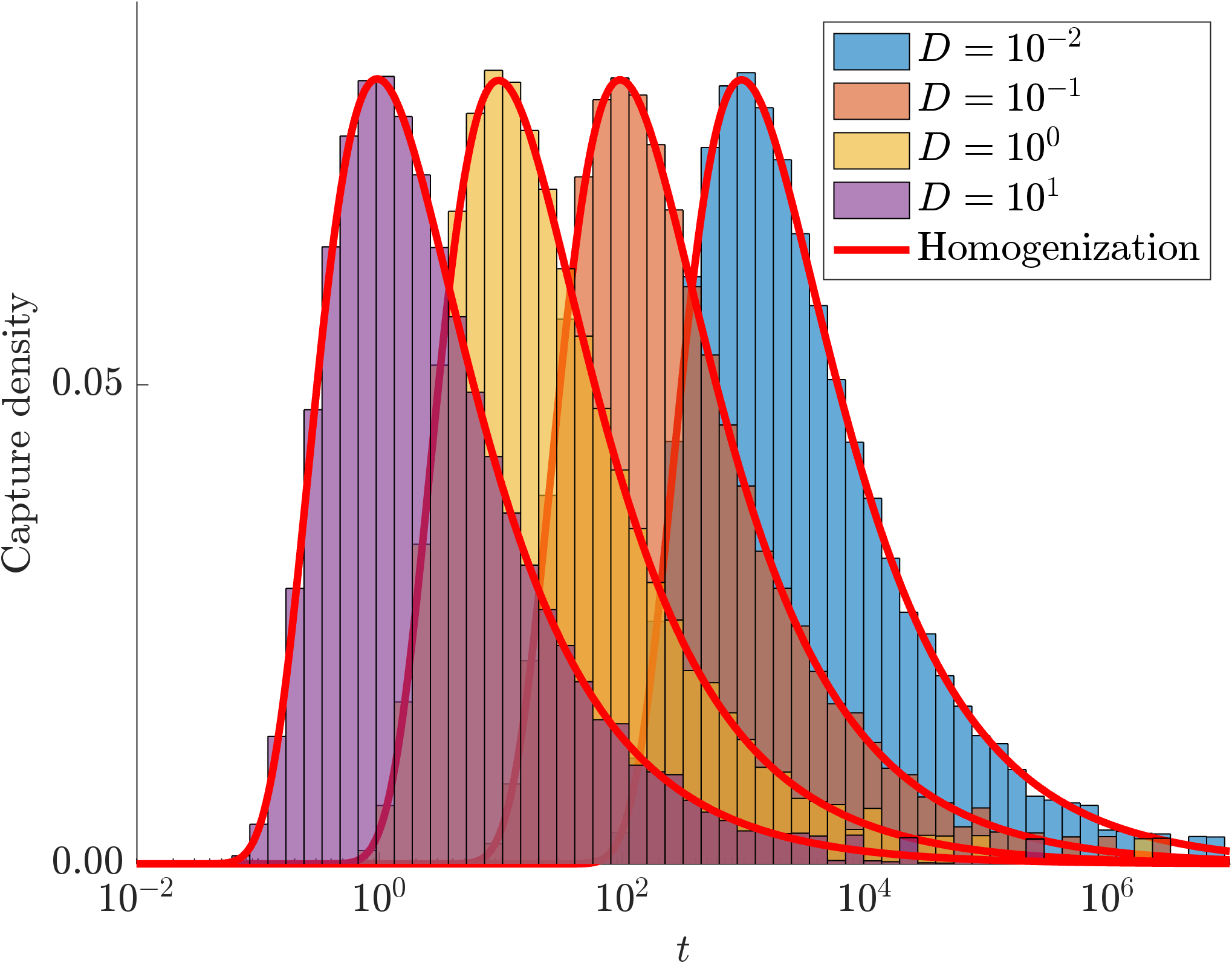} \label{fig:cubecdf}}\qquad
\subfigure[Relative error in homogenized CDF ($D=10^1$).]{\includegraphics[width = 0.45\textwidth]{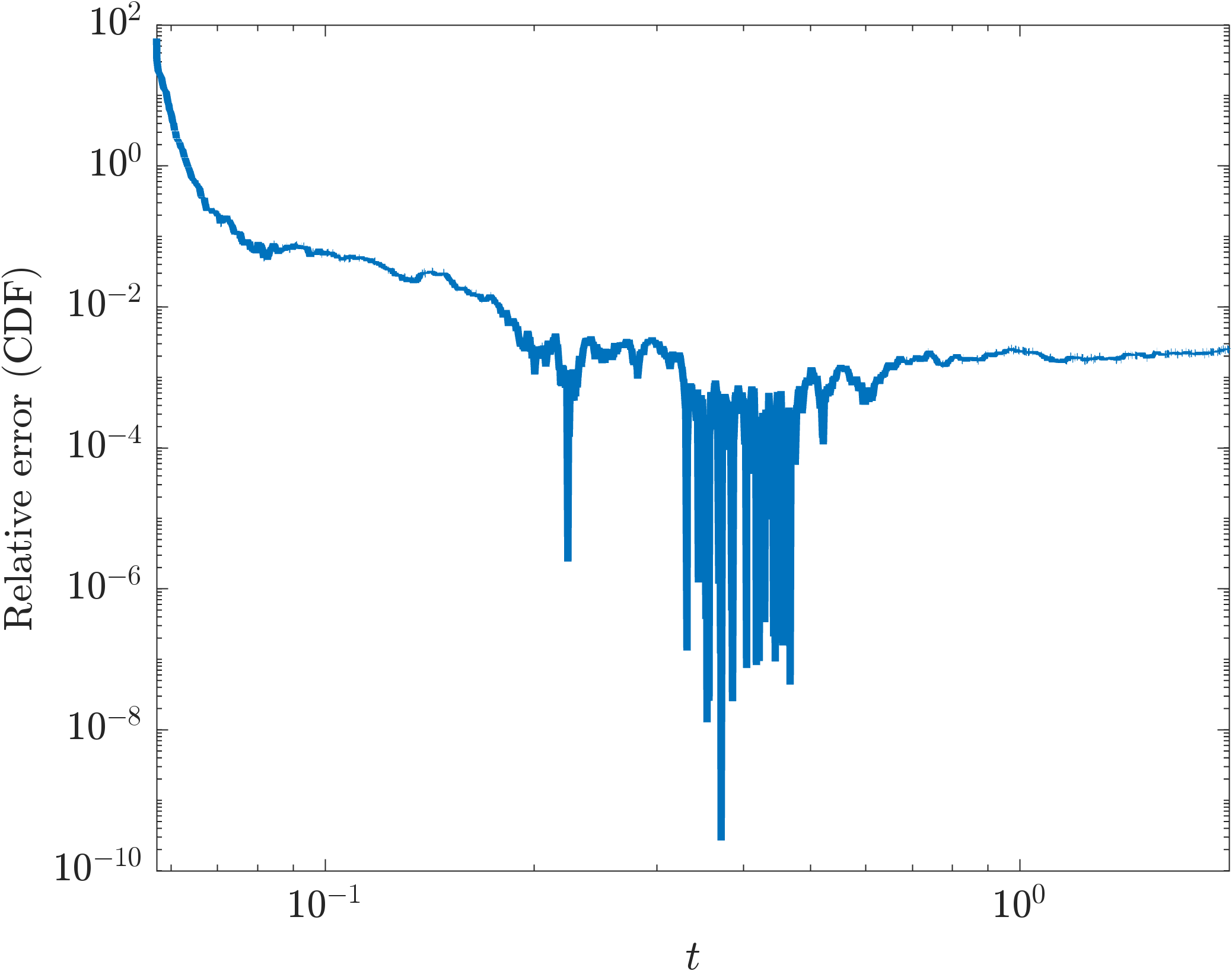} \label{fig:cubeerrorCDF}}
\caption{Arrivals to the unit cube. Panel (a): The approximation of the cube capacitance converges with the number of KMC trajectories, $M$. Particles are released at uniform points on a sphere of radius $R=5$. Replacement of the cube with a sphere of equivalent capacitance yields a highly accurate representation of the cumulative capture fractions \eqref{eq:CDF_Cube} (panel (b)) and capture densities (panel (c)) for various diffusivities. KMC simulations shown with $M = 10^7$ trajectories initialized at position $\bx_0=(0,0,5)$. In panel (d) we plot the relative error between the CDF \eqref{eq:CDF_Cube} from KMC trajectories and homogenization for $D=10^1$. We observe a small error, except when $t\ll1$ indicating homogenization may not accurately approximate early arrivals.}
\end{figure}

The integral equation formulation of \cite{HELSING2013445} estimated the value of $C_{\textrm{cube}}$ as
\begin{equation}\label{eq:CapExact}
C_{\textrm{cube}} = 0.66067815409957 \pm 10^{-13},
\end{equation}
which we use as a validation of our method. We remark that \cite{HWANG20101089} estimated $C_{\textrm{cube}}$ to a relative error of $10^{-7}$ with $M= 5\times 10^{13}$ walk-on-sphere trajectories in an approach similar to that presented here. Our first experiment is to verify convergence to the capacitance \eqref{eq:CapExact}. We sample $M = 10^8$ trajectories initiated at uniformly distributed points on a sphere of radius $R=5$. The estimate $C_{\textrm{KMC}} = R \KMCProb$ is formed from the KMC method where $\KMCProb$ is the fraction of particles which are captured. We obtain the estimate
\[
C_{\textrm{KMC}} = 0.6606454 \pm 1.7\times 10^{-4},
\]
which agrees with previously obtained values. Figure 
\ref{fig:cubeerror} studies the convergence to this capacitance as a function of the number $M$ of KMC trajectories. We once again use bootstrap resampling to estimate the error in our approximation (cf.~Fig.~\ref{fig:planarOnePore}) and convergence with an $\bigoh(M^{-\frac12})$ error scaling is evident.

\paragraph{Homogenizing the cube: a spherical approximation} A significant simplification can be obtained by replacing complex geometries with spherical ones with appropriately chosen radius $R$. A natural choice for $R$ is to choose a sphere that reproduces the capacitance of the original object. To verify the accuracy of such a process, we compute arrivals to the cube based on $M= 10^7$ KMC trajectories initialized at $\bx_0=(0,0,5)$ and plot the CDF of the equivalent spherical distribution \eqref{eq:CDFsphere}
with $R=C_{\textrm{cube}}$ and $R_0 = |\bx_0|$, given by
\begin{equation}\label{eq:CDF_Cube}
F(t) = \frac{R}{R_0} \text{erfc}\Big[ \frac{R_0-R}{2\sqrt{Dt}} \Big].
\end{equation}
In Fig.~\ref{fig:cubepdf}-\ref{fig:cubecdf} we plot comparisons between KMC simulation data and the homogenized expressions for capture density. These agree extremely well suggesting that an effective sphere condition is a useful simplification provided the radius is chosen appropriately.

Finally, in Fig. \ref{fig:cubeerrorCDF} we plot the relative error between the CDF from the KMC trajectories and the spherical homogenization for $D=10^1$. We observe a small error, except when $t\ll1$.  This can be understood by realizing that early arrivals represent nearly ballistic trajectories \cite{Lindsay2023a,Lindsay2023b} that directly impact the absorbing cube and therefore feel the nonspherical geometry at leading order.

\subsection{Arrival distribution to the unit sphere}\label{sec:sphere}
In this section we investigate and validate results for the unit sphere with small absorbing windows. As with previous scenarios, exact solutions to the time dependent problem with both absorbing and reflecting potions are few and far between. Hence we first compare with previously derived asymptotic results of the steady state problem, namely the capacitance and the splitting probabilities which provide a framework for comparison and validation. The asymptotic result for the capacitance allows us to replace the absorbing pores and reflecting complement with a homogenized Robin boundary condition which can be solved exactly yielding an asymptotic approximation for the capture PDF and CDF which we verify with our KMC method.

In particular, we consider the diffusion problem \eqref{eqn:IntroP}
where $\partial\Omega$ is the unit sphere and $\Gama = \cup_{k=1}^N\{\partial\Omega_k\}$ is the union of $N$ small non-overlapping locally circular pores with centers $\bx_k$ and radii $a_k$. The reflecting portion of the domain is given by $\Gamr = \partial\Omega\setminus\Gama$. In spherical coordinates, the location and extent of the pores are given by
\bsub\label{eqn:Pore_extent}
\begin{gather}
\label{eqn:Pore_extent_a} \bx_k = (\sin\theta_k\cos\psi_k,\sin\theta_k\sin\psi_k,\cos\theta_k),\\[4pt]
\label{eqn:Pore_extent_b} \partial\Omega_k \equiv \{  (\theta,\psi) \, | \, (\theta-\theta_k)^2 + \sin^2(\theta_k) (\psi-\psi_k)^2 \leq a_k^2\},
\end{gather}
\esub
so that $| \partial \Omega_k| \approx\pi a_k^2$. The are two special cases where we can connect with established results.

\paragraph{Homogenization} In \cite{LWB2017}, the capacitance of a unit sphere with $N$ circular pores of common radius $a$ was determined asymptotically as $a\to 0$. This was used to obtain an effective boundary condition in the limit $N\to\infty$, $a\to 0$ with the area fraction $\surfrac = \frac{N\pi a^2}{4 \pi} = \frac{Na^2}{4}$ fixed. For uniformly distributed pores, it was proposed that the mixed boundary conditions \eqref{eqn:IntroP_b} be replaced by the Robin condition \cite{LWB2017} 
\begin{equation}\label{eq:RobinBC}
D\nabla p \cdot \hn = \kappa p, \quad |\bx| = 1, \qquad \kappa = \frac{4D\surfrac}{\pi a}\left[1- \frac{4}{\pi} \sqrt{\surfrac} + \frac{a}{\pi}\log(4e^{-1/2}\sqrt{\surfrac})\right]^{-1},
\end{equation}
\al{which is valid when $|\bx_0|$ is sufficiently large and provided  $\sqrt{Dt} \gtrapprox \bigoh(1)$.}

The PDE \eqref{eqn:IntroP} together with the homogenized boundary condition \eqref{eq:RobinBC} allows for analytical solution. The distribution of the arrival times \cite{Lindsay2023b} is given by
\begin{equation}\label{eq:HomogSphereFlux}
\J(t) = \frac{\kappa}{R} 
e^{-\frac{\left(R -1\right)^2}{4 \D t}}
\left[ \frac{1}{\sqrt{\pi  \D  t}}
-\mathrm{erfc}\! \left(\beta \right)
 {\mathrm e}^{\beta^2}  \left(\kappa/D +1\right)\right],
\end{equation}
where $\beta= \frac{R -1}{2 \sqrt{Dt}}+\left(\kappa/D +1\right) \sqrt{Dt} $. The CDF of this distribution can then be calculated as
\begin{equation}\label{eq:HomogSphereCDF}
F(t) = \int_0^t \J(\tilt) \, d\tilt 
= \frac{1} {\left(1+D/\kappa \right) R}
 \left[ \mathrm{erfc}\! \left(\frac{R -1}{2  \sqrt{Dt}}\right)
 -\mathrm{erfc}\! \left(\beta \right)  e^{\beta^2} e^{-\frac{\left(R -1\right)^2}{4 \D t}}
 \right ] \ .
\end{equation}
We remark that 
\begin{equation}\label{eq:totalfractCapture}
\int_0^{\infty} \J(\tilt) \,d\tilt = \lim_{t \to \infty} F(t) = \frac{1}{(1+D/\kappa)R},
\end{equation}
so that the probability of capture is not unity, but inversely proportional to the initial distance to the sphere.

\begin{figure}[htbp]
\centering
\subfigure[Sphere with $N=51$ absorbing pores.]{\includegraphics[width = 0.35\textwidth]{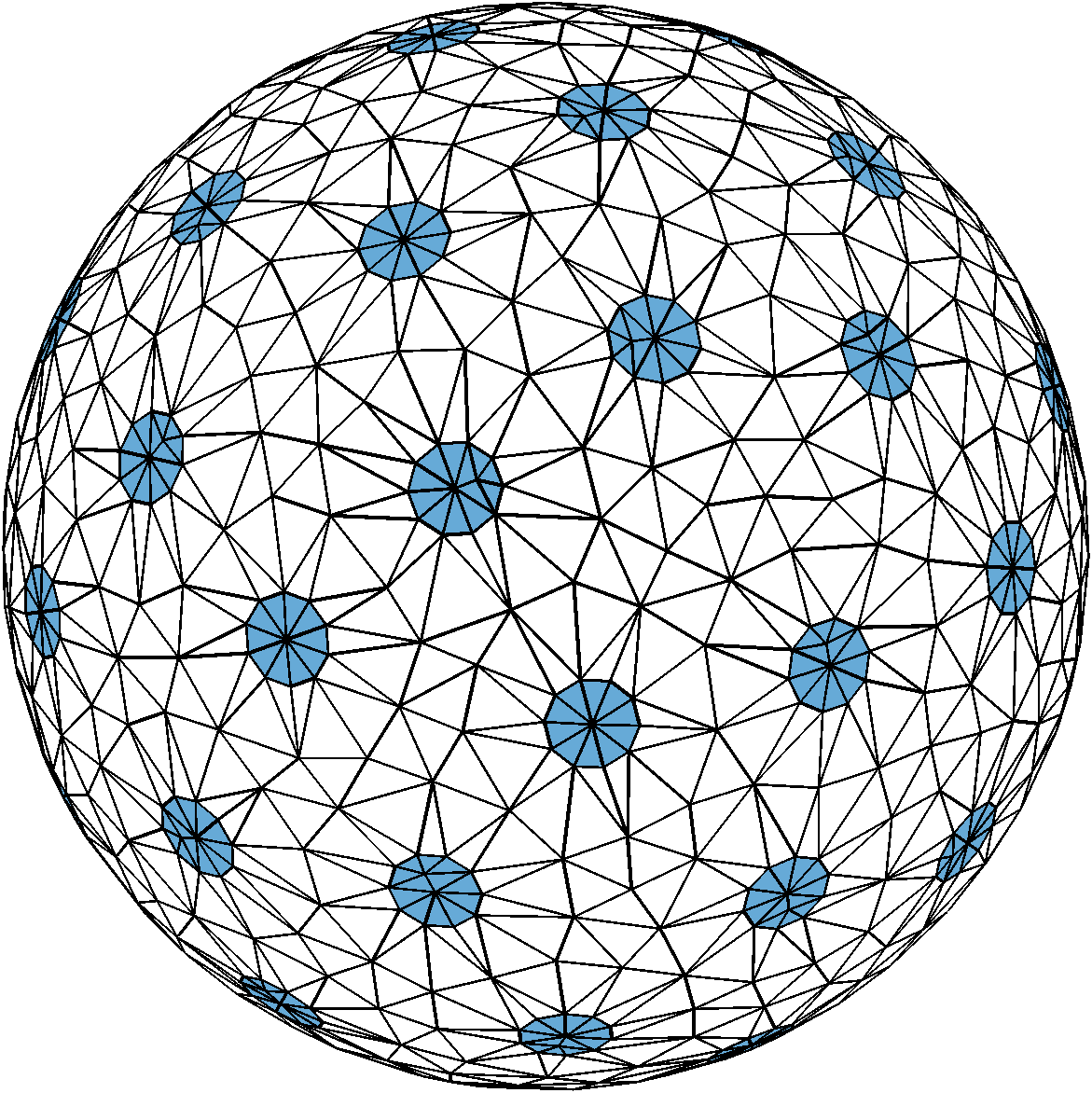} \label{fig:SphereDiff_a}}\hspace{2cm}
\subfigure[Capture time densities for $D = 10^{-2}, 10^{-1},10^{0}.$]{\includegraphics[width = 0.45\textwidth]{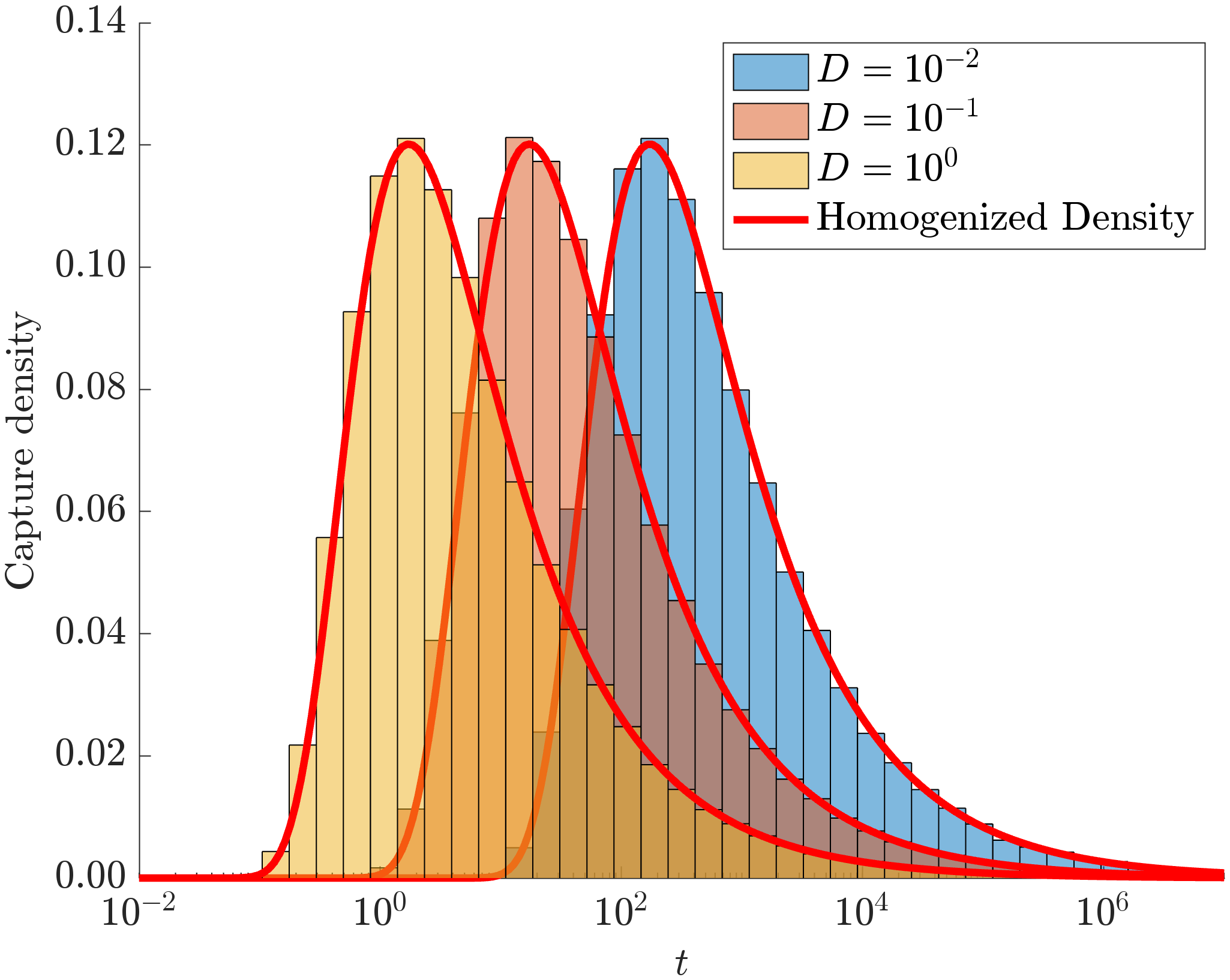}}\\
\subfigure[Total captured fraction with time for $\surfrac = 0.1$ and diffusivities $D = 10^{-2}, 10^{-1},10^{0}.$]{\includegraphics[width = 0.45\textwidth]{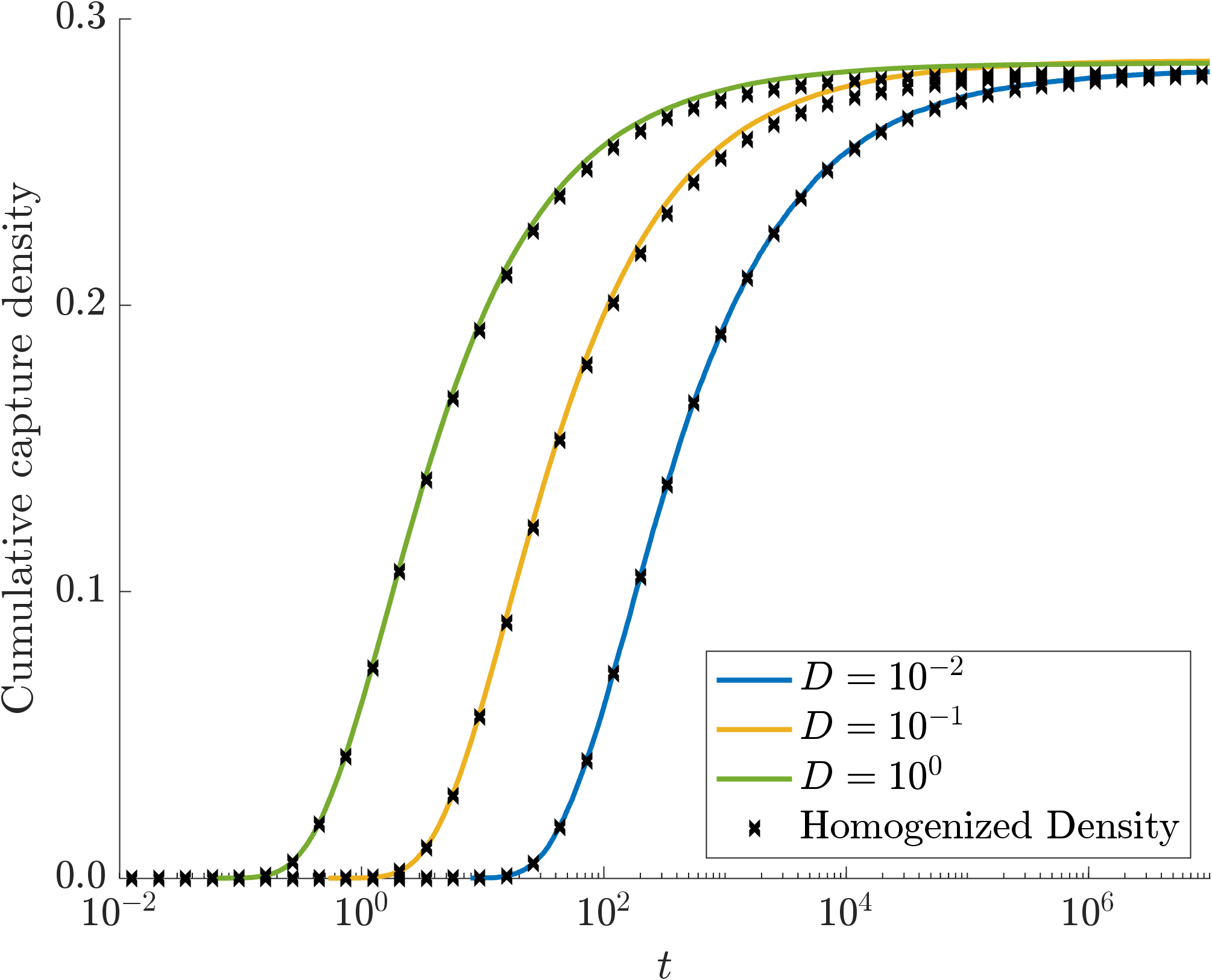}}\qquad
\subfigure[Captured fraction for $D=1$ and coverage fractions $\surfrac = 0.2,0.1,0.05$.]{\includegraphics[width = 0.45\textwidth]{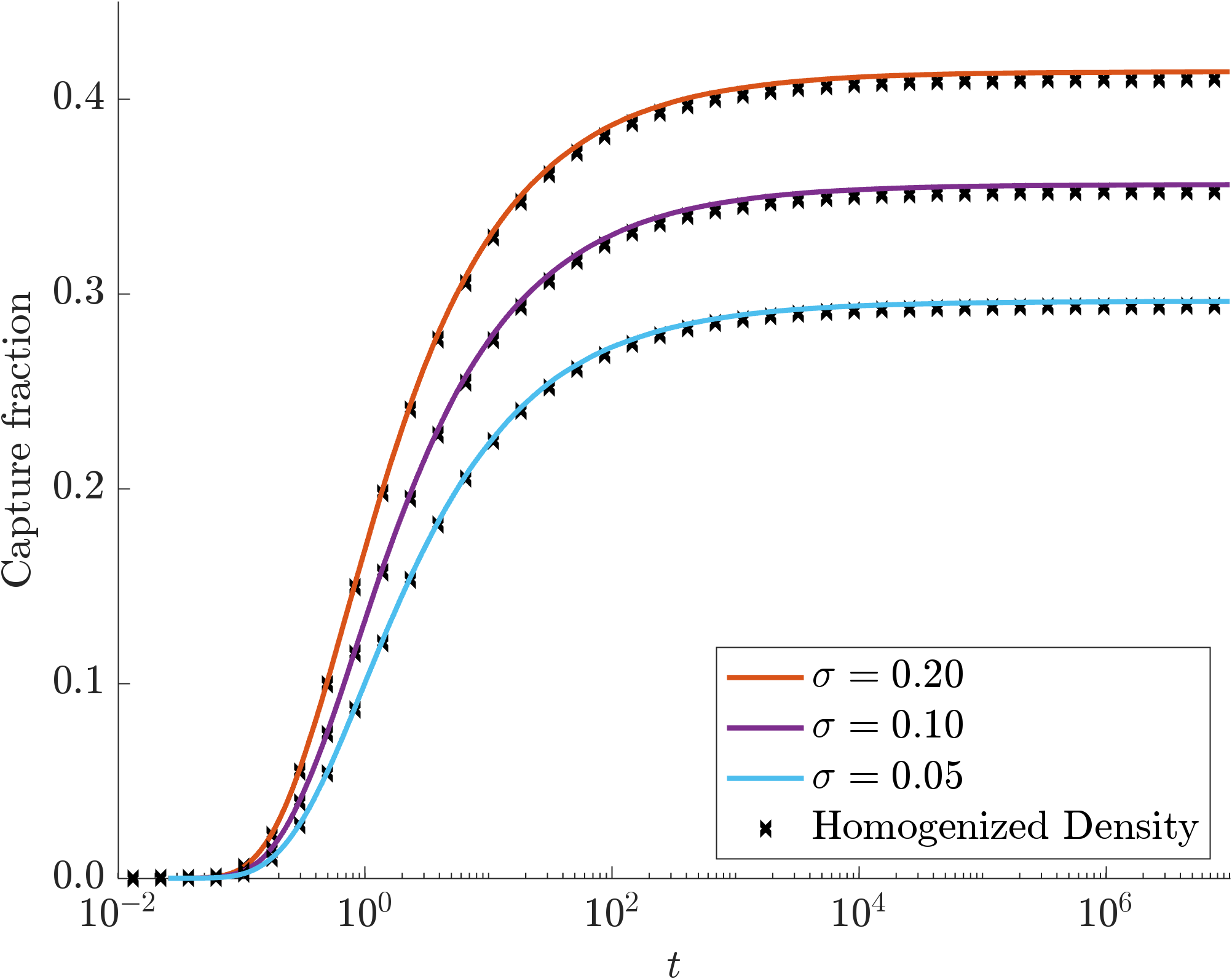}}
\caption{Comparison of homogenized theory with numerical results based on $M=10^6$ KMC trajectories initialized at $\bx_0 = (0,0,2.5)$. Panel (a): Schematic of the domain with $N=51$ pores with combined absorbing surface fraction of $\surfrac = 0.1$. Panel (b): Logarithmic histogram of the  arrival time distributions from KMC data and homogenized theory for $\surfrac = 0.1$ and diffusivities $D=10^{-2},10^{-1},10^{0}$.  Panel (c): Capture fraction over time predicted by homogenization theory and KMC for $\surfrac = 0.1$ and diffusivities $D=10^{-2},10^{-1},10^{0}$. Panel(d): Capture fraction over time predicted by homogenization theory and KMC for $D = 0.1$ and surface absorbing fraction $\surfrac=0.2, 0.1, 0.05.$ \label{fig:SphereDiff}}
\end{figure}

To validate both homogenization and the numerical method, we compare solutions of \eqref{eqn:IntroP} computed from the KMC algorithm with both the homogenized flux densities \eqref{eq:HomogSphereFlux} and capture fraction \eqref{eq:HomogSphereCDF}. In Fig.~\ref{fig:SphereDiff} we consider a sphere of unit radius with $N=51$ absorbing windows of common radius $a$ centered at the Fibonacci spiral points \cite{LWB2017} and perform two experiments based on $M=10^6$ KMC trajectories initialized at $\bx_0 = (0,0,2.5)$. Each of the circular pores in the mesh shown in Fig.~\ref{fig:SphereDiff_a} are represented by $10$ fixed points. The remaining mesh points are placed uniformly on the sphere and then repositioned by minimization of a repulsive discrete energy based on the reciprocal of pairwise distances. 

First we fix $D=1$ and vary the combined absorbing fraction $\surfrac= \frac{Na^2}{4}$ over values $\surfrac = 0.2,0.1,0.05$. Second, we fix the fraction $\surfrac=0.1$ and vary the diffusivity over values $D=10^{-2},10^{-1},10^{0}$. In each case we observe that the homogenized theory very closely matches full numerical simulations. We remark that the total probability of capture \eqref{eq:totalfractCapture} varies with the absorbing surface fraction $\surfrac$ while changes in the diffusivity control the timescale of equilibriation.  

\paragraph{Splitting probabilities} In the previous example, we used homogenization to describe the capture rate over the whole spherical surface. In this example, we consider the capture rate to individual receptors and compare with the splitting probabilities. These quantities $\{ \Q_k(\bx) \}_{k=1}^N$ describe the probability that a particle originating from $\bx$ is first captured at the receptor $\Omega_{k}$. They satisfy
\bsub\label{eqn:split}
\begin{gather}
 \Delta \Q_k = 0 , \quad \bx\in \mathbb{R}^3\setminus\Omega; \qquad \Q_k(\bx)\ \mbox{bounded as}\ |\bx|\to\infty;\\[4pt]
\Q_k  = \delta_{jk}, \quad \bx\in\partial\Omega_j, \quad j = 1,\ldots,N; \qquad  \nabla \Q_k  \cdot \hn = 0, \quad \bx\in\Gamr,
\end{gather}
\esub
We focus on the case where $\partial\Omega$ is the unit sphere and each receptor $\Omega_k$ is a circular patch of radius $a$ centered at point $\bx_k$ as described by \eqref{eqn:Pore_extent}. In \cite{LLM2020} it was shown that in the limit $a\to0$, a two term expansion for the solution of \eqref{eqn:split} is given by
\begin{equation}\label{eq:AsySplit}
\Q_k(\bx) = 4aG(\bx,\bx_k) + \frac{4a^2}{\pi} \Big[  \Big(\frac{3}{2} - \log (2a)\Big)G(\bx;\bx_k) - 4\pi \sum_{\substack{j=1\\ j\neq k}}^N G(\bx_j;\bx_k)G(\bx;\bx_j)\Big] + \littleoh(a^2),
\end{equation}
where $G(\bx,\bxi)$ is the Green's function of the Laplacian, exterior to the unit sphere. For $|\bxi|=1$, $G(\bx,\bxi)$ is given \cite{SurfaceGreen3D} by
\begin{equation}
G(\bx;\bxi) = \frac{1}{2\pi} \left[ \frac{1}{|\bx-\bxi|} - \frac{1}{2} \log \left( \frac{1-\bx\cdot\bxi + |\bx-\bxi|}{|\bx| - \bx\cdot\bxi} \right) \right].
\end{equation}

\begin{figure}[htbp]
\centering
\subfigure[Cumulative arrival fraction $m_k(t;\bx_0)/M$ and the splitting probabilities $\Q_k(\bx_0)$ given in \eqref{eq:AsySplit}.]{\includegraphics[width = 0.45\textwidth]{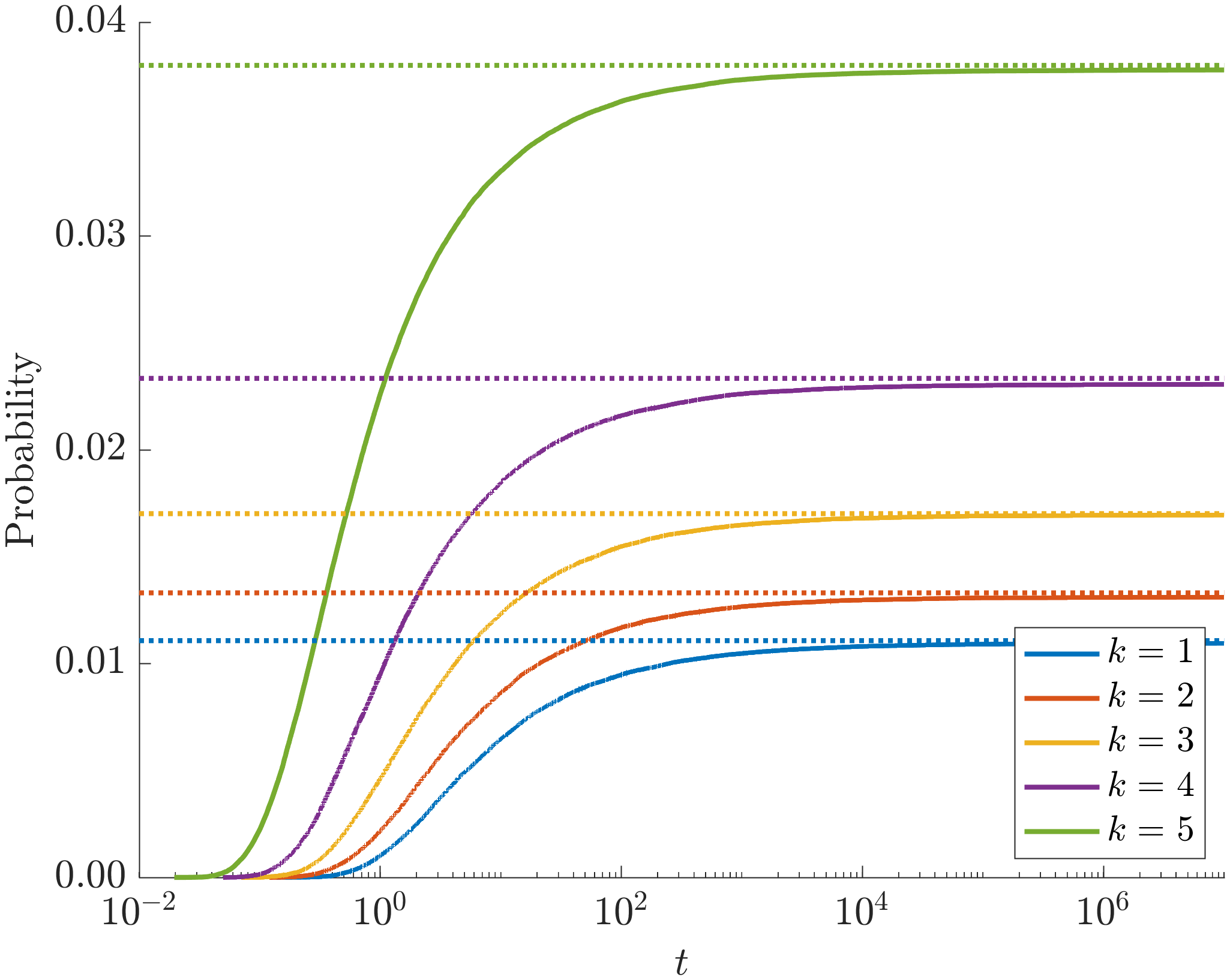} \label{fig:sphere_a}}
\qquad
\subfigure[Domain schematic and source location.]{ \includegraphics[width = 0.325\textwidth]{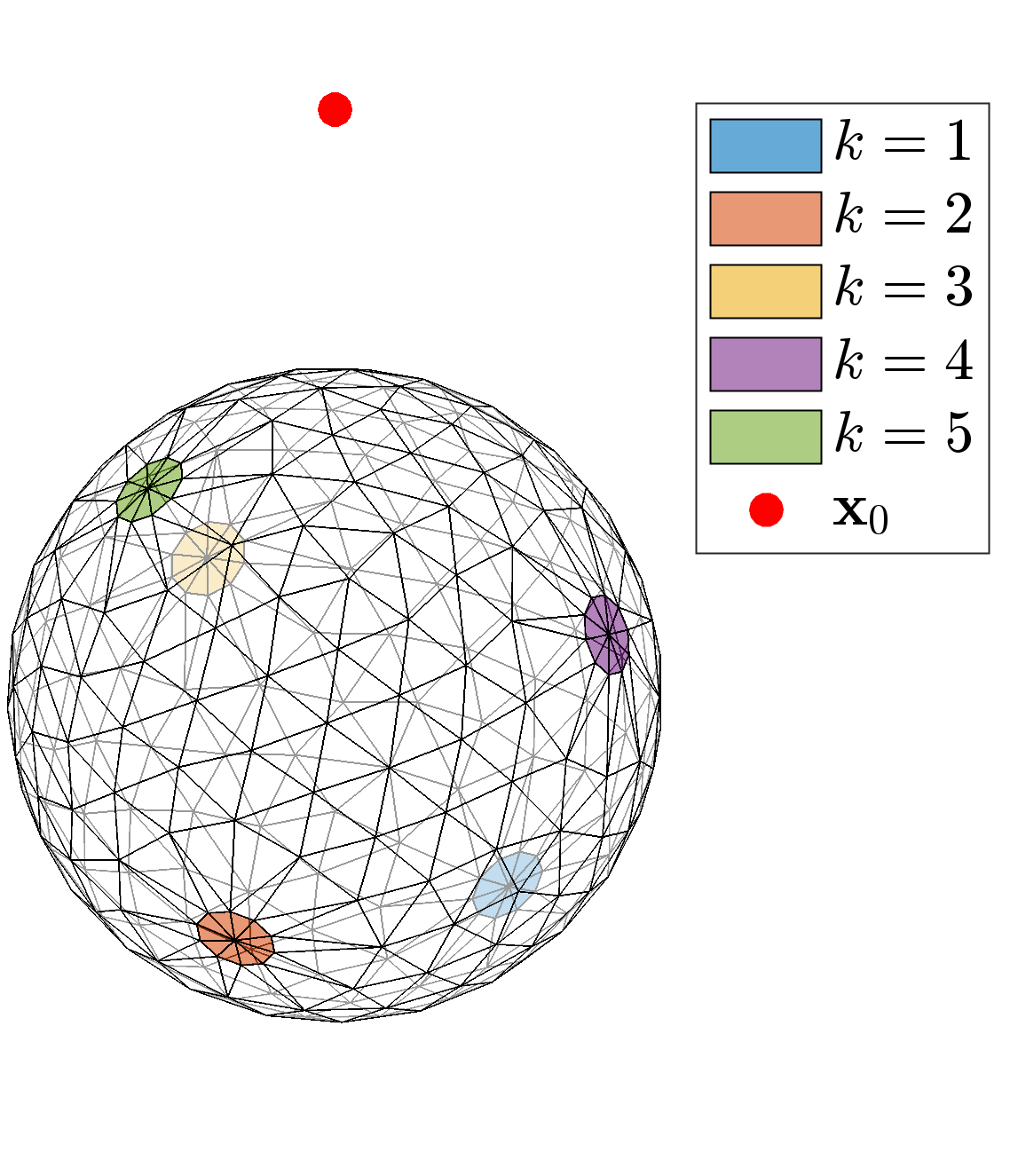} \label{fig:sphere_b}}
\caption{Convergence of individual KMC pore data to the splitting probabilities. Panel (a): The empirical CDF $m_k(t;\bx_0)/M$ for each of the $N=5$ pores. The limiting values agree with the asymptotic splitting probabilities \eqref{eq:AsySplit}. Panel (b): Schematic of the domain with initial condition $\bx_0=(0,0,2)$ and pores highlighted.\label{fig:sphere} }
\end{figure}

The splitting probabilities are equilibrium quantities, i.e.~they are fully determined when all particles originating from $\bx_0\in \Omega$ have arrived at a receptor or escaped to infinity. To describe the dynamic approach to these steady quantities using KMC, we consider $M$ particles originating at $\bx_0$ and calculate the number of particles $m_k(t;\bx_0)$ which have arrived at the $k^{th}$ receptor by time $t$. The fraction of particles absorbed at each receptor converge to the splitting probabilities, specifically,
\begin{equation}\label{eq:limitingSplit}
    \Q_k(\bx_0) = \lim_{t\to\infty} \frac{m_k(t;\bx_0)}{M}.
\end{equation}

As an example to demonstrate the convergence of the cumulative fluxes to the splitting probabilities, we consider a simple scenario of $N=5$ receptors centered at Fibonacci spiral points \cite{LWB2017} with common radius $a$ and absorbing surface fraction $\surfrac = \frac{Na^2}{4} = 0.02$. We initiate $M=10^6$ trajectories of diffusivity $D=2$ from the initial point $\bx_0 = (0,0,2)$ and calculate the fraction of particles captures at each receptor. In Fig.~\ref{fig:sphere} we show the domain and observe the limiting behavior described in \eqref{eq:limitingSplit}. As usual with exterior diffusion problems, the timescale of equilibriation is long. 

\begin{figure}[htbp]
    \centering
    \includegraphics[width = 0.45\textwidth]{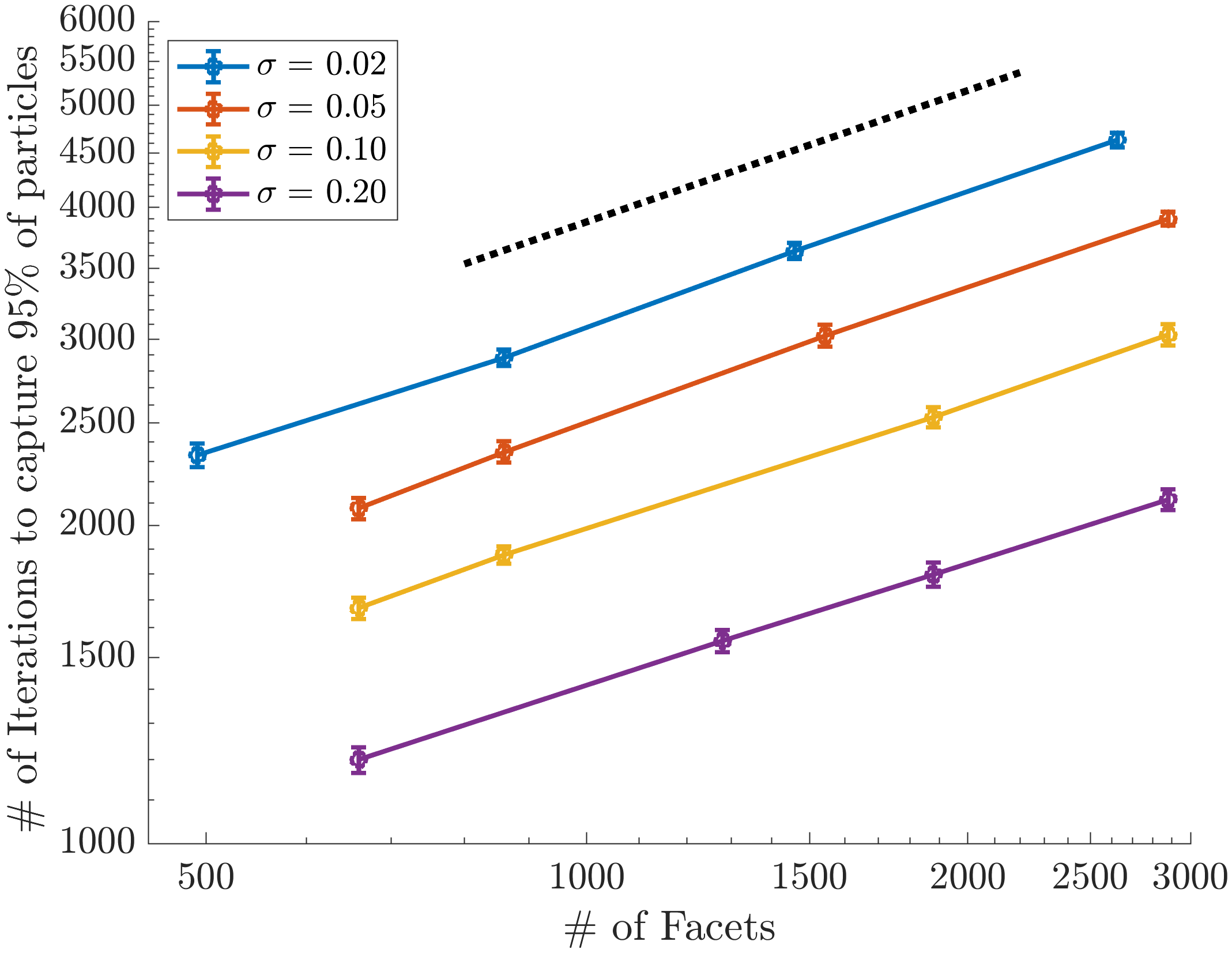}
    \caption{\al{Scaling of computational effort (number of iterations) with surface refinement (number of facets). Each data point represents the mean of $40$ repetitions with error bars indicating standard deviations. Surfaces with smaller absorbing fractions ($\sigma$) and more detailed refinements require increased iterations to capture particles. A best fit line (dashed) of slope $p=0.412$ is added for comparison which suggests the computational effort scales roughly as the square root ($p=1/2$) of the number of facets. \label{fig:scaling}}}
\end{figure}
\al{\paragraph{Computational scaling} Here we evaluate the computational scaling of the KMC algorithm in terms of the surface absorbing fraction and level of surface refinement. The KMC algorithm exactly simulates jumps on reflecting portions of the surface. A particle on a reflecting facet escapes roughly 20 \% of the time and usually impacts an adjacent facet, effectively taking a random walk on the sphere. Hence highly refined triangulations will require many more iterations to reach an absorbing site.} 

\al{To explore these effects on computational time, we simulated arrivals to a sphere with $M = 10^5$ particles initiated at $\bx_0 = (0,0,2.5)$. For $N=21$ identical circular absorbers centered at the Fibonacci spiral points, we varied the surface absorbing fraction $\sigma$ over values $\{0.01,0.05, 0.1, 0.2\}$ for various levels of surface refinement. We calculated the number of iterations (individual projection steps) required to capture $95\%$ of the initial particles. We observe that more iterations are required when either the surface is highly refined or has a small absorbing coverage fraction $\sigma$. In Fig.~\ref{fig:scaling} we observe that computational effort scales approximately with the square root of the number of facets (best fit slope $p=0.412$). The linear dimensions of a facet (such as the perimeter or the circumradius) scale like the reciprocal of the square root of the number of facets.  Berg and Purcell \cite[Eqn. (10)]{bergp} in their classic paper \emph{Physics of Chemoreception} (1977) argued that the number of jumps of length at least $\ell$ for a random walk on the surface of sphere leading to absorption scales as $1/\ell$ in the limit of small $\ell$; equating these two processes and lengthscales gives a heuristic explanation for the observed scaling.}


\subsection{Arrival distribution to a family of ellipsoidal geometries}\label{sec:Ellipsoid}

\begin{figure}[htbp]
    \centering
    \subfigure[Family of ellipsoidal geometries with two circular pores of unit radius centered around the vertical axis on planes $z\pm1$. The top pore is aligned towards the source (green dot) positioned at $\bx_0 = (0,0,z_0)$ for $z_0>1$.]{\includegraphics[width = 0.75\linewidth]{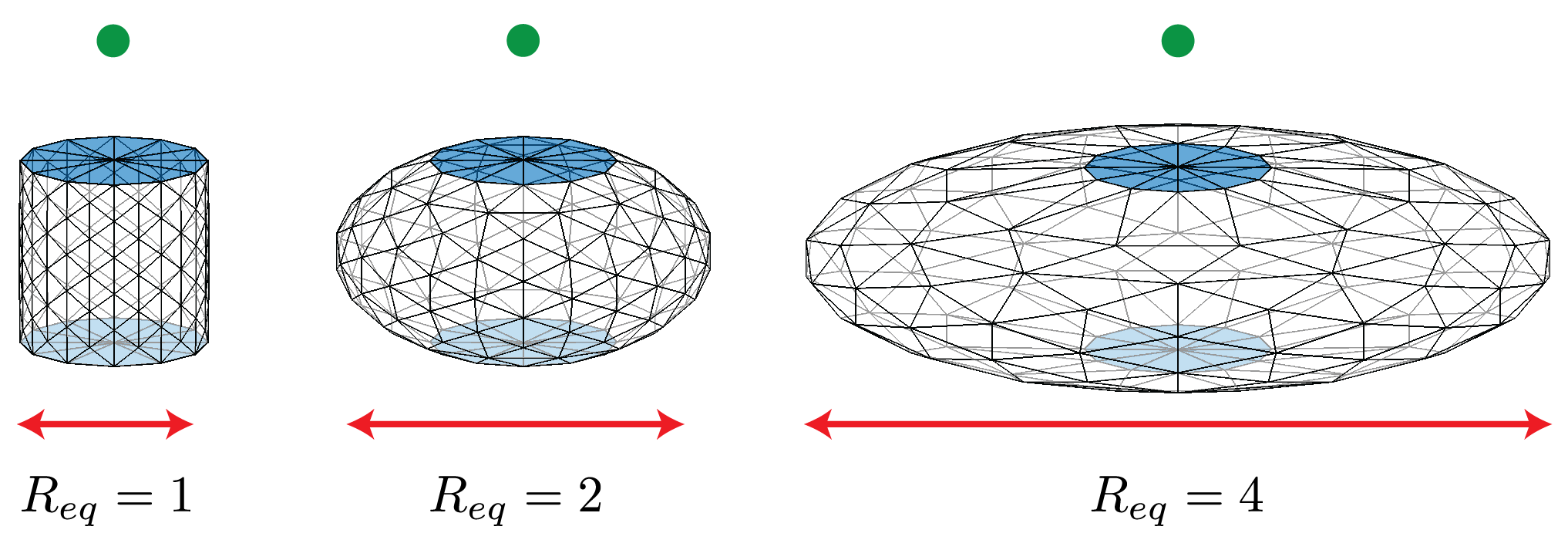}\label{fig:ellipsoid_a}}\\
    \subfigure[Splitting probability between pores.]{\includegraphics[width=0.315\linewidth]{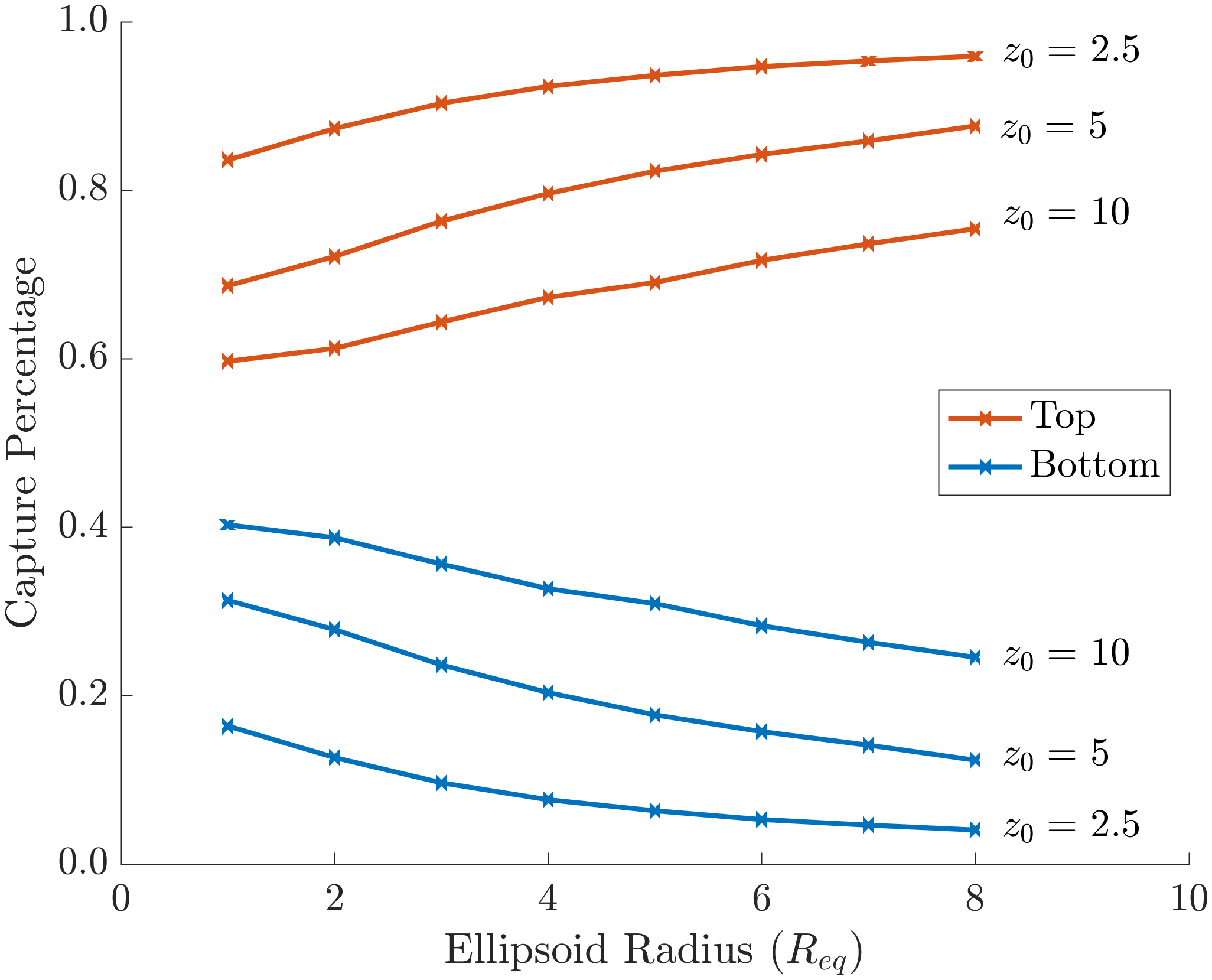}\label{fig:ellipsoid_b}}\quad
    \subfigure[Capture fraction dynamics for initial point $\bx_0 = (0,0,10)$.]{\includegraphics[width=0.315\linewidth]{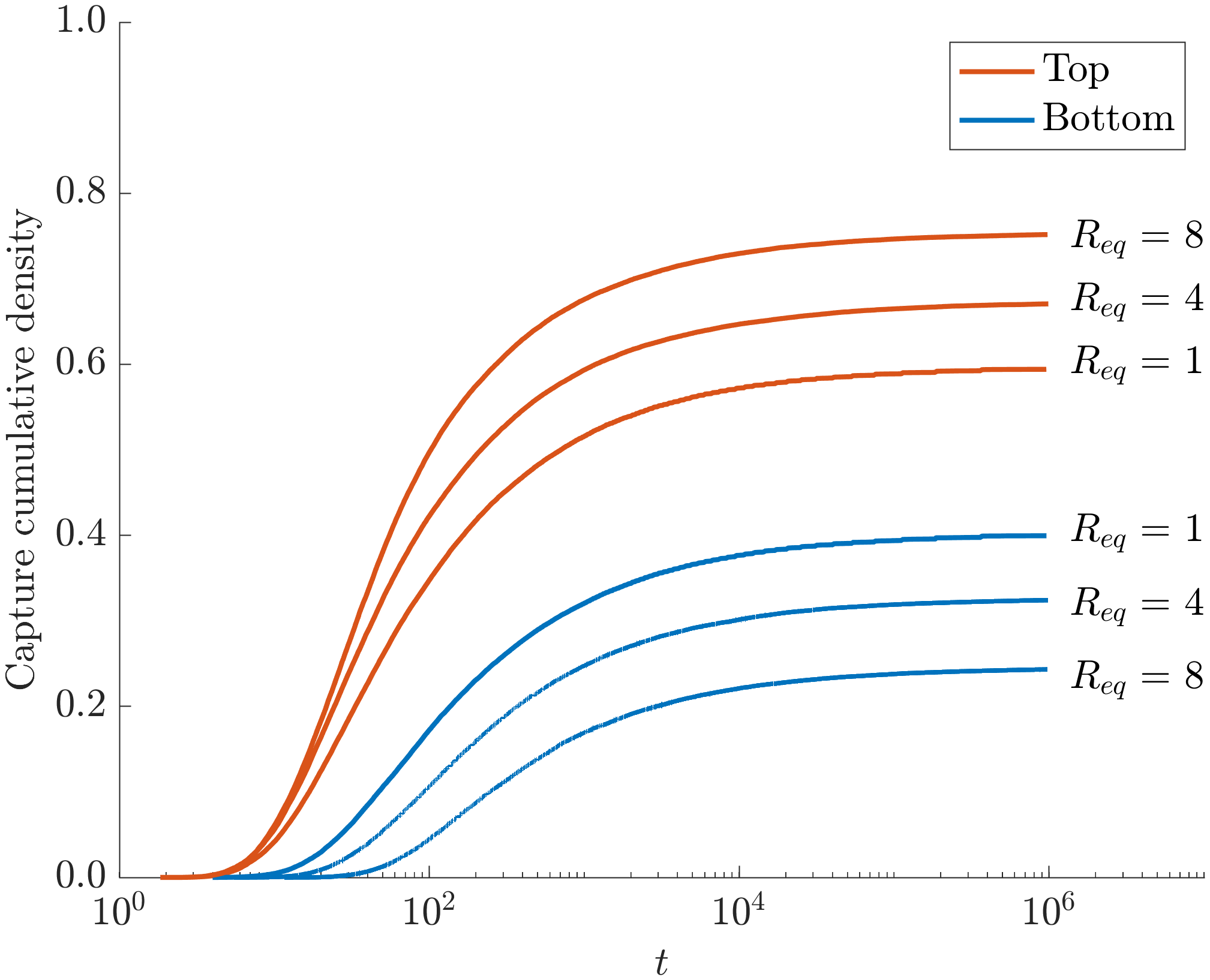}\label{fig:ellipsoid_c}}\quad
    \subfigure[\ajb{Differential flux densities} and CDF.]{\includegraphics[width=0.315\linewidth]{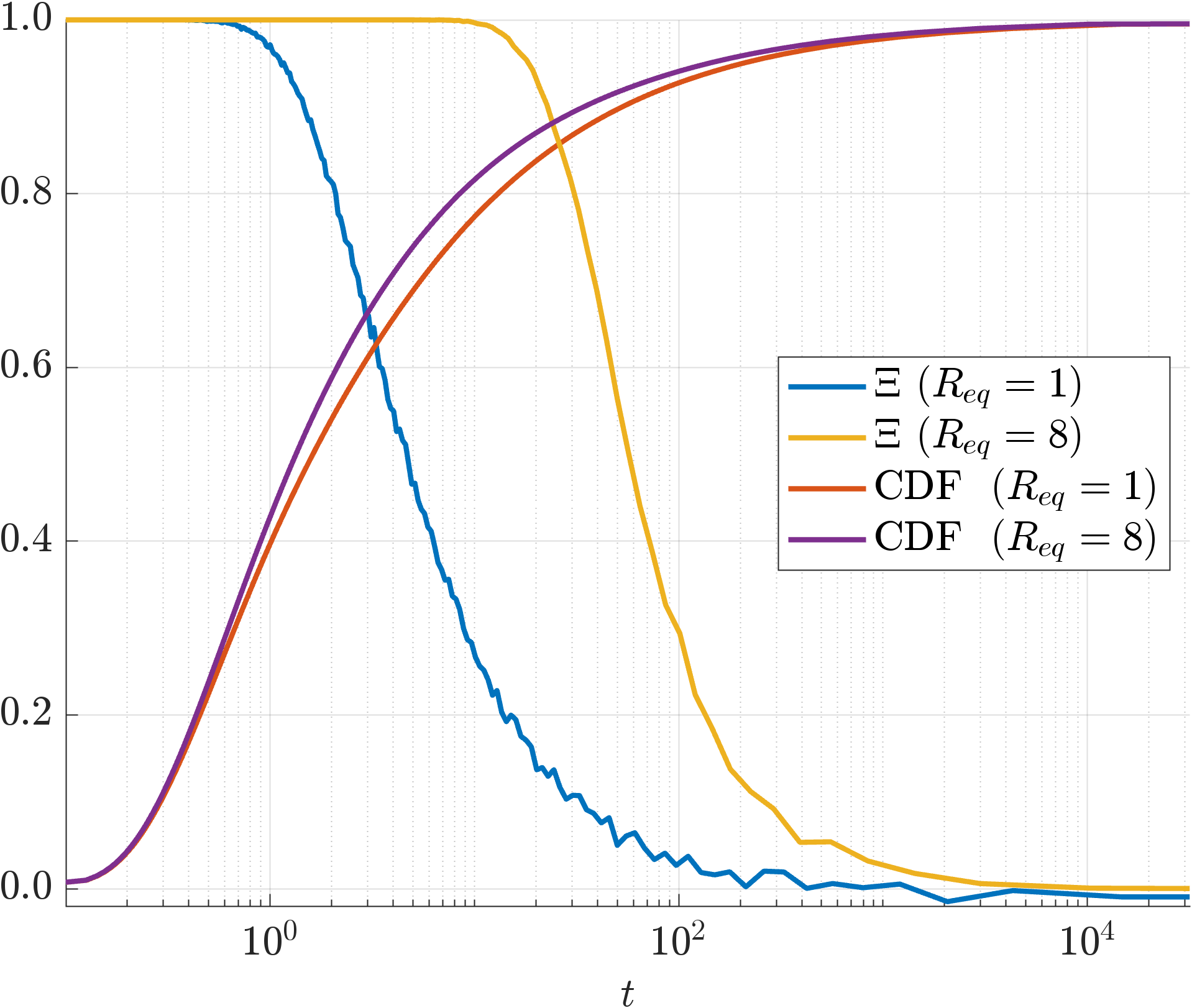}\label{fig:ellipsoid_d}}
    \caption{The role of geometry in directional sensing of source at $\bx_0 = (0,0,z_0)$. Panel (a): Schematic of non spherical domains where $R_{eq}\geq1$ parameterizes a shielding effect in which the geometry reduces the capture rate to the bottom pore and increases it to the top. Panel (b): Capture probabilities for a range of $R_{eq}$ for initial elevations $z_0$. Smaller values of $z_0$ and larger values of $R_{eq}$ are associated with increased capture rates at the top pore. Panel (c): The dynamic approach of \ajb{cumulative} fluxes to the capture probabilities for initial elevation $z_0 = 10$ and several $R_{eq}$ values. Panel (d): The scaled differential flux \ajb{densities} between pores \eqref{eq:bigXi} and the combined capture rate for $z_0 = 2.5$ and values $R_{eq}=1,8$.}
    \label{fig:ellipsoid}
\end{figure}

In this section we consider a family of ellipsoids with two circular pores located at the north and south poles. The pores are fixed with unit radius in the planes $z=\pm 1$. The surface joining them forms a skirt of radius $R_{eq}$ that varies between a cylinder ($R_{eq}=1$) and an oblate spheroid many times wider than the pore ($R_{eq}>1$).  We consider a source on the polar axis above the sphere as shown in Fig.~\ref{fig:ellipsoid_a}.

Differential receptor activity over distal sections of a cellular surface (ratiometric sensing) has been suggested as a mechanism for the inference of chemical cues \cite{Lew2019,Ismael2016,Lakhani2017,Bumsoo2019}. Cell morphology is frequently non-spherically, a fact which has been observed to modulate this process \cite{nakamura2024gradient,kaiyrbekov2024does}. To shed light on how cellular geometry might influence response to external cues, we calculate the splitting probabilities of arrivals at the top and bottom pores over the range $R_{eq} \in[1,8]$ as shown in Fig.~\ref{fig:ellipsoid_b}. We observe that both proximity to the source ($z_0$) and larger $R_{eq}$ give a stronger signal to the pore aligned with the source (Top). In Fig.~\ref{fig:ellipsoid_c} we focus on the particular source location $\bx_0 = (0,0,10)$ and determine the \ajb{cumulative} dynamic fluxes to each pore, confirming the convergence to the \ajb{splitting} probabilities displayed in Fig.~\ref{fig:ellipsoid_b}. As seen in previous examples, equilibration occurs over a long timescale that typically exceed those observed in biological examples of directional sensing \cite{servant2000}. To explore the possibility for directional inference before steady state, we consider the normalized differential flux \ajb{density} between the top ($\J_{\textrm{Top}}(t)$) and bottom ($\J_{\textrm{Bottom}}(t)$) pores,
\begin{equation}\label{eq:bigXi}
\Xi(t) = \frac{\J_{\textrm{Top}}(t) - \J_{\textrm{Bottom}}(t)}{\J_{\textrm{Top}}(t) + \J_{\textrm{Bottom}}(t)}.
\end{equation}
The quantity $\Xi(t)$ gives a dynamic measure of strong ($\Xi\approx 1$) or weak ($\Xi\approx 0$) directional information. Strong directional information is expected at shorter times, before particles have had a chance to \emph{thermalize}, i.e. explore a sufficient volume of parameter space to lose information about where they started.  In Fig.~\ref{fig:ellipsoid_d}, we plot $\Xi(t)$ for source location $\bx_0 = (0,0,2.5)$ and geometries $R_{eq} = 1$ and $R_{eq} = 8$ together with their respective cumulative capture fractions. We observe that $\Xi(t)\approx 1$ over a significantly extended timescale in the case $R_{eq}=8$ compared to $R_{eq}=1$ while the cumulative capture fraction is largely unchanged. This demonstrates a plausible new mechanism for improving the strength and duration of the directional signal. Specifically as the oblateness of the ellipsoid increases it both shields the lower pore, increasing the differential between the splitting probabilities, and lengthens the time for thermalization. 

The symmetry of the geometry gives an easy visualization of when thermalization occurs;
when a particles strikes the equatorial plane ($z=0$), it becomes equally likely to be captured by either pore.  The more of the equatorial plane that is shielded by the skirt of the oblate spheroid, the longer thermalization of the ensemble of particles takes.

\section{Discussion}
In this work we have presented and validated a numerical method for the simulation of three dimensional random walks to convex surfaces with absorbing and reflecting portions. In a reduced scenario of small non overlapping absorbers on the plane, we also derived a new matched asymptotic expansion \eqref{eq:MainPlanarResult} for the capture rate to individual absorbers.

In the case of convex three dimensional geometries we have validated the KMC method against several existing steady state results for simple geometries such as the cube (Sec.~\ref{sec:resultsCube}) and sphere (Sec.~\ref{sec:sphere}). Overall, our computations confirm that homogenization is a very effective method for describing dynamic fluxes to a target set, with errors largely confined to the $t\to 0^{+}$ regime. This is precisely the segment of the distribution that characterizes extreme statistics of Brownian motion \cite{Lawley2020a,Lawley2020} and hence care must be used when applying homogenization to such scenarios.

\al{The asymptotic results presented here for the plane suggest an important and achievable goal for future studies is to obtain equivalent asymptotic approximations for the dynamic fluxes to receptors arranged on the sphere and other general three dimensional geometries \cite{Tzou2023,Cheviakov2015}. This can be achieved through application of the Laplace transform coupled with matched asymptotic analysis informed by detailed local behavior of the exterior Helmholtz Green's function.}

Homogenization provides estimates for the capture rate across the whole surface by averaging out local variations. However, the sensing of chemical cues often necessitates a comparison between fluxes across individual pores. The method developed here can rapidly and accurately calculate the dynamics of these signals and we have demonstrated their applicability to directional sensing problems in spherical and non-spherical scenarios. In a family of ellipsoidal domains of varying eccentricity, we demonstrate a potential mechanism for promoting directional sensing through geometric shielding.

While we believe the machinery developed here provides a firm foundation for further investigation of capture problems in exterior domains, it is useful to report on several possible avenues for improving these methods. An advantage of Monte Carlo methods that we have not fully exploited is that particle trajectories are completely independent and can therefore fully leverage massively parallel computer configurations. In scenarios where the geometry requires representation by a very fine triangulation, it is essential to optimize the calculation of the signed distance function which can be accomplished with a tree based search \cite{CLNQ22,Wang-2022-dualocnn}.

\al{The KMC method developed here exactly resolves many of the challenging aspects of simulating diffusion, including escape, long time integration and surface flux singularities. However, an equivalent PDE method would be highly desirable for its high accuracy and ability to simultaneously resolve the solution in all space. Avenues for the development of such a method include boundary integral methods which are well suited to solving elliptic PDEs in exterior geometries \cite{BL2018,cherry2024boundaryintegralmethodsparticle,KAYE2020}.}

Another set of questions that arises naturally is that of extreme statistics an example of which would be resolving the fat tails that occur at large times in these simulations. The problem is that even with $M=10^8$ particles the representation in these tails is extremely sparse. A modification that would allow the exploration of these extremely unlikely events is Markov Chain Monte Carlo methods whereby particles with long survival times are branched (in a weighted fashion) \al{to} allow better resolution of these extremely rare events \cite{Driscoll2007}. 


\section*{Acknowledgments}
A.~E.~Lindsay was supported by NSF grant DMS-2052636.  
A.~J.~Bernoff was supported by Simons Foundation grant 317319. \al{We gratefully acknowledge the careful reading and helpful comments from the anonymous reviewers.}

\appendix

\section{First Arrival for a Point Source External to a Sphere}\label{sec:spherearrival}
In this appendix we will solve the diffusion equation 
(\ref{eqn:IntroP_a}) for the domain exterior to an absorbing unit sphere, $\cS$, at the origin with a point source on the polar axis at $\bx_0=(0,0,R)$ with $R>1$.  As the geometry is  axisymmetric and the absorbing surface is spherically symmetric we use spherical coordinates for which the resulting problem is separable. Let $p(x,y,z,t) = p(r,\theta,t)$ where $r$ is the distance to the origin and $\theta$ is the polar angle (cf.~Fig.~\ref{fig:app:angle_dist}). The axisymmetric diffusion equation in spherical coordinates is
\bsub\label{eq:SphereExt}
\begin{align}
\frac{\partial p}{\partial t} &=\D \left [ \frac{1}{r^2} \frac{\partial }{\partial r}\Big( r^2 \frac{\partial p}{\partial r} \Big) + \frac{1}{r^2 \sin\theta} \frac{\partial }{\partial \theta}\Big( \sin \theta \frac{\partial p}{\partial \theta} \Big)\right ], \qquad r>1, \quad \theta\in\al{ [0,\pi)};\\[5pt]
p(1,\theta,t) &= 0, \qquad \theta \in \al{ [0,\pi)}; \\[5pt]
p(r,\theta,0) &= \frac{1}{\tp R^2\sin\theta} \delta(r-R)\delta(\theta), \qquad r > 1,\quad \theta \in \al{[0,\pi)};\\[5pt]
\al{p(r,\theta,t)} & \al{\text{ smooth as } \theta \to 0^{+}\ \mbox{for} \ r\neq R.}
\end{align}
\esub
\ajb{The last two boundary conditions here deal with the coordinate singularity at $\theta=0$; in Cartesian coordinates the initial density is $p(x,y,z,0)= \delta(\bx-\bx_0)$, that is a unit mass at $\bx=\bx_0$.}

The solution to \eqref{eq:SphereExt} is sought via the Laplace transform, 
$$\al{\hat{p}(r,\theta;s)} = \int_0^{\infty} e^{-st} p(r,\theta,t) \, dt, $$ 
where $\al{\hat{p}(r,\theta;s)}$ solves the modified Helmholtz equation,
\bsub\label{eq:SphericalLaplace}
\begin{align}
\label{eq:SphericalLaplace_A} \D \left [ \frac{1}{r^2} \frac{\partial }{\partial r}\Big( r^2 \frac{\partial \hat{p}}{\partial r} \Big) + \frac{1}{r^2 \sin\theta} \frac{\partial }{\partial \theta}\Big( \sin \theta \frac{\partial \hat{p}}{\partial \theta} \Big) \right ] - s \hat{p} &= \frac{-1}{\tp R^2\sin\theta} \delta(r-R)\delta(\theta)\, \quad r>1, \quad \theta \in \al{ [0,\pi)}; \\[5pt]
\label{eq:SphericalLaplace_B} \al{\hat{p}(1,\theta;s)} & = 0, \quad \theta \in\al{ [0,\pi)};
\end{align}
\esub
The separable general solution of \eqref{eq:SphericalLaplace}, which is finite as $r\to \infty$, satisfies $\al{\hat{p}(1,\theta;s)=0}$, \al{continuous at $r=R$, and has no singularity at $\theta=0$ for $r\neq R$}, is given by the series
\begin{equation}
\al{\hat{p}(r,\theta;s)} = \left\{
\begin{array}{ll} 
\ds\sum_{n=0}^{\infty} A_n \left[ i_n(\sD r) -  \frac{i_n(\sD)}{k_n(\sD)} k_n(\sD r) \right]  P_n(\cos \theta),  & r < R;\\[10pt]
\ds\sum_{n=0}^{\infty} A_n \left[ \frac{i_n(\sD R)}{ k_n(\sD R)} -  \frac{i_n(\sD)}{k_n( \sD )} \right] k_n( \sD r) P_n(\cos \theta),  & r > R,\\
\end{array}
\right.
\label{eq:LegendreExpansion}
\end{equation}
where $\sD= \sqrt{s/D}$. The modified spherical Bessel functions $i_n(x)$ and $k_n(x)$ are defined \al{\cite[{\S 10.2}]{abramowitz1965handbook}} as\footnote{\ajb{
The reader is cautioned that the definition and normalization of modified spherical Bessel functions varies between sources. For reference, we note that this definition yields $i_0(x)=\frac{\sinh x}{x}$ and $k_0(x)=\frac{\pi}{2} \frac{e^{-x}}{x}$.
}}
\begin{equation}
i_n(x) = \sqrt{\frac{\pi}{2 x}} I_{n+\frac12}(x), \qquad k_n(x) = \sqrt{\frac{\pi}{2 x}} K_{n+\frac12}(x)\al{,}
\end{equation}
for modified Bessel functions $I_n(x)$ and $K_n(x)$. The functions $P_n(x)$ are the Legendre polynomials normalized \al{by} $P_n(1)=1$ which satisfy the orthogonality condition
\begin{equation}\label{eq:Porthog}
\int_{\theta=0}^\pi P_n(\cos\theta) P_m(\cos \theta ) \sin \theta \, d \theta =
\int_{x=-1}^1 P_n(x) P_m(x)\, dx 
=  \frac{2}{2n+1} \delta_{mn},
\end{equation}
where $\delta_{mn}$ is the Kronecker delta function. The constants $A_n$ are fixed by incorporating the Dirac source on the right hand side of \eqref{eq:SphericalLaplace_A}. We obtain a jump condition by substituting the expansion \eqref{eq:LegendreExpansion} into the governing equation
\eqref{eq:SphericalLaplace_A}, multiplying by $P_n(\cos \theta)$ and integrating over the sphere, then by integrating a small interval around the $\delta$-function,
$r\in(R-\varepsilon,R+\varepsilon)$ and passing to the limit $\varepsilon \to 0$ yielding
\begin{equation} \al{
\tp D \left ( \frac{2}{2n+1} \right ) 
A_n \sD R^2 \left[    \frac{i_n(\sD R)}{k_n(\sD R)} k_n'(\sD R) - i'_n( \sD R ) \right] = - 1.}
\end{equation}
The Wronskian identity yields
\[
k_n(x)\, i_n'(x) - k_n'(x)\, i_n(x)  = \al{\frac{\pi}{2}} \frac{1}{x^2},
\]
simplifies the constants $A_n$ to
\[
\al{A_n = \frac{1}{ \pi^2 D}  \frac{(2n+1)}{2} \, \sD \, k_n(\sD R).}
\]
The Laplace transform of flux density through the spherical surface is then
\begin{equation}
D \al{\hat{p}_r(1,\theta;s)} = \al{\alpha}D \ \sum_{n=0}^{\infty}  A_n \left[i'_n( \sD  ) -  \frac{i_n( \sD )}{k_n( \sD )} k_n'( \sD ) \right]  P_n(\cos \theta) = \rtp \sum_{n=0}^{\infty} \frac{(2n+1)}{2} \frac{k_n( \sD R)}{k_n( \sD )} P_n(\cos \theta).
\end{equation}
The Laplace transform of the total flux is now given by integrating over the surface of the sphere, remembering that $\sD = \sqrt{s/D}$, 
\begin{align}
    \nonumber \hat{\JS}(s) &= \iint_{\cS} D \al{\hat{p}_r(1,\theta;s)} \, dS = 2 \pi \int_0^{\pi} \left [ \rtp  \sum_{n=0}^{\infty} \frac{(2n+1)}{2} \frac{k_n( \sD R)}{k_n( \sD )} P_n(\cos \theta)\right ] \, \sin \theta \, d \theta \\
\label{eqn:AppenFluxLT} &=  \frac{k_0(R\sqrt{s/D})}{k_0(\sqrt{s/D})} = \frac{e^{-(R-1)\sqrt{s/D}}}{R},
\end{align}
where we have \al{applied the orthogonality properties of $P_n$ stated in \eqref{eq:Porthog} and applied the identity \cite[\S 10.2.17]{abramowitz1965handbook} $$k_0(x) = \sqrt{\frac{\pi}{2x}}K_{\frac12}(x) = \frac{\pi}{2} \frac{e^{-x}}{x}.$$}
The inverse Laplace transform of \eqref{eqn:AppenFluxLT} is the flux \ajb{density} (the PDF of arrival times) into the sphere, $\cS$, and given by \al{\cite[{\S 29.3.82}]{abramowitz1965handbook}}
\begin{equation}
\JS(t)= \frac{R-1}{2R\sqrt{\pi D}} \cdot t^{-3/2}  \exp\left[ \frac{-(R-1)^2}{4Dt}\right].
\end{equation}
The \al{cumulative distribution function} (CDF) of this distribution can then be calculated as
\begin{equation}
P_T(t) = \int_0^t \JS(\tilt) \, d\tilt = \frac{1}{R} \mbox{erfc} \left[\frac{R-1}{2\sqrt{D t}} \right].
\label{eq:CDFsphere}
\end{equation}
We remark that 
$$\int_0^{\infty} \JS(\tilt) \, d\tilt = 
\lim_{t \to \infty} P_T(t) =
R^{-1},$$
so that the probability of capture is not unity, but inversely proportional to the initial distance to the sphere.

The KMC propagator described in Sec.~\ref{subsubsec:KMCReinsert} first determines if a particle escapes to infinity or impact the sphere at some $t_*$. If the particle impacts the sphere one needs to determine the probability density of where it strikes conditioned on $t_*$. This distribution is uniform in the azimuthal angle $\phi \in [0,2\pi)$. \ajb{The probability distribution of the polar angle $\theta$, condition on arrival at time $t_*$, and with the azimuthal dependence integrated out} is 
\begin{align}
\nonumber \rho_{\theta}(\theta;t_*) &= \al{\int_{\phi = 0}^{2\pi}  \frac{ D p_r(1,\theta;t_{\ast})}{\J(t_{\ast})} \,d\phi =
\frac{ \tp  Dp_r(1,\theta\al{;}t_*)}{\JS(t_*) }} \\[5pt]
&= \frac 12 
+\frac{1}{{\JS(t_*)}}\sum_{n=1}^{\infty} \frac{(2n+1)}{2} \chi_n(t_*) P_n(\cos\theta)\, ,\label{eq:JointSp}
\end{align}
where $\chi_n(t)$ are given by their inverse Laplace transform
\begin{equation}\label{eq:chiN}
\chi_n(t) = \frac{1}{2\pi i }\int_{s=\gamma -i\infty}^{\gamma +i\infty} \hat{\chi}_n(s) e^{st}\, \, ds, \qquad \hat{\chi}_n(s) = \frac{k_n(\alpha R)}{k_n(\alpha)}, \qquad \alpha = \sqrt{\frac{s}{D}}.
\end{equation}
\al{In \eqref{eq:chiN} $\mathcal{B}=\{ \gamma + i \eta \ | \ \eta \in\mathbb{R}\}$ is the Bromwich contour chosen so that the vertical line $\gamma=\mbox{Re}(\mathcal{B})$} lies to the right of any singularities in the integrand, in this case a branch cut along the negative real axis, and decaying for a ray at angle $ \textrm{Arg} (s) \in \left ( -\frac{\pi}{2}, \frac{\pi}{2} \right ) $.  An exact real integral can be obtained by deforming the contour to be a hairpin around the branch cut on the negative real axis \cite{CLNQ22,Lindsay2023b}. This yields that  
 \begin{equation}\label{eqn:chi_n}
\chi_n(t) = \frac{2D}{\pi } \int_{w=0}^\infty 
 \left [\frac{ j_n(w) y_n(wR) -y_n(w) j_n(wR) } {[j_n(w)]^2+[y_n(w)]^2 } \right ]
  w e^{-w^2 Dt} \, dw,
\end{equation}
\al{where $j_n(w)$ and $y_n(w)$ are spherical Bessel functions of the first and second kinds respectively. The scaling of $\chi_n(t)$ as $t\to\infty$ can be derived by considering that the main contribution arises when $w^2Dt = \mathcal{O}(1)$, or $w\ll1$. We can therefore estimate (see also \cite[Eq.~(4.5)]{CLNQ22}) that 
\begin{equation} \label{chinlarget}
\chi_n(t) \sim \frac{D R^n ( 1 - R^{-1-2n})}{2^{2n+1}\Gamma(n+ \frac12) }(Dt)^{-n-\frac32} \quad \mbox{as} \quad t\to\infty.
\end{equation}
}
In practice, we obtain $\chi_n(t)$ by inverting the Laplace transform \eqref{eq:chiN}, rather than evaluating the real integral \eqref{eqn:chi_n}. A simple and effective method is to deform the contour $\mbox{Re}(s) = \gamma$ into the left hand plane and apply numerical integration based on $M$ quadrature points. A highly optimized choice, called the Talbot contour \cite{TW2014,Abate2006}, is given by 
\begin{equation}\label{eqn:talbot}
  s(\varphi) = \frac{M}{t}\left(
    -0.6122 + 0.5017\varphi\cot(0.6407\varphi) + 0.2645i\varphi\right),
  \qquad \varphi \in (-\pi,\pi).
\end{equation}
The integral \eqref{eq:chiN} is then evaluated along \eqref{eqn:talbot} using the midpoint rule. \al{The convergence rate of this quadrature has been theoretically established at $\mathcal{O}(e^{-1.358M}$) and in practice machine precision is achieved for around $M= 24$ quadrature points \cite{Weideman2015}.}

Having the CDF for the arrival angle, $P_\Theta (\theta_* ;t_*) $, will allow us to eventually sample the arrival point for a KMC particle. Integrating \eqref{eq:JointSp} yields
\begin{align}
\nonumber 
P_\Theta (\theta_* ;t_*) 
& = \int_0^{\theta_*} \rho_{\theta}(\theta;t_*)  \sin\theta\,d\theta \, ,\\
\nonumber
&=\int_0^{\theta_*} 
\left[\frac 12 
+\frac{1}{{\JS(t_*)}}\sum_{n=1}^{\infty} \frac{(2n+1)}{2} \chi_n(t_*) P_n(\cos\theta)
\right ] \sin\theta \, d\theta \, ,\\
\nonumber & = \frac{1}{2} [1- \cos \theta_*] + 
\frac{1}{\JS(t_*)} \sum_{n=1}^{\infty} \frac{(2n+1)}{2} \chi_n(t_*) 
\left [ \int_{\cos\theta_*}^1 P_n(x) \,dx \right ] , \\
&= \frac{1}{2} [1- \cos \theta_*]+ 
\frac{1}{2}\sum_{n=1}^{\infty} \frac{\chi_n(t_*)}{\JS(t_*)}  [ P_{n-1}(\cos\theta_*) - P_{n+1}(\cos\theta_*)],
\label{eq:CDFSphereArrival}
\end{align}
\al{where we have applied the recurrence relation
\[
(2n+1) P_n(x) = \frac{d}{dx}[P_{n+1}(x) - P_{n+1}(x)].
\]}
This distribution is graphed in Fig.~\ref{fig:app:angle_dist}.
The first term corresponds to a  uniform distribution on the spherical surface; when $t_* \gg 1$, $\chi_n(t_*) \to 0$  (cf. \ref{chinlarget}) and this term dominates.

\section{Asymptotic analysis of capture time distribution}\label{sec:asy}

In this section we describe new asymptotic results for full time dependent arrival statistics and splitting probabilities for well separated absorbers on a plane.

\subsection{Asymptotic description of arrival times to planar regions with absorbing traps}\label{sec:asyplane}

As a confirmation of our results, we derive in this section an approximate distribution of arrival times based on an asymptotic analysis in the limit of vanishing pore radius. A novelty of this result is that it provides an explicit expressions for the full time dependent dynamics of the arrivals to individual absorbing sites, a significant extension on previous results relating to steady state capture rates \cite{sneddon1966mixed,Strieder08,Strieder12,berez2006,muratov}

The first step in this process to apply the Laplace transform $\hat{p}(\bx;s) = \int_0^{\infty} p(\bx,t) e^{-st}dt$ to \eqref{eqn:IntroP} which gives \cite{CLNQ22} the modified Helmholtz problem
\bsub\label{eqn:HelmholtzP}
\begin{gather}
\label{eqn:HelmholtzP_a}  D\Delta \hat{p} - s \, \hat{p} = -\delta(\bx-\bx_0), \qquad \bx \in \Om;\\[5pt] 
\label{eqn:HelmholtzP_b}  \hat{p} = 0 \quad \mbox{on} \quad \bx \in \Gama; \qquad 
\nabla  \hat{p} \cdot \hn = 0 \quad \mbox{on} \quad \bx \in \Gamr.
\end{gather}
\esub
We consider the absorbing set $\Gama$ to be $N$ non-overlapping absorbing pores in the plane $\pOm$ \al{where $\Omega$ is the upper half-space $z\geq0$}. For planar pores centered at coordinates $\bx_j= (x_j,y_j,0)$, we have
\begin{equation}\label{eq:domains}
\Gamma_a = \bigcup_{j=1}^N  \A_j \left (\frac{\bx -\bx_j} {\eps} \right ), \qquad \bx = (x,y,0).
\end{equation}
Here $\A_j$ is  a closed bounded set representing the shape of the $j^\text{th}$ pore and $\eps$ is a scale factor which yields a parameterized family of homothetic shapes centered at $\bx_j$. The asymptotic result for the distribution will be obtained in terms of the capacitances and locations of the individual pores on the plane. The following analysis adapts techniques \cite{BL2018} for the solution of Laplace's equation to the determination of capacitance to modified Helmholtz problem \eqref{eqn:HelmholtzP}.

A key tool in the analysis is knowledge of the modified Helmholtz Green's function $G_h$ satisfying 
\bsub\label{eqn:HelmholtzGreens}
\begin{gather}
\label{eqn:HelmholtzGreens_A} D \Delta G_h - s G_h = -\delta(\bx-\bxi), \qquad \bx \in \Om;\\[5pt] 
\label{eqn:HelmholtzGreens_B}  \int_{\Om}G_h\, d\bx = \frac{1}{s}; \qquad 
   \partial_z G_{h} = 0 \quad  \text{on} \quad z=0.
\end{gather}

The method of images provides the explicit Green's function
\begin{equation}\label{eqn:HelmholtzGreens_C}
G_h(\bx,\bxi) = \frac{1}{4\pi D} \left[ \frac{e^{-\alpha |\bx-\bxi|}}{ |\bx-\bxi|} + \frac{e^{-\alpha |\bx-\bxi^{\ast}|}}{ |\bx-\bxi^{\ast}| } \right], \qquad \alpha = \sqrt{\frac{s}{D}},
\end{equation}
\esub
where for $\bxi=(\xi_1,\xi_2,\xi_3)$, $\bxi^{\ast}=(\xi_1,\xi_2,-\xi_3)$ is the image point obtained from reflection in the plane $z=0$. For the case $\xi_3=0$, the local behavior near the singularity takes the form
\begin{equation}\label{eqn:HelmholtzGreens_c} 
G_h \sim \frac{1}{2\pi D}\frac{1}{|\bx-\bxi|} + R_h +\littleoh(1) \quad \textrm{as} \quad \bx\to\bxi; \qquad R_h = -\frac{\al{\alpha}}{2\pi D} ,
\end{equation}
where $R_h$ is referred to as the regular part of $G_h$. 

We now outline the key steps in obtaining an asymptotic solution of \eqref{eqn:HelmholtzP} as $\eps\to0$. Full details of the analysis can be found in the related works \cite{BL2018,LWB2017,LBS2018,Lindsay2023a}. We first examine a boundary layer solution of \eqref{eqn:HelmholtzP} near each pore through the stretched co-ordinates
\[ 
w_j(\by) = \hat{p}(\bx_j + \eps \by;s), \qquad \by  = (\sone,\stwo,\eta) \equiv  (\bx-\bx_j)/\eps,
\]
where $\eta$ is a scaled distance to the plane $z=0$ and $\vs = (\sone,\stwo)\in\mathbb{R}^2$ is a coordinate system on the plane centered at the $j^{th}$ pore. In this boundary layer, the solution for $w_j$ satisfies a static Laplace problem up to terms of $\bigoh(s\eps^2)$ which we assume to be negligible. \al{In addition, we assume that $|\bx_0-\bx_j| = \mathcal{O}(1)$ as $\eps\to0$ so that the effect of the initial condition does not arise in the local problem near the pore.} If $\wdisk_j$ solves the static single pore problem
\bsub\label{wc}
\begin{gather}
\Delta_{\by}  \wdisk_j = 0\,, \qquad \eta>0\,, \quad \vs \in\mathbb{R}^2;
\\[5pt]
\wdisk_j(\vs,0) = 1, \quad \vs \in \A_j\,;
\qquad \partial_{\eta} \wdisk_j (\vs,0)= 0, \quad 
\vs  \notin \A_j \, , \\[5pt]
\lim_{\rho \to \infty} \wdisk_j(\vs,\eta) =0 \, , 
\end{gather}
where the Laplacian, $\Delta_{\by}$, is in the stretched $\by=(\sone,\stwo,\eta)$ coordinates, 
then
\begin{equation}\label{eq:inner}
w_j = A_j( 1- \wdisk_j) + \eps B_j (1- \wdisk_j) + \bigoh(\eps^2),
\end{equation}
where $A_j$ and $B_j$ may depend on the Laplace transform parameter $s$. If we define $\rho \equiv \left | \by \right | = \sqrt{|\vs|^2+\eta^2} $ as the distance to the pore center,
the far-field behavior (cf.~\cite{belyaev}) is
\begin{equation}\label{wc_farfield}
\wdisk_j \sim  \frac{c_j}{\rho} + \bigoh \left(\frac{1}{\rho^2}\right) + \cdots \,,
\qquad \mbox{as} \quad \rho \to\infty \,.
\end{equation}
\esub
where $c_j$ is the \emph{capacitance} of the $j^\text{th}$ pore. Integrating equation \eqref{eq:inner} over a hemisphere of radius $R$ and passing to the limit as $R\to\infty$, we can obtain the transformed flux \ajb{density}
\begin{equation}
\hat{\J}_j = D\int_{\mathcal{A}_j} \partial_{\eta} w_j |_{\eta=0}\, d \vs = 2\pi D [c_jA_j + \eps c_j B_j + \bigoh(\eps^2)].
\end{equation}

For the case where $\A_j$ is the circular pore $|\vs| \leq a_j$, \eqref{wc} corresponds to the well-known electrified disk problem
(cf.~\cite{sneddonelements}) with solution
\begin{equation}\label{wc:sol}
\wdisk_j(\vs,\eta)  = \frac{2}{\pi}\sin^{-1}\left[
\frac{ 2a_j} 
{\sqrt{  \left (  |\vs | + a_j \right )^{2} + \eta^2} + \sqrt{ \left ( |\vs| - a_j \right )^{2} + \eta^2}}
\right], \qquad c_j = \frac{2a_j}{\pi}.
\end{equation}
The far-field of the boundary layer solution results provides a matching condition on the solution of \eqref{eqn:HelmholtzP} of the form
\begin{equation}
\hat{p} \sim A_j + \eps \left[ - \frac{c_j A_j }{|\bx - \bx_j|} + B_j  \right] - \eps^2  \frac{c_j B_j }{|\bx - \bx_j|} + \cdots, \qquad \bx\to\bx_j,
\end{equation}
where the $A_j$ and $B_j$ are constants to be determined. The expansion of \eqref{eqn:HelmholtzP} as $\eps\to0$ is therefore of the form $\hat{p}= \hat{p}_0 + \eps \hat{p}_1+ \eps^2 \hat{p}_2+\bigoh(\eps^3)$ where at leading order
\bsub\label{eqn:HelmholtzP0}
\begin{gather}
\label{eqn:HelmholtzP0_a}  D\Delta \hat{p}_0 - s\,  \hat{p}_0 = -\delta(\bx-\bx_0), \qquad \bx \in \Om\setminus\{ \bx_0\};\\[5pt] 
\label{eqn:HelmholtzP0_b}  \hat{p}_0 \sim A_j + \cdots \quad \mbox{as} \quad \bx \to \bx_j; \qquad 
  \hat{p}_{0z} = 0 \quad \text{on} \quad z=0.
\end{gather}
\esub
The solution of \eqref{eqn:HelmholtzP0} is expressed as $\hat{p}_0 = G_h(\bx,\bx_0;\alpha)$ where $G_h$ satisfies
 \eqref{eqn:HelmholtzGreens}. Matching $\hat{p}_0= G_h(\bx,\bx_0;\alpha)$ with the local condition \eqref{eqn:HelmholtzP0_b} gives
\begin{equation}\label{eqn:Aj}
A_j = G_h(\bx_j,\bx_0;\alpha),
\end{equation}
which quantifies the influence of the initial condition on the absorption rate at $\bx_j$. At the following order, the equation for $\hat{p}_1$ is 
\bsub\label{eqn:HelmholtzP1}
\begin{gather}
\label{eqn:HelmholtzP1_a}  D\Delta \hat{p}_1 - s\,  \hat{p}_1 = 0, \qquad \bx \in \Om \setminus \{ \bx_1,\ldots,\bx_N\};\\[5pt] 
\label{eqn:HelmholtzP1_b}  \hat{p}_1 \sim \frac{-c_jA_j}{|\bx - \bx_j |} + B_j + \cdots \quad \mbox{as} \quad \bx \to \bx_j; \qquad 
  \hat{p}_{1z} = 0 \quad \mbox{on} \quad  z=0 ,
\end{gather}
which can be expressed as
\begin{equation}\label{eqn:P1}
\hat{p}_1 = -2\pi D\sum_{j=1}^N c_j A_j G_h(\bx,\bx_j;\alpha).
\end{equation}
\esub
Applying the local condition \eqref{eqn:HelmholtzP1_b} in \eqref{eqn:P1} yields that
\begin{equation}
\label{eqn:Bj}
B_j = - 2\pi D \Big[c_j A_j R_h + \sum_{\substack{k =1\\ k\neq j} }^N c_k A_k G_h(\bx_j, \bx_k;\alpha) \Big].
\end{equation}
The terms $B_j$ quantify the effect of inter-pore competition on the rate of capture at the $j^\text{th}$ pore. The equation for $\hat{p}_2$ is now
\bsub\label{eqn:HelmholtzP2}
\begin{gather}
\label{eqn:HelmholtzP2_a} D \Delta \hat{p}_2 - s\,  \hat{p}_2 = 0, \qquad \bx \in \Om \setminus \{ \bx_1,\ldots,\bx_N\};\\[5pt] 
\label{eqn:HelmholtzP2_b}  \hat{p}_2 \sim \frac{-c_jB_j}{|\bx - \bx_j |} + \cdots \quad \mbox{as} \quad \bx \to \bx_j; \qquad 
  \hat{p}_{2z} = 0 \quad \mbox{on} \quad   z=0 ,
\end{gather}
which has solution
\begin{equation}\label{eqn:P2_sol}
\hat{p}_2 = -2\pi D\sum_{j=1}^N c_jB_j G_h(\bx,\bx_j;\alpha).
\end{equation}
\esub
At this point we end the asymptotic procedure, though this recursive matching process can be continued to obtain higher order corrections if desired. Such corrections would describe the influences of pore shape and heterogeneity (differences in capacitances) on the capture rate. Details and further examples of the higher order matching process can be found in \cite{LTK2015,belyaev,LBS2018}.

The survival probability $S(t)= \int_{\Omega} p(\bx,t) \, d\bx$ then has Laplace transform $\hat{S}$ given by
\begin{equation}\label{eqn:P2}
\hat{S} = \frac{1}{s} - \frac{2\pi\eps D}{s} \sum_{j=1}^N \Big(c_j A_j + \eps c_j B_j \Big) +\bigoh(\eps^3),
\end{equation}
where we have used the identity $\int_{\Omega}G_h(\bx;\bxi,\alpha) d\bx = s^{-1}$ from \eqref{eqn:HelmholtzGreens_B}. It is easier to work with \al{the} PDF of the capture time distribution, $\J(t)$, and its Laplace transform, $\hat{\J}(s)$. \al{We calculate that 
\begin{equation}
\frac{d}{dt} S(t) 
= \int_{\Omega} p_t(\bx,t) \, d\bx 
=\int_{\Omega} D \Delta p \, d\bx 
= - D\int_{\partial\Omega} \nabla p(\bx,t) \cdot \hat{\bf n}\, dS = -\J(t),
\end{equation}
where the negative sign is due to the normal at the boundary pointing into the bulk.\footnote{\ajb{Strictly speaking, the divergence theorem needs to be applied on a closed surface. As $p(\bx,t)$ is exponentially small far from the source one can justify this calculation by considering a hemisphere whose base contains the absorbing pores and letting its radius tend to infinity.}}
In Laplace space this implies that 
$$
\hat{\J}(s) = S(0)  -s \hat{S}(s) = 1 -s \hat{S}(s).
$$}
Equation \eqref{eqn:P2} now yields
\begin{equation}\label{eqn:Chat}
\hat{\J}(s)=  {2\pi\eps D}\sum_{j=1}^N \Big(c_j A_j + \eps c_j B_j \Big) +\bigoh(\eps^3),
\end{equation}

Substitution \eqref{eqn:Aj}  and \eqref{eqn:Bj} into \eqref{eqn:Chat} yields
\begin{align*}
\hat{\J}(s)  &=  2\pi  D  \,\sum_{j=1}^N \eps c_j G_h(\bx_j,\bx_0;\alpha) \left [
1-   2 \pi \eps D  c_j R_h \right ]  \\
&-  (2\pi  D)^2  \,\sum_{j=1}^N   \sum_{\substack{k =1\\ k\neq j} }^N
 \eps^2 c_j c_k G_h(\bx_k,\bx_0;\alpha) G_h(\bx_j,\bx_k;\alpha)
 +\bigoh(\eps^3)
, \\
 &= \sum_{j=1}^N \eps c_j \frac{e^{-\alpha |\bx_j-\bx_0|}}{|\bx_j-\bx_0|} \left [1+ \eps c_j \alpha \right ]
-  \eps^2    \sum_{j=1}^N  \sum_{\substack{k =1\\ k\neq j} }^N
c_j c_k \frac{e^{-\alpha \left(|\bx_k-\bx_0| +|\bx_j-\bx_k|\right ) }}{|\bx_k-\bx_0| |\bx_j-\bx_k|} +\bigoh(\eps^3), \\
&=\sum_{j=1}^N \eps c_j \frac{e^{-\alpha R_j}}{R_j} \left [1+ \eps  c_j \alpha \right ]
-  \eps^2    \sum_{j=1}^N  \sum_{\substack{k =1\\ k\neq j} }^N
c_j c_k \frac{e^{-\alpha \left(R_k +d_{jk}\right ) }}{R_k d_{jk}} +\bigoh(\eps^3),
 \end{align*}
where for the $N$ circular traps at $\bx_j=(x_j,y_j,0)$ of capacitance $\eps c_j$  and a source at $\bx_0 = (x_0,y_0,z_0)$
we have defined $R_j = |\bx_j-\bx_0|$ and $d_{jk} = |\bx_j-\bx_k| $. The inverse Laplace transform of this expression \al{\cite[{\S 29.3.82}, {\S 29.3.87}]{abramowitz1965handbook}} yields
\begin{align}\label{eq:asyCap}
\J(t)   &=  \eps \frac{D}{2 \sqrt{\pi}} 
\sum_{j=1}^N c_j
e^{-\frac{R_j^2}{4Dt} }  \left [  
\left  (1- \frac { \eps c_j}{ R_j} \right )  \frac 1 {(D t)^{3/2} } 
+   
\frac{\eps c_j R_j} {2}    \frac 1 {(D t)^{5/2} } \right ] \\
&-  \eps^2  \frac{D}{2 \sqrt{\pi}}   \sum_{j=1}^N  \sum_{\substack{k =1\\ k\neq j} }^N c_j c_ke^{-\frac{(R_k+d_{jk})^2}{4Dt} }
\left [ \left  (\frac 1 { R_k}+ \frac 1 {d_{jk}} \right ) 
 \frac 1 {(D t)^{3/2} } \right ] +\bigoh(\eps^3). \nonumber
\end{align}
The CDF of this distribution, defined as $F(t) = \int_0^t \J(\tau) \, d\tau$,  can then be calculated as
\begin{equation}\label{asy:CDF}
F(t) =  \eps \sum_{j=1}^N
\frac{c_j}{ R_j} \mbox{erfc} \left (\frac{R_j}{2 \sqrt{Dt}} \right )
+  \eps^2  \sum_{j=1}^N \Big[ \frac{c_j^2}{R_j}
\frac{e^{-\frac{R_j^2}{4Dt} } }{\sqrt{\pi Dt}}  - \sum_{\substack{k =1\\ k\neq j} }^N
 \frac{c_j c_k}{d_{jk} R_k } \mbox{erfc} \left (\frac{R_k+d_{jk}}{2 \sqrt{Dt}} \right ) \Big]+\bigoh(\eps^3).
\end{equation}
We remark that the capture probability is given by
\begin{equation}\label{eq:final_fullC}
\int_0^{\infty} \J(\tilt) \, d \tilt  = \lim_{t\to\infty}F(t) = \eps\sum_{j=1}^N \frac{c_j} {R_j}  - \eps^2 \sum_{j=1}^N  \sum_{\substack{k =1\\ k\neq j} }^N \frac{ c_j c_k} {d_{jk} R_k}  +\bigoh(\eps^3),
\end{equation}
which is not unity. This reflects the fact that a fraction of particles, inversely proportional to the initial distance to the cluster of traps, escape to infinity.

For a single trap centered at the origin, $\bx_1=(0,0,0)$, of capacitance $\eps c $ and a source at
$\bx_0 = (x_0,y_0,z_0)$ with $|\bx_0|=R$, this reduces to
\begin{equation}\label{eq:asyCaptureOne}
\J(t)   =  \eps c \frac{D}{2 \sqrt{\pi}} 
e^{-\frac{R^2}{4Dt} }  \left [  
\left  (1- \frac { \eps c}{ R} \right )  \frac 1 {(D t)^{3/2} } 
+   
\frac{\eps  c R} {2}    \frac 1 {(D t)^{5/2} } \right ] +\bigoh(\eps^3).
\end{equation}
For a circular trap of radius $\eps a$, the capacitance is given by $c=2a/\pi$.
The CDF of this distribution is
\begin{equation}
\label{eq:asyIntCaptureOne}
F(t) = \frac{\eps c}{ R} \mbox{erfc} \left (\frac{R}{2 \sqrt{Dt}} \right )
+ \frac{(\eps c)^2}{R}
\frac{e^{-\frac{R^2}{4Dt} } }{\sqrt{\pi Dt}} +\bigoh(\eps^3).
 \end{equation}
Also, $\int_0^{\infty} \J(\tau) \, d \tau  = {\eps c}/{ R}+\bigoh(\eps^3) $ so that the probability of capture is not unity, but inversely proportional to the initial distance to the trap center. 

\subsection{Asymptotic determination of splitting probabilities}\label{sec:asysplit}

Here we give a brief calculation of the planar splitting probabilities $\{ \Q_k(\bx) \}_{k=1}^N$ which solve
\begin{gather}
\Delta \Q_k = 0 \qquad \bx \in \Om; \qquad \Q_k(\bx) \quad \mbox{bounded as } \ |\bx|\to\infty;\\[5pt] 
\Q_k = \delta_{jk} \quad \mbox{on} \quad \bx \in \Omega_j,\quad  j = 1,\ldots,N; \qquad 
\nabla \Q_k \cdot \hn = 0 \quad \mbox{on} \quad \bx \in \Gamr ,
\end{gather}
where the definition of the planar regions is given in \eqref{eq:domains}. The relevant surface Green's function $G_s$ for the Laplacian is given by
\begin{equation}\label{eqn:SurfaceGreens}
\Delta G_s = 0, \quad \bx \in \Om; \qquad  \nabla G_s \cdot \hn = -\delta(\bx-\bxi) \quad \mbox{on} \quad \bx \in \partial\Omega\setminus\{\bxi\} , \qquad G_s(\bx,\bxi) = \frac{1}{2\pi|\bx-\bxi|}.
\end{equation}
The analysis is very similar to the previous section and so we provide an abbreviated account. The inner solution near the $j^{th}$ pore (equivalent of \eqref{eq:inner}) takes the form 
\[
\Q_k(\bx_j + \eps \by) \sim \delta_{jk} \bar{w}_j(\by) + \eps A_{jk} (1 - \bar{w}_j(\by))+ \cdots, \qquad j = 1,\ldots,N,
\]
and constants $A_{jk}$. In light of the large argument behavior $ \bar{w}_j(\by) \sim c_j/|\by|$ as $|\by|\to\infty$ for $\by = (\bx-\bx_j)/\eps$, we have that
\begin{equation}
\Q_k(\bx) \sim \eps\Big[ \frac{c_k}{|\bx-\bx_k|} + A_{jk} \Big]- \eps^2 \Big[ \sum_{j=1}^N \frac{A_{jk}c_j}{|\bx-\bx_j|} + \cdots\Big], \qquad \bx\to\bx_k.
\end{equation}
This motivates the expansion $\Q_k(\bx) = \eps\Q_k^{(1)}(\bx)+ \eps^2 \Q_k^{(2)}(\bx)+ \cdots $ where we form equations for $\Q_k^{(1)}(\bx)$ and $\Q_k^{(2)}(\bx)$. Solving these equations yield
\bsub
\begin{gather}
\Q_k^{(1)}(\bx) = 2\pi c_k G_s(\bx;\bx_k), \qquad A_{jk} =  2\pi (1-\delta_{jk}) c_j G_s(\bx_j;\bx_k);\\[5pt]
\Q_k^{(2)}(\bx) = -2\pi \sum_{j=1}^N  c_j A_{jk} G_s(\bx;\bx_j) = -4\pi^2 \sum_{\substack{j =1\\ j\neq k} }^N c_j c_k G_s(\bx_j;\bx_k) G_s(\bx;\bx_j).
\end{gather}
\esub
The series $\Q_k(\bx)  = \eps \Q_k^{(1)}(\bx) + \eps^2\Q_k^{(2)}(\bx) + \cdots$ now yields the final result for the splitting probability
\begin{equation}\label{eq:split_final}
\Q_k(\bx) = \frac{\eps c_k}{|\bx-\bx_k|} - \eps^2 \sum_{\substack{j =1\\ j\neq k} }^N \frac{c_j c_k}{|\bx_j - \bx_k | |\bx- \bx_j |}  + \mathcal{O}(\eps^3).
\end{equation}
An addition verification of the result \eqref{eq:split_final} can be obtained by summing over $1\leq k\leq N$ to arrive at 
\[
\sum_{k=1}^N \Q_k(\bx) = \eps \sum_{k=1}^N \frac{ c_k}{R_k} - \eps^2  \sum_{k=1}^N\sum_{\substack{j =1\\ j\neq k} }^N \frac{c_j c_k}{d_{jk} R_j } + \mathcal{O}(\eps^3), \qquad \begin{array}{rcl} R_j&=& | \bx - \bx_j| \\[2pt] d_{jk} &=& |\bx_j - \bx_k | \end{array}
\]
which agrees with \eqref{eq:final_fullC}. 

\bibliographystyle{siam}
\bibliography{newbib}

\end{document}